\documentclass[10pt,a4paper]{article}

\usepackage[utf8]{inputenc}
\usepackage{amsfonts}
\usepackage{amsmath} 	
\usepackage{amssymb} 	
\usepackage{appendix}
\usepackage{bbm}        
\usepackage{booktabs}   
\usepackage{cancel}     
\usepackage{caption}    
\usepackage{cite}       
\usepackage{colortbl}   
\usepackage{eso-pic} 	
\usepackage{float}      
\usepackage{fullpage} 	
\usepackage{hyperref}	
\usepackage{graphicx} 	
\usepackage{moreverb}	
\usepackage{pgf, tikz}	
\usepackage{pgfplots}	
    \usetikzlibrary{arrows,shapes,positioning}
    \usetikzlibrary{patterns}
    \usepgfplotslibrary{colormaps}   
    \pgfplotsset{compat=1.15}
\usepackage[figuresleft]{rotating} 
\usepackage{stmaryrd}
\usepackage{subcaption} 

\DeclareMathOperator{\spn}{span}

\definecolor{mesh_patch_center}{rgb}{0.73725,0.74118,0.13333} 
\definecolor{mesh_patch_x1}{rgb}{0.83922,0.15294,0.15686} 
\definecolor{mesh_patch_x2}{rgb}{0.17255,0.62745,0.17255} 
\definecolor{mesh_patch_x3}{rgb}{0.54902,0.33725,0.29412} 
\definecolor{mesh_patch_x4}{rgb}{0.58039,0.40392,0.74118} 
\definecolor{mesh_patch_x5}{rgb}{0.12157,0.46667,0.70588} 
\definecolor{mesh_patch_x6}{rgb}{1.000000,0.49804,0.0549} 
\definecolor{mesh_patch_x7}{rgb}{0.49804,0.49804,0.49804} 
\definecolor{mesh_patch_x8}{rgb}{0.89020,0.46667,0.76078} 

\makeatletter
\let\@fnsymbol\@arabic
\makeatother

\begin{document}

\title{Local cubic spline interpolation for Vlasov-type equations on a multi-patch geometry.
}
\author{Pauline Vidal 
            \thanks{Corresponding author: Pauline Vidal \\
                Max-Planck-Institut für Plasmaphysik, Garching, Germany
                Email: pauline.vidal@ipp.mpg.de  
                ORCID: \url{https://orcid.org/0009-0008-7233-6351}},
        \and Emily Bourne
            \thanks{Emily Bourne \\
                SCITAS, EPFL, CH-1015 Lausanne, Switzerland
                ORCID: \url{https://orcid.org/0000-0002-3469-2338}},
        \and Virginie Grandgirard 
            \thanks{Virginie Grandgirard \\
                CEA, IRFM, 13108 Saint-Paul-lez-Durance Cedex, France
                ORCID: \url{https://orcid.org/0000-0001-7821-9107}},
        \and Michel Mehrenberger 
            \thanks{Michel Mehrenberger \\
                Aix Marseille Univ, CNRS, I2M, Marseille, France
                ORCID: \url{https://orcid.org/0000-0001-6376-409X}},
        \and Eric Sonnendrücker 
            \thanks{Eric Sonnendrücker \\
                Max-Planck-Institut für Plasmaphysik, Garching, Germany,
                Technische Universität München, Zentrum Mathematik, Garching, Germany 
                ORCID: \url{https://orcid.org/0000-0002-8340-7230}}.
}
\date{Last update: January, 8th, 2026}

\maketitle

\begin{abstract}
We present a semi-Lagrangian method for the numerical resolution of Vlasov-type equations on multi-patch meshes.  
Following N. Crouseilles et al. [\textit{A parallel Vlasov solver based on local cubic spline interpolation on
patches.} Journal of Computational Physics (2009)], 
we employ a local cubic spline interpolation with Hermite boundary conditions between the patches.
The derivative reconstruction is adapted to cope with non-uniform meshes as well as non-conforming situations.
In the conforming case, the constraint of the number of points for each patch, found in previous studies, 
is removed; however, a small global system must now be solved.
In that case, the local spline representations coincide with the corresponding global spline reconstruction.
Alternatively, we can choose not to apply the global system and the
derivatives can be approximated.
The influence of the most distant points diminishes as the number of points per patch increases. 
For uniform per patch configurations, a study of the explicit and asymptotic behavior of this influence has been led. 
The method is validated using a two-dimensional guiding-center model with an O-point.
All the numerical results are carried out in the Gyselalib++ library.


\end{abstract}

\paragraph{Acknowledgements}
This work was carried out using the Gyselalib++ library. 
The project has received funding from the EoCoE-III project
(Energy-oriented Centre of Excellence for Exascale HPC applications).
This work has been carried out within the framework of the EUROfusion Consortium, 
funded by the European Union via the Euratom Research and Training Program 
(Grant Agreement No. 101052200 — EUROfusion). Views and opinions expressed are
those of the author(s) only and do not necessarily reflect those of the European Union 
or the European Commission. Neither the European Union nor the European Commission can be
held responsible for them.
Emily Bourne's salary was paid for by the EUROfusion Advanced Computing Hub (ACH).

\section{Introduction}
\label{section_introduction}
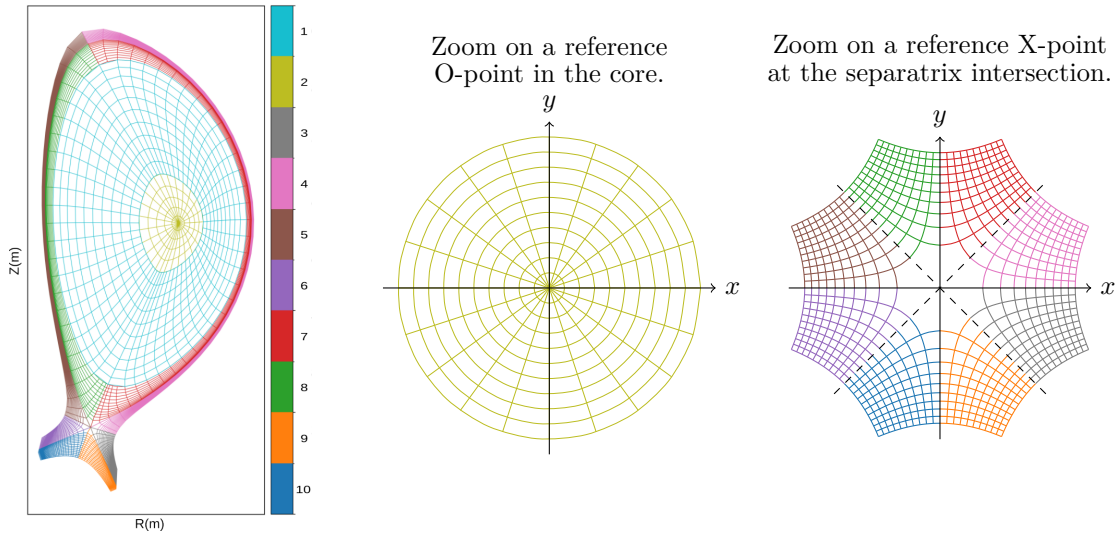
\begin{figure}[h]
    \centering
    \includegraphics[width= 4cm]{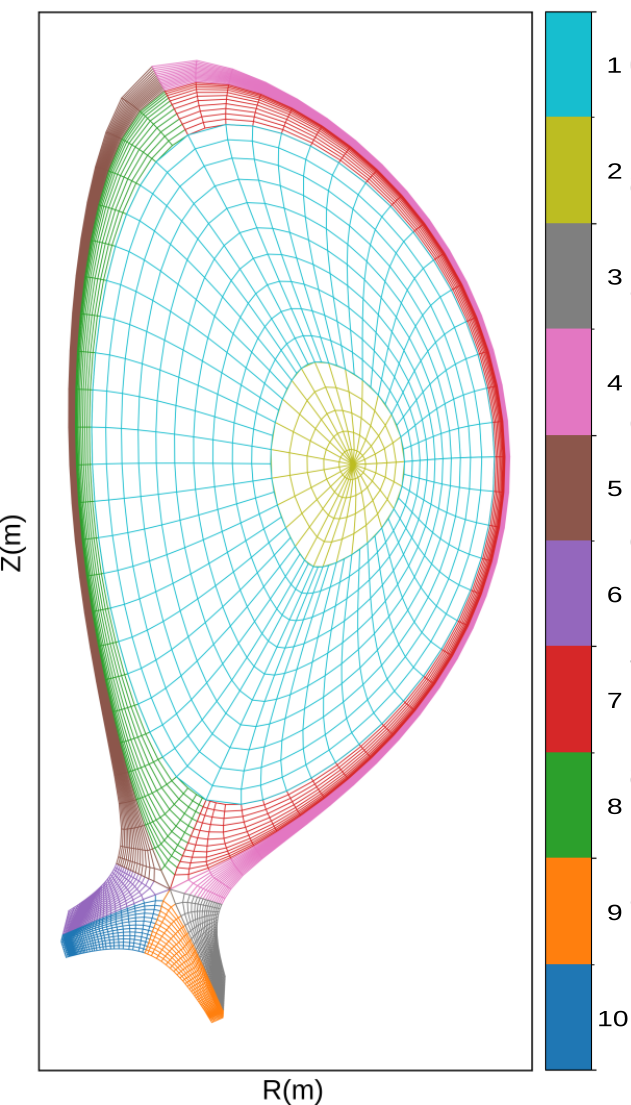}
    \hspace*{0.5cm}
    \begin{tikzpicture}[xscale = 0.4, yscale = 0.4]
        \node at (0, 8) {Zoom on a reference} ; 
        \node[below] at (0, 7.75) {O-point in the core.} ; 
        \draw[white] (-6,-8)--(0,0); 
        \def\R{5}; 
        \def\Theta{6.28318530717958}; 
        \foreach \r in {0, 0.5, ..., 5}{
            \draw[mesh_patch_center, line width = 0.3, domain=0:\Theta, smooth, variable=\th, samples=20]
                plot ({\r*cos(deg(\th))}, {\r*sin(deg(\th))});
        }
        \foreach \th in {0, 0.314159, ..., \Theta}{
            \draw[mesh_patch_center, line width = 0.3, domain=0:\R, smooth, variable=\r, samples=10]
                plot ({\r*cos(deg(\th))}, {\r*sin(deg(\th))});
        }
        \draw[->] ({-\R*1.1},0) -- ({\R*1.1},0) node[right] {$x$};
        \draw[->] (0,{-\R*1.1}) -- (0,{\R*1.1}) node[above] {$y$};
    \end{tikzpicture}
    \begin{tikzpicture}[xscale = 0.4, yscale = 0.4]
        \node at (0, 8) {Zoom on a reference X-point} ; 
        \node[below] at (0, 7.75) {at the separatrix intersection.} ; 
        \draw[white] (-6,-8)--(0,0); 
        \def\sgnx{1}
        \def\sgny{1}
        \foreach \u in {1, 2,...,10} { 
            \draw[mesh_patch_x1, domain=0:10, smooth, variable=\v, samples=20]
            plot ({\sgnx*sqrt(sqrt(\u*\u + \v*\v) - \u)}, {\sgny*sqrt(sqrt(\u*\u + \v*\v)+\u)});
        }
        \foreach \v in {1, 2,...,10} { 
            \draw[mesh_patch_x1, domain=0:10, smooth, variable=\u, samples=20]
            plot({\sgnx*sqrt(sqrt(\u*\u + \v*\v) - \u)}, {\sgny*sqrt(sqrt(\u*\u + \v*\v)+\u)});
        }

        \def\sgnx{-1}
        \def\sgny{1}
        \foreach \u in {1, 2,...,10} { 
            \draw[mesh_patch_x2, domain=0:10, smooth, variable=\v, samples=20]
            plot ({\sgnx*sqrt(sqrt(\u*\u + \v*\v) - \u)}, {\sgny*sqrt(sqrt(\u*\u + \v*\v)+\u)});
        }
        \foreach \v in {1, 2,...,10} { 
            \draw[mesh_patch_x2, domain=0:10, smooth, variable=\u, samples=20]
            plot({\sgnx*sqrt(sqrt(\u*\u + \v*\v) - \u)}, {\sgny*sqrt(sqrt(\u*\u + \v*\v)+\u)});
        }

        \def\sgnx{-1}
        \def\sgny{1}
        \foreach \u in {-10, -9, ..., -1} { 
            \draw[mesh_patch_x3, domain=0:10, smooth, variable=\v, samples=20]
            plot ({\sgnx*sqrt(sqrt(\u*\u + \v*\v) - \u)}, {\sgny*sqrt(sqrt(\u*\u + \v*\v)+\u)});
        }
        \foreach \v in {1, 2,...,10} { 
            \draw[mesh_patch_x3, domain=-10:0, smooth, variable=\u, samples=20]
            plot({\sgnx*sqrt(sqrt(\u*\u + \v*\v) - \u)}, {\sgny*sqrt(sqrt(\u*\u + \v*\v)+\u)});
        }
        
        \def\sgnx{-1}
        \def\sgny{-1}
        \foreach \u in {-10, -9, ..., -1} { 
            \draw[mesh_patch_x4, domain=0:10, smooth, variable=\v, samples=20]
            plot ({\sgnx*sqrt(sqrt(\u*\u + \v*\v) - \u)}, {\sgny*sqrt(sqrt(\u*\u + \v*\v)+\u)});
        }
        \foreach \v in {1, 2,...,10} { 
            \draw[mesh_patch_x4, domain=-10:0, smooth, variable=\u, samples=20]
            plot({\sgnx*sqrt(sqrt(\u*\u + \v*\v) - \u)}, {\sgny*sqrt(sqrt(\u*\u + \v*\v)+\u)});
        }

        \def\sgnx{-1}
        \def\sgny{-1}
        \foreach \u in {1, 2,...,10} { 
            \draw[mesh_patch_x5, domain=0:10, smooth, variable=\v, samples=20]
            plot ({\sgnx*sqrt(sqrt(\u*\u + \v*\v) - \u)}, {\sgny*sqrt(sqrt(\u*\u + \v*\v)+\u)});
        }
        \foreach \v in {1, 2,...,10} { 
            \draw[mesh_patch_x5, domain=0:10, smooth, variable=\u, samples=20]
            plot({\sgnx*sqrt(sqrt(\u*\u + \v*\v) - \u)}, {\sgny*sqrt(sqrt(\u*\u + \v*\v)+\u)});
        }

        \def\sgnx{1}
        \def\sgny{-1}
        \foreach \u in {1, 2,...,10} { 
            \draw[mesh_patch_x6, domain=0:10, smooth, variable=\v, samples=20]
            plot ({\sgnx*sqrt(sqrt(\u*\u + \v*\v) - \u)}, {\sgny*sqrt(sqrt(\u*\u + \v*\v)+\u)});
        }
        \foreach \v in {1, 2,...,10} { 
            \draw[mesh_patch_x6, domain=0:10, smooth, variable=\u, samples=20]
            plot({\sgnx*sqrt(sqrt(\u*\u + \v*\v) - \u)}, {\sgny*sqrt(sqrt(\u*\u + \v*\v)+\u)});
        }

        \def\sgnx{1}
        \def\sgny{-1}
        \foreach \u in {-10, -9, ..., -1} { 
            \draw[mesh_patch_x7, domain=0:10, smooth, variable=\v, samples=20]
            plot ({\sgnx*sqrt(sqrt(\u*\u + \v*\v) - \u)}, {\sgny*sqrt(sqrt(\u*\u + \v*\v)+\u)});
        }
        \foreach \v in {1, 2,...,10} { 
            \draw[mesh_patch_x7, domain=-10:0, smooth, variable=\u, samples=20]
            plot({\sgnx*sqrt(sqrt(\u*\u + \v*\v) - \u)}, {\sgny*sqrt(sqrt(\u*\u + \v*\v)+\u)});
        }

        \def\sgnx{1}
        \def\sgny{1}
        \foreach \u in {-10, -9, ..., -1} { 
            \draw[mesh_patch_x8, domain=0:10, smooth, variable=\v, samples=20]
            plot ({\sgnx*sqrt(sqrt(\u*\u + \v*\v) - \u)}, {\sgny*sqrt(sqrt(\u*\u + \v*\v)+\u)});
        }
        \foreach \v in {1, 2,...,10} { 
            \draw[mesh_patch_x8, domain=-10:0, smooth, variable=\u, samples=20]
            plot({\sgnx*sqrt(sqrt(\u*\u + \v*\v) - \u)}, {\sgny*sqrt(sqrt(\u*\u + \v*\v)+\u)});
        }

        \draw[->] (-5,0) -- (5,0) node[right] {$x$};
        \draw[->] (0,-5) -- (0,5) node[above] {$y$};

        \draw[dashed] (-3.5,-3.5) -- (3.5,3.5);
        \draw[dashed] (3.5,-3.5) -- (-3.5,3.5);  
    \end{tikzpicture}
    \caption{\label{JET_SOLEDGE} 
    On the left, an example of global domain split into 10 patches. 
    Mesh for the ITER tokamak based on data from SOLEDGE3X and generated by K. Obrejan
    (\url{https://soledge3x.onrender.com}, article \cite{SOLEDGE3X}).
    At the center, a zoom on the core with
    a reference mapping with O-point.
    On the right, a zoom on the separatrix ($--$) intersection forming
    the X-point with a reference mapping.
    }
\end{figure}

\paragraph{}
In numerical tokamak simulations with semi-Lagrangian schemes,
we align the mesh lines on the flux surfaces of the magnetic field,
since the variations of the different physical quantities along 
the magnetic field lines are small \cite{Latu2018}. 
These mesh lines separate the poloidal plane into different zones (see Fig. \ref{JET_SOLEDGE}). 
In the core, the mesh lines converge to a central point called the O-point.
On the edge, the separatrix (the boundary between the closed and open magnetic field lines) forms 
an X-point at its intersection. The mesh lines diverge from the X-point around it. 
These zones raise different numerical difficulties that are interesting to address using meshes adapted 
to these specificities. 
The multi-patch approach is an effective tool for this geometry type.
In this approach, the multi-patch geometry corresponds to a global domain split 
into different subdomains called patches.
 
The multi-patch approach has interesting optimization properties. For instance, 
computations on uniform meshes can be more optimized than on non-uniform meshes
(see \cite{Qian2013} for an example of optimization on uniform meshes).
Non-uniform meshes can reduce memory constraints, but they can also potentially increase the 
calculation requirements \cite{Bourne2023Splines}.
So, if we want to locally refine an area of a global domain,
we can split it into multiple patches to be uniform piecewise and remain optimized.
Moreover, with a 2D basis formed from the tensor product of 1D function bases, the multi-patch 
approach allows us to deal with non-conforming meshes (as is the case in Fig. \ref{JET_SOLEDGE}). 
Studies have already been conducted on non-conforming refinement on a multi-patch
geometry for finite element methods (FEM) with B-splines
\cite{Gu2018, Güçlü2022, Campos2024}.
Splitting a domain into independent subdomains is also interesting for parallelization. 
We can cite the work in \cite{Crouseilles2009, Kormann2019} on MPI parallelization 
for Vlasov-Poisson equations with a domain decomposition and B-splines.
However, parallelizing operators on a poloidal plane with both an O-point and an X-point may be
difficult here to achieve a good load balance. This will therefore not be studied in this article. 

Once the domain is split into multiple patches, we are interested in simulating particle motion. 
In a tokamak, the time evolution of the particle distribution functions is governed by the Vlasov equation,
a 6D kinetic equation (3D in space and 3D in velocity).
It can be averaged over the cyclotron motion. This gives the gyrokinetic equation which 
describes a charge distribution function in a 5D phase space with 3 dimensions in space and 2 dimensions 
in velocity \cite{Latu2018, Grandgirard2016, Grandgirard2006},
\begin{equation}
    \partial_t f(t,\mathbf{x}, v_{\|}, \mu) 
    + \mathbf{v}_{GC}(t,\mathbf{x}, v_{\|}, \mu)  \cdot\nabla f (t,\mathbf{x}, v_{\|}, \mu) 
    + a_{GC}(t,\mathbf{x}, v_{\|}, \mu)  \partial_{v_{\|}} f (t,\mathbf{x}, v_{\|}, \mu)
    = 0,
    \label{gyrokinetic_eq}
\end{equation}
with $f$ the gyro-center distribution function of plasma particles, 
$\mathbf{v}_{GC}$ the gyrokinetic velocity,
$a_{GC}$ the gyrokinetic acceleration defined on position $\mathbf{x}\in\mathbb{R}^3$, 
$ v_{\|}\in\mathbb{R}$ the parallel velocity, and $ \mu\in\mathbb{R}^+$ the modified magnetic moment. 
To solve the Vlasov equation with the gyrokinetic approximation, we use 
an operator splitting method to combine different advection equations along each dimension.  
The main difficulty arises at special points (the O-point and the X-point) which are located 
in the poloidal cross-section. 
So, to study these poloidal cross-section singularities,
we focus in this paper on the following 2D model with guiding-center equations
where the advection is in the poloidal cross-section, 
\begin{equation}
    \left\{
        \begin{aligned}
            & \partial_t \rho + \mathbf{v}_{GC}\cdot\nabla\rho = 0, \\
            & \Delta \phi = -\rho,
        \end{aligned}
    \right.
    \label{guiding_center_eq}
\end{equation}
with $\rho$ the density of the plasma particles.
In the Poisson equation in \eqref{guiding_center_eq}, 
$\phi$ is the electric potential of the electric field $\mathbf{E} = - \nabla \phi$.
Additionally, the magnetic field is supposed to be a unit vector, $\mathbf{B} = \mathbf{e_z}$, 
so the drift velocity 
$ \mathbf{v}_{GC} =\frac{\mathbf{E}\times \mathbf{B}}{|B^2|} = \mathbf{E} \times \mathbf{e_z}$.

The numerical method and the results that will be presented are implemented in the Gyselalib++ library 
\cite{Gyselalibxx} which is a C++ rewrite of the GYSELA code \cite{Grandgirard2016}.

In the library, the Poisson equation on the poloidal cross-section in \eqref{guiding_center_eq} 
is solved with a FEM \cite{Zoni2019}.
The advection equation is solved with a backward semi-Lagrangian method (BSL) \cite{Sonnendrücker1999}. 
The basis functions used in the FEM and the interpolation are cubic B-splines.
They have good performance and interpolation qualities \cite{Unser1991}. 
So, we will study the advection on multiple patches solved with a BSL on cubic B-splines 
to implement the method in Gyselalib++.

In 2019, E. Zoni and Y. Güçlü,
introduced methods to solve guiding-center equations on O-point geometries \cite{Zoni2019}.
They present methods to solve Poisson equations with an O-point by modifying the basis 
functions around the O-point. The advection equation is solved in the physical domain by linearizing 
the advection field around the O-point. 
In 2023, E. Bourne et al.
studied different Poisson solvers \cite{Bourne2023Poisson} and implemented the Spline FEM solver in Gyselalib++.
We use the work of these two articles as a starting point and adapt it to multi-patch geometries with O-point. 

The main difficulty in the BSL method on a multi-patch domain lies in the interpolation 
and dealing with the patch connections.  
We manage this issue here by using Hermite boundary conditions on the local splines defined
on each patch. 
These conditions confer a $\mathcal{C}^1$ regularity to the global piecewise spline representation.
However, to impose the interface derivatives, we need to compute them. 
There are different methods to compute derivatives. 
For instance, we can apply an epsilon method \cite{Seibold2012}.
It consists of adding two interpolation points at $\pm \varepsilon$ of the interfaces and 
computing the derivatives using the function values at these points. 
To use this method, we then need to add two interpolation points for each derivative that
we want to compute. 
Another idea specific to the backward semi-Lagrangian method presented by E. Bourne in her doctoral thesis 
\cite{BournePhD} is to advect the derivatives simultaneously with the function.
Starting from the relation in the backward semi-Lagrangian method, 
$f(t, \mathbf{X}(t; s, \mathbf{x})) = f(s,\mathbf{x}), \forall s \in \mathbb{R}, \mathbf{x}\in \Omega$,
a relation for the derivatives can be directly deduced: 
$\partial_{\mathbf{x}} (f(t, \mathbf{X}(t; s, \mathbf{x})) )
= \partial_{\mathbf{x}} \mathbf{X}(t; s, \mathbf{x}) \partial_{\mathbf{x}} f(s,\mathbf{x}), 
\forall s \in \mathbb{R}, \mathbf{x}\in \Omega.$
However, this method additionally requires the derivatives of the characteristics to be determined
$\partial_{\mathbf{x}} \mathbf{X}(t; s, \mathbf{x})$. 
 
This article is devoted to the generalization of a formula developed by N. Crouseilles, 
G. Latu and E. Sonnendrücker, in \cite{Crouseilles2009}.
It consists of calculating the derivatives from the function values
for parallel computation purposes (see also \cite{Unser1991,MendozaPhD}).
It is also interesting for the multi-patch approach 
(even without parallelization considerations). 
This formula was developed for uniform meshes.
We use Hermite polynomials to handle non-uniform meshes and non-conforming situations.
The generalization introduced in this paper allows us to calculate the
derivatives of an equivalent global spline and thus reach $\mathcal{C}^2$ 
continuity in the conforming case. The use of derivatives equal to those of an equivalent global spline 
comes at the cost of the resolution of a small global system, but the derivatives can also be approximated.
The weight of the values at the most distant points in the calculation of the derivative decreases 
as the number of points per patch increases.

We can draw a parallel between this interface derivatives approach and the 
hybridizable discontinuous Galerkin (HDG) approach \cite{Carrero2006, Cockburn2016, Cockburn2009}. 
The two approaches add degrees of freedom and split the domain to locally solve a problem. 
The added interface derivatives in the multi-patch approach can be seen as elements 
of the numerical traces in the static condensation method.

This paper is organized as follows. Some reminders of the backward semi-Lagrangian method and 
the B-splines can be found in Section \ref{section_Semi-Lagrangian_global}. 
Section \ref{section_Semi-Lagrangian_1D_multi-patch} presents the generalized method that we 
suggest. To advect with a backward semi-Lagrangian method on 1D multi-patch geometries, 
we develop a formulation to calculate the derivatives of an equivalent global spline.
We also study an explicit formulation and its asymptotic behavior in the uniform per patch case. 
Section \ref{section_Semi-Lagrangian_2D_multi-patch} describes the application of the method on 2D 
multi-patch geometries.
Section \ref{section_Numerical_results} discusses some numerical results on spline interpolations,
advections, and the guiding-center equations run with the Gyselalib++ code on multi-patch geometries
with O-point. 
We apply the method only to geometries with an O-point using the work in \cite{Zoni2019,BournePhD}, 
but not to geometries with X-points as further studies are needed. 

\section{Semi-Lagrangian scheme on a global domain}
\label{section_Semi-Lagrangian_global}
\paragraph{}
We briefly recall the backward semi-Lagrangian method on a global domain without multi-patch in the 
first subsection \ref{subsection_BSL_method}. 
For the interpolation step, we use B-splines which are described in the second subsection 
\ref{subsection_Bsplines_interpolation}. 
These B-splines are defined on a logical domain and then mapped to a physical domain. These terms are 
defined in subsection \ref{subsection_logical_physical_domains}. 

\subsection{Backward semi-Lagrangian method}
\label{subsection_BSL_method}
\paragraph{}
The hyperbolic advection equation in  \eqref{guiding_center_eq} is solved with a backward semi-Lagrangian 
method. 
One advantage of semi-Lagrangian methods is that they are less restrictive on the time step, as there 
is no Courant–Friedrichs–Lewy (CFL) condition \cite{CFL}.
However, the method is not structure-preserving and, we may observe diffusion in long time simulation. 
The method is detailed in \cite{Sonnendrücker1999}.
For a typical advection equation as follows, 
\begin{equation}
    \partial_t f(t,\mathbf{x}) + \mathbf{A}(t,\mathbf{x})\cdot\nabla f(t,\mathbf{x}) = 0, \\
    \label{classical_advection_equation}
\end{equation}
with a given advection field $\mathbf{A}$, the method consists of adopting the Lagrangian point of view. 
We define the characteristics as the solution $\mathbf{X}$ of the following ordinary differential equation (ODE), 
\begin{equation*}
    \left\{
        \begin{aligned}
            & \partial_s \mathbf{X}(s; t, \mathbf{x}) =  \mathbf{A}(s, \mathbf{X}(s; t, \mathbf{x})), 
                &\qquad s\geq 0, \mathbf{x}\in\Omega,\\
            & \mathbf{X}(t; t, \mathbf{x}) = \mathbf{x}, 
                &\qquad \mathbf{x}\in\Omega.
        \end{aligned}
    \right. 
\end{equation*}

Thanks to the property of conservation along the characteristics, we have the following relation for 
the solution of \eqref{classical_advection_equation}, with fixed $t\in \mathbb{R}, \mathbf{x}\in\Omega$,
\begin{equation*}
    f(t, \mathbf{X}(s; t, \mathbf{x})) = f(t, \mathbf{x}), \qquad \forall t\in \mathbb{R}.
\end{equation*}

Knowing the function at $s = t-\Delta t$, we can deduce it at $t$. We need to solve backward
the equations of the characteristics for every interpolation point. The solutions of the ODE
are called the feet of the characteristics.
We interpolate the function known at $t-\Delta t$ at the feet of the characteristics.
Each function value corresponds to the function value at $t$ at the interpolation point $\mathbf{x}$. 
To adapt the method to a multi-patch geometry, the main difficulty is in the interpolation.

\subsection{B-splines interpolation}
\label{subsection_Bsplines_interpolation}
\paragraph{}
The chosen interpolation method for the function at $t$ is an interpolation on cubic B-splines. 
These basis functions are well described in \cite{NURBS, Patrikalakis2002, Farin2014, Unser1991}. 
In 1D space, the spline $s$ is given by 	
\begin{equation*}
    s(x) = \sum_{k=0}^{N_c+2} c_k b_k(x), \qquad x\in\Omega
\end{equation*}
with $N_c$ the number of cells, 
$\{b_k\}_{k = 0, ..., N_c+2}$ the basis functions, called the B-splines, 
and $\{c_k\}_{k = 0, ..., N_c+2}$ coefficients to determine
(also called \textit{control points}).

We call \textit{break points}, $\{x_{b,i}\}_{i=0}^{N_c}$ the boundaries of the cells.
There are $N_c+1$ break points. 
We add $N_c+6$ points, $\{k_i\}_{i=0}^{N_c+5}$ called \textit{knots} to build the B-splines. 
The central knot sequence is defined on the break points; the additional points are placed at 
regular intervals outside the domain (this is a common choice in the uniform case) or are defined 
with repetitions of the first and the last break points 
(see \ref{Appendix_knots_and_Bsplines} for more details). 
The latter case is referred to as an \textit{open knots sequence}.

Whatever knot sequence we choose, this system has $N_c+3$ degrees of freedom to fix. 
We then introduce interpolation points $\{x_i\}_{i\in N}$.
In the Gyselalib++ code, three interpolation types are implemented.

\paragraph{Interpolation with Hermite boundary conditions.} If the function values and its 
derivatives are known on the boundaries, a way to fix the $N_c+3$ degrees of freedom is to 
define $N_c+1$ interpolation points $\{x_i\}_{i=0}^{N_c}$ 
and fix the two last degrees of freedom with the derivatives at the boundaries as follows.

\begin{equation*}
    \left\{
        \begin{aligned}
            & f_i = \sum_{k=0}^{N_c+2} c_k b_k(x_i), 	
                \qquad i = 0, ..., N_c,\\ 
            & (\partial_x f)_i = \sum_{k=0}^{N_c+2} c_k b'_k(x_i), 	
                \qquad i = 0, N_c. \\ 			
        \end{aligned}
    \right. 
\end{equation*}

\paragraph{Interpolation using periodic boundary conditions.} If the derivatives are unknown on the boundaries, 
but the domain is periodic, then we can define $N_c+1$ interpolation points and use the periodicity 
property to fix all the degrees of freedom. 
\begin{equation*}
    f_i = \sum_{k=0}^{N_c+2} c_k b_k(x_i), 	 \qquad i = 0, ..., N_c.
\end{equation*}

\paragraph{Interpolation points as closure conditions.} 
If we do not know the derivatives at the boundaries and the domain 
is not periodic, we add more interpolation points to fix all the degrees of freedom. We generally 
use Greville points which are defined from the knot sequence as follows, (here $d = 3$).
In this paper, we use open knot sequences.
\begin{equation}
    x_i = \frac{1}{d}\sum_{j=i+1}^{i+d} k_{j}.
    \label{Greville_points}
\end{equation}

Using the Greville points as interpolation points is usually considered as an optimal solution 
for accuracy and efficiency, \cite{Johnson2005, Greville1964, Greville1969, Farin2014}. 
In the uniform case, the Greville points (except at the boundaries) are the break points for cubic splines.
They would be cell midpoints for even degree splines. 
With $N_c+6$ knots, it gives $N_c+3$ interpolation points $\{x_i\}_{i=0}^{N_c+2}$. 
This fixes all the degrees of freedom.
\begin{equation*}
    f_i = \sum_{k=0}^{N_c+2} c_k b_k(x_i), 	\qquad i = 0, ..., N_c+2 
\end{equation*}

Once the degrees of freedom are fixed, we can determine the B-spline coefficients $\{c_k\}_{k=0}^{N_c+2}$. 
They define a spline representation of the function at $t - \Delta t$ that we can evaluate to determine the 
function at $t$ in the backward semi-Lagrangian method, \ref{subsection_BSL_method}. 
This spline representation does not depend on the knot sequence used to define the B-splines. 

\paragraph{Remark.}
Using a uniform knot sequence provides identical B-splines except for a support shift. 
It is then easier to optimize the code. So, it is more interesting to split a global non-uniform 
domain into several uniform patches to optimize it.

\subsection{Logical and physical domains}
\label{subsection_logical_physical_domains}
\begin{figure}[h]
    \centering
    \begin{tikzpicture}[xscale = 0.75, yscale = 0.75]
        \begin{scope}[xshift=-1.5cm, yshift = 0cm]
            \draw[lightgray] (0,0) grid [step = 0.5] (3,3);
            \draw[thick] (0,0) grid [step = 3] (3,3);
            
            \draw[->] (0, 0) -- (3.3, 0) node [right]{$r$}; 
            \draw[->] (0, 0) -- (0, 3.3) node [above]{$\theta$}; 

            \node[above] at (1.5, 3.7){Logical domain $\widehat{\Omega}$}; 
        \end{scope}
                    
        \begin{scope}[xshift=9cm, yshift = 1.5cm]
            \def\Rmax{1.6} 
            \def\numTheta{10} 
            \def\a{1.5} 
            \def\b{1.75} 
            
            \foreach \r in {0,0.2666667,...,\Rmax} {
                \draw[gray] (0,0) ellipse ({\a-\r*(\a/\Rmax)} and {\b-\r*(\b/\Rmax)});
            }
            
            \foreach \theta in {0,40,...,360} {
                \def\radius{\a*\b / sqrt(\b*\b*cos(\theta)*cos(\theta) + \a*\a*sin(\theta)*sin(\theta))} 
                \draw[gray] (0,0) -- (\theta:{\radius});
            }
            
            \draw[black] (0,0) ellipse ({\a} and {\b});
            \draw[black] (0,0) -- (0:\a);
            
            \draw[->] (0, 0) -- ({\a*1.1}, 0) node [right]{$x$}; 
            \draw[->] (0, 0) -- (0, {\b*1.1}) node [above]{$y$}; 

            \node[above] at (0, 2.3){Physical domain $\Omega$}; 
        \end{scope}
        
        \draw[->] (3, 1.5) -- (6, 1.5); 
        \node at (4.5,2) {$\mathcal{F}$}; 
    \end{tikzpicture}
    \caption{\label{Logical_physical_domain}
        Logical and physical domains. We note $(r,\theta)$ the curvilinear coordinates 
        on the logical domain and $(x,y)$ the Cartesian coordinates on the physical domain. }
\end{figure}
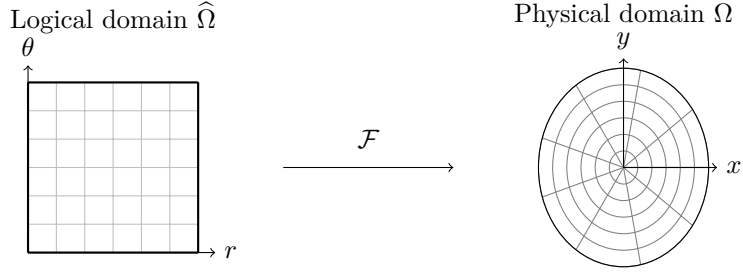

\paragraph{}
Let's introduce some notations. 
The B-splines are defined on a rectangle-shape domain called the \textit{logical domain} 
and denoted $\widehat{\Omega}$. 
The 2D B-spline basis is a tensor product of two 1D B-spline bases. This logical domain is mapped
to the \textit{physical domain} noted $\Omega$. 
An example of a physical domain and its associated logical domain is shown in Fig. \ref{Logical_physical_domain}. 
The mapping from the logical domain to the physical domain
and converting curvilinear coordinates into Cartesian coordinates is noted 
$\mathcal{F}: \widehat{\Omega} \mapsto \Omega$.

As the B-splines are defined on the logical domain, all the interpolations are performed in the logical 
domain. 

\section{Semi-Lagrangian scheme on a 1D multi-patch domain}
\label{section_Semi-Lagrangian_1D_multi-patch}
\paragraph{}
When we split a global domain into multiple patches, we need to be aware of the regularity at 
the patch interfaces. We start by defining the multi-patch domain in \ref{subsection_def_multipatch}. 
We explain our choice of a $\mathcal{C}^1$ regularity in \ref{subsection_local_spline_regularity_interface}.
The $\mathcal{C}^1$ regularity implies fixing the interface derivatives. 
The authors in \cite{Crouseilles2009} already gave us a formula 
to compute the derivatives with a precision of $10^{-6}$ for a uniform global mesh. Subsection 
\ref{subsection_computation_interface_derivatives} presents a new generalization of this formula 
to arbitrary meshes and any precision. 
We end this section with the algorithm for a backward semi-Lagrangian method adapted to a multi-patch mesh, 
\ref{subsection_BSL_multi_patch_1D}.

\subsection{Definition of multi-patch domain}
\label{subsection_def_multipatch}
\paragraph{}
We define a multi-patch domain as a physical domain split into different subdomains 
$\Omega = \bigcup_{k}\Omega^k$ 
(as the colors show in Fig. \ref{JET_SOLEDGE}). The interiors of these subdomains are disjoint. 
The domain recovery is only on the edges. 
We build equivalent logical domains $\widehat{\Omega}^k$ mapping 
to each physical subdomain. 
These logical domains are called \textit{patches}. 
The coordinates in the different patches are not necessarily equivalent. The logical domains are not 
necessarily defined in the same dimensions.
On each patch, we define independent spline representations of the function. We called these splines 
\textit{local splines}. For simplicity, all the local splines will be cubic piecewise polynomials. 

\subsection{Local spline regularity at the interfaces}
\label{subsection_local_spline_regularity_interface}
\paragraph{}
The local splines are defined independently of each other. However, we want the solution of
the guiding-center model, \eqref{guiding_center_eq} to be continuous in the global domain.  
We then need to impose a regularity on the interfaces between two patches to connect the local 
splines one another. 
The first idea is to impose a $\mathcal{C}^0$ regularity.

\subsubsection{$\mathcal{C}^0$ regularity}
\paragraph{}
If we impose a $\mathcal{C}^0$ regularity, then we do not need to define derivatives on the boundaries
of each local patch. We can define local splines with interpolation points as closure conditions.
As mentioned previously in 
\ref{subsection_Bsplines_interpolation}, optimal interpolation points are the Greville points. 
However, despite Greville points, 
we have noticed that using only a $\mathcal{C}^0$ regularity can be unstable. 

We start from the results in the article \cite{Crouseilles2011}.
In this article, they study a 1D advection on a periodic domain with a constant advection field. 
The domain $\widehat{\Omega}$ is split into $N_p$ patches. 
Unlike this article, 
we solve the advection with a backward semi-Lagrangian method \cite{Sonnendrücker1999}.
We choose to work with a uniform mesh for simplicity. 
The interpolation is realized with piecewise polynomial basis functions of degree $d=3$ over $N_c$ cells
for each patch,
\begin{equation}
    \varphi_j(x) = \sum_{k = 0}^d \varphi_{j,k} x^k, 
    \qquad j = 0, ..., N_c-1.
\end{equation}

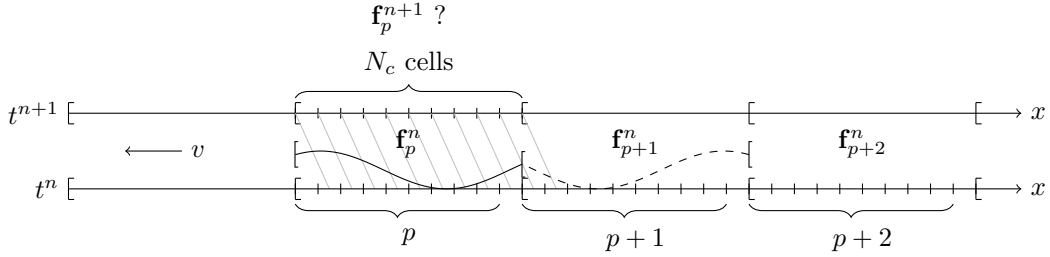
\begin{figure}[h]
    \centering
    \begin{tikzpicture}[xscale = 1.5]
        \draw[->] (0,0) -- (8.4,0) node[right]{$x$}; 
        \foreach \x in {0, 2, 4, 6, 8}{
            \draw (\x, -4pt) --(\x, +4pt); 
            \draw (\x, -4pt) --(\x+0.05, -4pt); 
            \draw (\x, +4pt) --(\x+0.05, +4pt); 
        }
        \foreach \x in {2, 2.2, ..., 4}{
            \draw (\x, -2pt) --(\x, +2pt); 
        }
        \draw [decorate, decoration={brace, amplitude=5pt}] (2,0.2) -- (4,0.2);
        \node at (3,0.4) [above]{$N_c$ cells};

        \foreach \x in {2, 2.2, ..., 4}{
            \draw[line width=0.2, lightgray] (\x,0) -- (\x+0.3,-1); 
        }

        \begin{scope}[yshift=-1cm]
            \draw[->] (0,0) -- (8.4,0) node[right]{$x$}; 
            \foreach \x in {0, 2, 4, 6, 8}{
                \draw (\x, -4pt) --(\x, +4pt); 
                \draw (\x, -4pt) --(\x+0.05, -4pt); 
                \draw (\x, +4pt) --(\x+0.05, +4pt);
            }
            \foreach \x in {2, 2.2, ..., 8}{
                \draw (\x, -2pt) --(\x, +2pt); 
            }

            \draw[black, domain=2:4, smooth, variable=\x, samples=50]
                plot ({\x}, {0.25+0.25*cos(\x*180*0.9)});	
            \draw[black, dashed, domain=2:4, smooth, variable=\x, samples=50]
                plot ({8-\x}, {0.25+0.25*cos(\x*180*0.9)});

            \node[black] at ({2}, {0.25+0.25*cos(2*180*0.9)}) {[}; 
            \node[black, dashed] at ({4}, {0.25+0.25*cos(4*180*0.9)}) {[}; 
            \node[black] at ({6}, {0.25+0.25*cos(2*180*0.9)}) {[}; 

            \draw [decorate, decoration={brace, amplitude=5pt}] (3.8,-0.2) -- (2,-0.2);
            \draw [decorate, decoration={brace, amplitude=5pt}] (5.8,-0.2) -- (4,-0.2);
            \draw [decorate, decoration={brace, amplitude=5pt}] (7.8,-0.2) -- (6,-0.2);
            \node at (3,-0.4) [below]{$p$};
            \node at (5,-0.4) [below]{$p+1$};
            \node at (7,-0.4) [below]{$p+2$};
        \end{scope}

        \node at (3,  1.25) {$\mathbf{f}^{n+1}_p$ ? }; 
        \node at (3, -0.4) {$\mathbf{f}^n_p$}; 
        \node at (5, -0.4) {$\mathbf{f}^n_{p+1}$}; 
        \node at (7, -0.4) {$\mathbf{f}^n_{p+2}$}; 

        \node at (0,0) [left]{$t^{n+1}$};
        \node at (0,-1) [left]{$t^n$};	
        \draw[<-] (0.5,-0.5) -- (1.,-0.5) node[right]{$v$}; 
    \end{tikzpicture}
    \caption{
        \label{fig_advection_Greville_C0}
        Constant advection on a periodic domain split into $N_p$ patches. 
        The function on patch $p$ is advected with a constant velocity $v<0$. 
        The feet of the characteristics land in the patch $p$ and $p+1$. To determine the 
        function at $t^{n+1}$, we interpolate 
        the function representation at $t^{n}$ using the information from patch $p$, 
        $p+1$ and $p+2$. 
    }
\end{figure}

The $\mathcal{C}^{d-1}$ continuity condition of the basis functions at the cell interfaces gives $d(N_c-1)$
constraints. 
The coefficients $\{\varphi_{j,k}\}_{j= 0, ..., N_c-1, k = 0, ..., d}$ are fixed 
with these constraints and the interpolation at the $N_c + d$ points of the function $f$ we want to advect in 
\eqref{classical_advection_equation}. 
To avoid duplication, each patch stores the $N_c + 2$ first function values and the 
last one is taken from the next patch.
An illustration of the advection step on these patches is given in Fig. \ref{fig_advection_Greville_C0}.
The function values $\mathbf{f}^{n+1}$ at $t^{n+1} = t^n + \Delta t^n$ are deduced from the function values 
$\mathbf{f}^{n}$ at $t^n$ for a given time discretization $\{t^n\}_{n\geq 0}$. 
We study the case of the constant velocity of the advection field $A = v$ 
such that $v<0$ and the displacement $|v\Delta t|$ is smaller than a patch length. 
We note $\mathbf{f}^{n}_p$ the set of function values on the patch $p$ of the function at $t^n$. 

The function values on patch $p$ at $t^{n+1}$ depend on the values of the function at $t^{n}$ on patches 
$p$, $p+1$ and $p+2$. As the advection is constant and the domain 
periodic, we can write $\mathbf{f}^{n+1} = \mathbf{A} \mathbf{f}^{n}$, with $\mathbf{A}$ a block matrix, 
\begin{equation}
    \mathbf{A} = 
    \begin{bmatrix}
        \mathbf{A}_0 & \mathbf{A}_1 & \mathbf{A}_2 & 0 & \dots & 0 \\
        0 & \mathbf{A}_0 & \mathbf{A}_1 & \mathbf{A}_2 & \ddots &  \vdots \\
        \vdots & \ddots & \ddots & \ddots & \ddots & 0 \\
        0 & \dots & 0 & \mathbf{A}_0 & \mathbf{A}_1 & \mathbf{A}_2 \\
        \mathbf{A}_2 & 0 & \dots & 0 & \mathbf{A}_0 & \mathbf{A}_1\\
        \mathbf{A}_1 & \mathbf{A}_2 & 0 & \dots & 0 & \mathbf{A}_0 \\
    \end{bmatrix}.
\end{equation}

It is a circulant matrix, we can therefore apply the Fourier transform to diagonalize it, 
\begin{equation}
    \left(\mathbb{R}^{Nc+2}\right)^{N_p} \rightarrow \left(\mathbb{C}^{Nc+2}\right)^{N_p}, 
    \quad
    (\mathbf{f}_0, \mathbf{f}_1, ..., \mathbf{f}_{N_p-1}) \mapsto 
        (\hat{\mathbf{f}}_0, \hat{\mathbf{f}}_1, ..., \hat{\mathbf{f}}_{N_p-1})
\end{equation}
with the relations
\begin{equation}
    \hat{\mathbf{f}}_p = \sum_{j=0}^{N_p-1} e^{-2i\pi p \ j/N_p} \mathbf{f}_j,
    \qquad
    \hat{\mathbf{A}}_p = \sum_{j=0}^{N_p-1} e^{-2i\pi p \ j/N_p} \mathbf{A}_j.
\end{equation}

However, it can be shown that for $k = 1$, $N_p = 5$ and $N_c = 3$, 
there exists a displacement $|v\Delta t| = 405/100000$, such that the maximum of 
the absolute value of the eigenvalues of $\hat{\mathbf{A}}_k$ is strictly superior to 1 
with an excess of $10^{-7}$. 
This counterexample shows that we do not have the eigenvalues 
strictly inferior to 1. We may encounter other examples with a higher excess, 
so a bigger instability.
So, this case shows that applying a $\mathcal{C}^0$ regularity can be unstable.

Thus, we need a stronger regularity between the patches. 

\subsubsection{$\mathcal{C}^1$ regularity}
\paragraph{}
As a $\mathcal{C}^0$ regularity may sometimes lead to instabilities,
we enforce a $\mathcal{C}^1$ regularity at the interfaces. 
We impose the equality of the derivatives and the spline representations of the function
at the interfaces. We then apply an interpolation with Hermite boundary conditions 
on the local splines to ensure this. 
However, we do not know the derivatives at the interfaces. 
To build the spline representations, we need to determine them. 
 
We suggest in this article 
a new generalization of a formula developed in \cite{Crouseilles2009}
which consists of calculating the derivatives from twenty function values. 
This formula is given for uniform meshes. 
We use Hermite polynomials to extend it to arbitrary meshes. 
The new generalization introduced in this paper allows us to calculate 
derivatives equal to those of an equivalent global spline and even reach $\mathcal{C}^2$ 
continuity in the conforming case.

\subsection{Computation of interface derivatives}
\label{subsection_computation_interface_derivatives}
\subsubsection{Assumptions}
\paragraph{}
Before starting the study, we describe the assumptions that are necessary for this work. 
The study conducted here supposes that the spline representations use cubic B-splines. 

Another remark concerns the interpolation points. As we will see, to define Hermite polynomials 
between two interpolation points, the spline representation between the points needs to be 
polynomial and not piecewise polynomial.
This means that each break point also needs to be an interpolation point. 
While Greville points are found at break points when uniform cells are used, this is not the case generally.
In order to fulfil this condition, we cannot use Greville points as interpolation points 
in the non-uniform case.

\subsubsection{Application of Hermite polynomials}
\paragraph{}
The authors in \cite{Crouseilles2009}
wanted to parallelize the computations in the semi-Lagrangian method using a domain decomposition. 
They split a global spline into local splines 
and used \eqref{Crouseilles_formula} to be as close as possible to the initial global spline. 
They developed a formula to compute the derivative $s'$ of the spline representation 
\eqref{spline_representation}
from the function values $\{f_{i+k}\}_k$ at an interpolation point $x_i$: 
\begin{equation}
    s'(x_i) = \sum_{j=-10}^{10} \frac{\widetilde{\omega}_j}{\Delta x} f_{i+j},
    \label{Crouseilles_formula}
\end{equation}
with given weights $\{\widetilde{\omega}_j\}_j$ ($\widetilde{\omega}_0 = 0$) 
(see Table 1 in \cite{Crouseilles2009}),
and $\Delta x$ the uniform length between the interpolation points.

This formula is interesting for our multi-patch geometry. 
It is an approximation of the spline derivative with a precision of $10^{-6}$. 
The coefficients $\{\widetilde{\omega}_k\}_k$ are given for a uniform mesh. 
However, it does not allow us to obtain better precision or deal with non-uniform cases. 
If we use this formula to compute the derivatives at the interface of two patches, 
both have to have the same uniformity. In practice, we would like to locally refine 
some patches. 
So, some interfaces will connect zones with two different cell lengths.
In the following sections, we will generalize this formula to non-uniform meshes 
and extend it to be able to use it with any precision.

\subsubsection{Relation between derivatives at three consecutive points}
\label{subsubsection_Relation_between_derivatives_at_three_consecutive_points}
\paragraph{}
To generalize the formula \eqref{Crouseilles_formula}, we define Hermite polynomials. 
The cubic Hermite polynomials interpolate the function values
and the function derivatives at the boundaries. 
They are defined as follows for $x\in[a,b]$, $a,b\in \mathbb{R}$,
\begin{equation}
    P (a + x(b-a)) 
        = f(a) H_0(x) 
        + f(b) H_1(x)
        + f'(a) (b-a) K_0(x)
        + f'(b) (b-a) K_1(x)
    \label{Hermite_polynomials}
\end{equation}
with the basis functions
\begin{equation}
    H_0(x) = (1-x)^2 (1+2x), 
    \qquad
    H_1(x) = x^2 (3-2x), 
    \qquad
    K_0(x) = (1-x)^2 x, 
    \qquad
    K_1(x) = x^2 (x-1). 
    \label{Hermite_basis_functions}
\end{equation}

\begin{figure}[h]
    \centering
    \begin{tikzpicture}[xscale = 1.5, yscale = 1.5]
        \draw (0,0) -- (2,0); 
        \foreach \x in {0,2}{
            \draw (\x, 0) -- (\x, 1); 
        }
        \def\coef{120}; 
        \draw[domain=0:2, smooth, variable=\x, samples=50]
                plot ({\x}, {0.25+0.25*cos(\x*\coef)});
        \foreach \x in {0, 2}{
            \draw[fill=black] ({\x}, {0.25+0.25*cos(\x*\coef)}) circle (0.75pt); 
        }

        \node at (1,1){$P_L$}; 
        \node at (0,-0.25){ $x^L_{i-1}$}; 
        \node at (2,-0.25){ $x^L_{i}$};
        \node at (-0.5,0.75){\footnotesize $s(x^L_{i-1})$}; 
        \node at (-0.5,1.00){\footnotesize $s'(x^L_{i-1})$}; 

        \node at (2.75,0.75){\footnotesize $s(x^L_{i})$ = $s(x^R_{i})$}; 
        \node at (2.75,1.00){\footnotesize $s'(x^L_{i})$ = $s'(x^R_{i})$};  

        \begin{scope}[xshift=1.5cm]
            \draw (2,0) -- (4,0); 
            \foreach \x in {2,4}{
                \draw (\x, 0) -- (\x, 1); 
            }
            \draw[domain=2:4, smooth, variable=\x, samples=50]
                plot ({\x}, {0.25+0.25*cos(\x*\coef)});	
            \foreach \x in {2, 4}{
                \draw[fill=black] ({\x}, {0.25+0.25*cos(\x*\coef)}) circle (0.75pt); 
            }

            \node at (3,1){$P_R$};
            \node at (2,-0.25){ $x^R_{i}$}; 
            \node at (4,-0.25){ $x^R_{i+1}$}; 
            \node at (4.5,0.75){\footnotesize $s(x^R_{i+1})$}; 
            \node at (4.5,1.00){\footnotesize $s'(x^R_{i+1})$}; 
        \end{scope}

    \end{tikzpicture}
    \caption{\label{2_Hermite_polynomials} Two Hermite polynomials: $P_L$ defined on $[x^L_{i-1}, x^L_{i}]$
    and $P_R$ defined on $[x^R_{i}, x^R_{i+1}]$ with functions connected across patches between 
    $x_i^L$ and $x_i^R$.}
\end{figure}
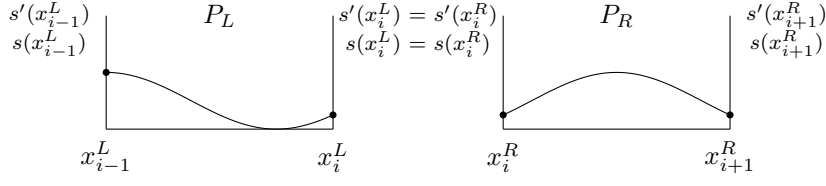

We start by defining two Hermite polynomials $P_{L}$ and $P_{R}$ on two domains defined using 
three consecutive interpolation points (see Fig. \ref{2_Hermite_polynomials}). 
The polynomials are fixed by the boundary conditions 
($s(x^L_{i-1})$, $s'(x^L_{i-1})$, $s(x^L_{i})$, $s'(x^L_{i})$ 
and $s(x^R_{i})$, $s'(x^R_{i})$, $s(x^R_{i+1})$, $s'(x^R_{i+1})$). 
Here, $s$ can represent the same spline on the two cells, or two different local splines. 
To stay as general as possible, we do not impose that the interpolation points 
$x^L_{i}$ and $x^R_{i}$ have the same value. 
This is particularly the case when we define functions on different logical domains 
which connect via the physical domain,but do not share coordinate systems.
The indices of $x^L_{i}$ and $x^R_{i}$ could also be different, e.g. $x^L_{i^L}$ and $x^R_{i^R}$.
But to lighten the notation, we will only write $i$. 

We also note in this section the $x$-dimension the dimension which we split into multiple patches. 
This reasoning can be applied equally well to both the $r$-dimension and the $\theta$-dimension of the patches.

We impose a $\mathcal{C}^2$ continuity between the two polynomials. Imposing it enforces the Hermite 
polynomials to be equal to an equivalent global cubic spline defined on the whole domain 
$[x^L_{i-1}, x^L_{i}]\bigcup  [x^R_{i}, x^R_{i+1}]$. 
Moreover, it gives us a relation between the derivatives.
Indeed, if we impose
\begin{equation*}
    P''_L(x^L_i) = P''_R(x^R_i),
\end{equation*}
with 
\begin{equation*}
    \left\{
    \begin{aligned}
        & P''_L(x^L_{i})
            = \frac{1}{(x^L_{i}-x^L_{i-1})^2} \left[
                s(x^L_{i-1}) H''_0(1) 
            + s(x^L_{i}) H''_1(1)
            + s'(x^L_{i-1}) (x^L_{i}-x^L_{i-1}) K''_0(1)
            + s'(x^L_{i}) (x^L_{i}-x^L_{i-1}) K''_1(1) 
            \right]\\
        & P''_R(x^R_{i}) 
            = \frac{1}{(x^R_{i+1}-x^R_{i})^2} \left[
                s(x^R_{i}) H''_0(0) 
            + s(x^R_{i+1}) H''_1(0)
            + s'(x^R_{i}) (x^R_{i+1}-x^R_{i}) K''_0(0)
            + s'(x^R_{i+1}) (x^R_{i+1}-x^R_{i}) K''_1(0) 
            \right]\\
    \end{aligned}
    \right.
\end{equation*}
and write $\Delta x^L_{i} = x^L_{i}-x^L_{i-1}$ and $\Delta x^R_{i} = x^R_{i+1}-x^R_{i}$, 
the relation $P''_L(x^L_i) = P''_R(x^R_i)$ becomes
\begin{equation*}
    \begin{aligned}
        & \frac{1}{(\Delta x^L_{i})^2}\left[
                6 s(x^L_{i-1})
            - 6 s(x^L_{i})
            + 2 \Delta x^L_{i} s'(x^L_{i-1}) 
            + 4 \Delta x^L_{i} s'(x^L_{i}) 
        \right]\\
        & = \frac{1}{(\Delta x^R_{i})^2}\left[
            - 6 s(x^R_{i})
            + 6 s(x^R_{i+1})
            - 4 \Delta x^R_{i} s'(x^R_{i}) 
            - 2 \Delta x^R_{i} s'(x^R_{i+1}) 
            \right]. \\
    \end{aligned}
\end{equation*}
Reordered knowing that $s'(x^L_{i}) = s'(x^R_{i})$ and $s(x^L_{i}) = s(x^R_{i})$ and that 
the values of the spline at the interpolation points are the function values, the relation is 
\begin{equation*}
    \begin{aligned}
        s'(x^L_{i})
        =& \frac{3}{2} \frac{1}{\Delta x^R_{i} +\Delta x^L_{i}} 
        \left[
                \frac{\Delta x^L_{i}}{\Delta x^R_{i}}  f^R_{i+1}
            + \left( 
                \frac{\Delta x^R_{i}}{\Delta x^L_{i}}  
                - \frac{\Delta x^L_{i}}{\Delta x^R_{i}} 
                \right) f^R_{i}
            - \frac{\Delta x^R_{i}}{\Delta x^L_{i}}  f^L_{i-1}
        \right] 
        \\ &
        -  
        \frac{1}{2} \frac{1}{\Delta x^R_{i} +\Delta x^L_{i}} 
        \left[
                \Delta x^L_{i} s'(x^R_{i+1}) 
            + \Delta x^R_{i} s'(x^L_{i-1}) 
        \right].
    \end{aligned}
\end{equation*}

We then obtain a relation between the derivatives at three consecutive interpolation points. 
This relation gives the exact derivative of the equivalent global spline defined on the whole 
domain $[x^L_{i-1}, x^L_{i}]\bigcup  [x^R_{i}, x^R_{i+1}]$.
For clarity, we introduce the coefficients $\gamma_i$, $\alpha_i$ and $\beta_i$ such that,
\begin{equation}
    s'(x^L_{i}) = \gamma_i + \alpha_i s'(x^R_{i+1}) + \beta_i s'(x^L_{i-1}) 
    \label{relation_3_interpol_points}
\end{equation}
with 
\begin{equation*}
    \left\{
    \begin{aligned}
        & \gamma_i = 
            \frac{3}{2} \frac{1}{\Delta x^R_{i} +\Delta x^L_{i}} 
            \left[
                \frac{\Delta x^L_{i}}{\Delta x^R_{i}}  f^R_{i+1}
                + \left( 
                    \frac{\Delta x^R_{i}}{\Delta x^L_{i}}  
                    - \frac{\Delta x^L_{i}}{\Delta x^R_{i}} 
                \right) f^R_{i}
                - \frac{\Delta x^R_{i}}{\Delta x^L_{i}}  f^L_{i-1}
            \right],
        \\ 
        & \alpha_i = -\frac{1}{2} \frac{\Delta x^L_{i}}{\Delta x^R_{i} +\Delta x^L_{i}}, \\
        & \beta_i =  -\frac{1}{2} \frac{\Delta x^R_{i}}{\Delta x^R_{i} +\Delta x^L_{i}}. \\
    \end{aligned}
    \right.
\end{equation*}

\paragraph{Remark 1.}
It is interesting to remark that the coefficients $\alpha_i$ and $\beta_i$ depend only on the 
proportionality of the distances between the interpolation points. The coefficient $\gamma_i$
is a linear combination of the function values and is inversely proportional to the distance 
between the interpolation points.

\paragraph{Remark 2.}
Note also that we wanted $\mathcal{C}^1$ regularity, and this formula gives us an improved 
$\mathcal{C}^2$ regularity.

\subsubsection{Extension of the relation}
\paragraph{}
Let's take the calculations a step further to 
extend this formula to a relation between derivatives further from each other. 
We first go in the right direction
and establish a relation between the derivatives at $x_{i-1}$, 
$x_{i}$ and $x_{i + N^R}$ with  $N^R\in\mathbb{N}^*$.
To simplify the notation, we remove the superscripts $\cdot^L$ and $\cdot^R$, 
knowing that the value of the points can still be different on the left and on the right 
but the values of the function and its derivatives are the same. 
We can apply a recursion using the relation  \eqref{relation_3_interpol_points} 
at different interpolation points. 

The notation of type $c^i_{n^l,n^r}$ 
corresponds to a relation between the derivative at $s'(x_{i})$, 
the derivative at $n^l$ points to the left $s'(x_{i-n^l})$ (i.e. spaced by $n^l$ cells), and 
the derivative at $n^r$ points to the right $s'(x_{i+n^r})$ (i.e. spaced by $n^r$ cells).

We can prove the following relation for 
the derivatives at $x_{i-1}$, $x_{i}$ and $x_{i + n}$ with $n\in \mathbb{N}^*$, 
\begin{equation*}
    s'(x_{i}) = c^i_{1,n} + a^i_{1,n} s'(x_{i+n}) + b^i_{1,n} s'(x_{i-1}).
\end{equation*}
with $c^i_{1,n}$, $a^i_{1,n}$ and $b^i_{1,n}$, coefficients depending on the chosen 
interpolation point $x_i$ and the number of cells $n$.

Indeed, for $n = 1$, it corresponds to the relation \eqref{relation_3_interpol_points}, 
and the coefficients are
\begin{equation}
    c^i_{1,1} = \gamma_i, 
    \qquad
    a^i_{1,1} = \alpha_i, 
    \qquad
    b^i_{1,1} = \beta_i.
    \label{coeff_forward_1}
\end{equation}
For $n = 2$, we use the relation \eqref{relation_3_interpol_points} at $x_{i}$
and $x_{i+1}$,
and inject the equation \eqref{relation_3_interpol_points} at $x_{i+1}$ into the equation 
\eqref{relation_3_interpol_points} at $x_{i}$ to get
\begin{equation*}
    \begin{aligned}
        & s'(x_{i})
            = \frac{1}{1 - \alpha_i  \beta_{i+1}} \left[ \gamma_i + \alpha_i \gamma_{i+1} \right]
            + \frac{\alpha_i \alpha_{i+1}}{1 - \alpha_i  \beta_{i+1}} s'(x_{i+2}) 
            + \frac{\beta_i}{1 - \alpha_i  \beta_{i+1}} s'(x_{i-1}),
    \end{aligned}
\end{equation*}
so the coefficients are 
\begin{equation}
    c^i_{1,2} = \frac{1}{1 - \alpha_i  \beta_{i+1}} \left[ \gamma_i + \alpha_i \gamma_{i+1} \right], 
    \qquad
    a^i_{1,2} = \frac{\alpha_i \alpha_{i+1}}{1 - \alpha_i  \beta_{i+1}},
    \qquad
    b^i_{1,2} = \frac{\beta_i}{1 - \alpha_i  \beta_{i+1}}. 
    \label{coeff_forward_2}
\end{equation}

For $n\geq2$, we apply the same method with the three following relations
\begin{equation*}
    \left\{
    \begin{aligned}
        & s'(x_{i})
            = c^i_{1,n-1} + a^i_{1,n-1} s'(x_{i+n-1}) + b^i_{1,n-1} s'(x_{i-1}), \\
        & s'(x_{i})
            = c^i_{1,n,1} + a^i_{1,n} s'(x_{i+n}) + b^i_{1,n} s'(x_{i-1}), \\
        & s'(x_{i+n}) 
            = \gamma_{i+n} + \alpha_{i+n} s'(x_{i+n+1}) + \beta_{i+n} s'(x_{i+n-1}). \\
    \end{aligned}
    \right.
\end{equation*}

Injecting the third into the second one gives us, 
\begin{equation*}
    s'(x_{i})
        = c^i_{1,n} 
        + a^i_{1,n} \gamma_{i+n} 
        + a^i_{1,n} \alpha_{i+n} s'(x_{i+n+1}) 
        + a^i_{1,n} \beta_{i+n} s'(x_{i+n-1})
        + b^i_{1,n} s'(x_{i-1}), 
\end{equation*}
and using the first relation to replace $s'(x_{i+n-1})$, 
\begin{equation*}
    \begin{aligned}
        s'(x_{i})
        =& \frac{1}{1 -  \beta_{i+n} \frac{a^i_{1,n}}{a^i_{1,n-1}} } 
        \left[  
                c^i_{1,n} 
            + a^i_{1,n} \gamma_{i+n} 
            - a^i_{1,n} \beta_{i+n} \frac{1}{a^i_{1,n-1}} c^i_{1,n-1}
        \right] \\
        &
        + \frac{a^i_{1,n} \alpha_{i+n}}{1 -  \beta_{i+n} \frac{a^i_{1,n}}{a^i_{1,n-1}} } s'(x_{i+n+1}) 
        + \frac{1}{1 -  \beta_{i+n} \frac{a^i_{1,n}}{a^i_{1,n-1}} } 
        \left[
            b^i_{1,n}
            - \beta_{i+n} \frac{a^i_{1,n} }{a^i_{1,n-1}} b^i_{1,n-1}
        \right] s'(x_{i-1}).    
    \end{aligned}
\end{equation*}

The coefficients at $n\geq2$ are then given by 
\begin{equation}
    \left\{
    \begin{aligned}
        & c^i_{1,n+1} = \frac{1}{1 -  \beta_{i+n} \frac{a^i_{1,n}}{a^i_{1,n-1}} } 
            \left[  
                c^i_{1,n} 
                + a^i_{1,n} \gamma_{i+n} 
                - \beta_{i+n} \frac{a^i_{1,n}}{a^i_{1,n-1}} c^i_{1,n-1}
            \right] \\
        & a^i_{1,n+1} = \frac{a^i_{1,n} \alpha_{i+n}}{1 -  \beta_{i+n} \frac{a^i_{1,n}}{a^i_{1,n-1}} }, \\
        & b^i_{1,n+1} = \frac{1}{1 -  \beta_{i+n} \frac{a^i_{1,n}}{a^i_{1,n-1}} } 
            \left[
                b^i_{1,n}
                - \beta_{i+n} \frac{a^i_{1,n} }{a^i_{1,n-1}} b^i_{1,n-1}
            \right].  \\
    \end{aligned}
    \right.
    \label{coeff_forward_n}
\end{equation}

Applying this recursion, we find the desired relation between the derivatives at $x_{i-1}$, 
$x_{i}$ and $x_{i + N^R}$,
\begin{equation*}
    s'(x_{i}) = c^i_{1,N^R} + a^i_{1,N^R} s'(x_{i+N^R}) + b^i_{1,N^R} s'(x_{i-1}).
\end{equation*}

As a second step, we handle the points to the left and 
establish a relation between the derivatives at $x_{i-N^L}$, 
$x_{i}$ and $x_{i + N^R}$ with $N^L \in\mathbb{N}^*$.
We can apply the same method going to the left with a new initial condition. 
We write the relation linking
the derivatives at $x_{i-m}$, $x_{i}$ and $x_{i+N^R}$ with $m\in \mathbb{N}^*$, 
\begin{equation*}
    s'(x_{i}) = c^i_{m,N^R} + a^i_{m,N^R} s'(x_{i+N^R}) + b^i_{m,N^R} s'(x_{i-m})
\end{equation*}
with $c^i_{m,N^R}$, $a^i_{m,N^R}$ and $b^i_{m,N^R}$ coefficients depending on 
the chosen interpolation point $x_i$, 
the number of right cells $N^R$, and the number of left cells $m$.

For $m = 1$, the coefficients are those given by the previous formula with the coefficients 
\eqref{coeff_forward_n}. 
For $m = 2$, using equations \eqref{relation_3_interpol_points} and \eqref{coeff_forward_n},
we have
\begin{equation}
    c^i_{2,N^R} = \frac{1}{1 - b^{i}_{1,N^R}  \alpha_{i-1}}
        \left[
            c^{i}_{1,N^R} + b^{i}_{1,N^R}  \gamma_{i-1}
        \right], 
    \qquad
    a^i_{2,N^R} = \frac{a^{i}_{1,N^R}}{1 - b^{i}_{1,N^R} \alpha_{i-1}}, 
    \qquad
    b^i_{2,N^R} = \frac{b^{i}_{1,N^R}   \beta_{i-1}}{1 - b^{i}_{1,N^R}  \alpha_{i-1}}. 
    \label{coeff_backward_2}
\end{equation}

By recursion, we can establish that for $m\geq2$, 
\begin{equation}
    \left\{
    \begin{aligned}
        & c^{i}_{m+1,N^R} = \frac{1}{1 - \alpha_{i-m} \frac{b^{i}_{m,N^R}}{b^{i}_{m-1,N^R}}}
            \left[
                c^{i}_{m,N^R}
                + b^{i}_{m,N^R}  \gamma_{i-m}
                - \alpha_{i-m} \frac{b^{i}_{m,N^R}}{b^{i}_{m-1,N^R}} c^{i}_{m-1,N^R}
            \right],
        \\
        & a^{i}_{m+1,N^R} = \frac{1}{1 - \alpha_{i-m} \frac{b^{i}_{m,N^R}}{b^{i}_{m-1,N^R}}}
            \left[
                a^{i}_{m,N^R} 
                - \alpha_{i-m} \frac{b^{i}_{m,N^R}}{b^{i}_{m-1,N^R}} a^{i}_{m-1,N^R} 
            \right],
        \\
        & b^{i}_{m+1,N^R} = \frac{b^{i}_{m,N^R} \beta_{i-m} }
            {1 - \alpha_{i-m} \frac{b^{i}_{m,N^R}}{b^{i}_{m-1,N^R}}}.
        \\
    \end{aligned}
    \right.
    \label{coeff_backward_m}
\end{equation}

With these two recursions \eqref{coeff_backward_m} for $m = N^L$ and \eqref{coeff_forward_n} for $n = N^R$, 
we have now a relation between further derivatives at the 
interpolating points $x_{i-N^L}$, $x_{i}$ and $x_{i + N^R }$ for a given number of cells
on the right and on the left $N^R,N^L\in\mathbb{N^*}$ of $x_i$,
\begin{equation}
    s'(x_{i}) = c^i_{N^L,N^R} + a^i_{N^L,N^R} s'(x_{i+N^R}) + b^i_{N^L,N^R} s'(x_{i-N^L}).
    \label{relation_NM}
\end{equation}

This equation gives the exact derivatives of an equivalent global spline defined with the 
same interpolation points and on the same cells.

\subsubsection{Application to the multi-patch geometry}
\paragraph{}
In the previous section, we have just established a relation between three derivatives
at three given interpolation points.
The interesting aspect of this formula is that it can be directly applied to a multi-patch geometry 
to get a relation between the derivatives at the interfaces.  

\begin{figure}[h]
    \centering
    \begin{tikzpicture}[xscale=1.5, yscale=1.5]
        \draw (0,0) -- (10,0); 
        \foreach \x in {0, 0.25, ..., 10}{
            \draw (\x, -2pt) -- (\x, 2pt); 
        }
        \foreach \x in {0, 2, 5, 7, 10}{
            \draw (\x, 0) -- (\x, 0.5); 
        }
        \node at (1   , 0.75){patch 1};
        \node at (3.5 , 0.75){patch 2};
        \node at (6   , 0.75){patch 3};
        \node at (8.5 , 0.75){patch 4};
        \def\coef{120}; 
        \draw[dashed, domain=0:2, smooth, variable=\x, samples=50]
                plot ({\x}, {0.25+0.25*cos(\x*\coef)});
        \draw[domain=2:5, smooth, variable=\x, samples=50]
                plot ({\x}, {0.25+0.25*cos(\x*\coef)});
        \draw[dotted, domain=5:7, smooth, variable=\x, samples=50]
                plot ({\x}, {0.25+0.25*cos(\x*\coef)});
        \draw[dashdotted, domain=7:10, smooth, variable=\x, samples=50]
                plot ({\x}, {0.25+0.25*cos(\x*\coef)});

        \node at (0 , -0.3){$\mathcal{X}_{0}$};
        \node at (2 , -0.3){$\mathcal{X}_{1}$};
        \node at (5 , -0.3){$\mathcal{X}_{2}$};
        \node at (7 , -0.3){$\mathcal{X}_{3}$};
        \node at (10, -0.3){$\mathcal{X}_{4}$};

        \draw[|->] (2 , -0.6) -- (5 , -0.6); 
        \node[below] at (3.5 , -0.6) {\footnotesize $N^R_1$ points}; 
        \draw[|->] (2 , -0.6) -- (0 , -0.6); 
        \node[below] at (1 , -0.6) {\footnotesize $N^L_1$ points};

        \draw[|->] (5 , -0.95) -- (7 , -0.95); 
        \node[below] at (6 , -0.95) {\footnotesize $N^R_2$ points}; 
        \draw[|->] (5 , -0.95) -- (2 , -0.95); 
        \node[below] at (3.5 , -0.95) {\footnotesize $N^L_2$ points};  

        \draw[|->] (7 , -1.3) -- (10 , -1.3); 
        \node[below] at (8.5 , -1.3) {\footnotesize $N^R_3$ points}; 
        \draw[|->] (7 , -1.3) -- (5 , -1.3); 
        \node[below] at (6 , -1.3) {\footnotesize $N^L_3$ points};  

        \node[above] at (0 , 0.5) {\footnotesize $s'(\mathcal{X}_0)$};
        \node[above] at (2 , 0.5) {\footnotesize $s'(\mathcal{X}_1)$};
        \node[above] at (5 , 0.5) {\footnotesize $s'(\mathcal{X}_2)$};
        \node[above] at (7 , 0.5) {\footnotesize $s'(\mathcal{X}_3)$};
        \node[above] at (10, 0.5) {\footnotesize $s'(\mathcal{X}_4)$};
    \end{tikzpicture}
    \caption{\label{multi-patch_split_4} 
        Example of a 1D multi-patch of the global domain split into $N_p = 4$ patches.
        The spline representation is composed of $N_p = 4$ local spline representations. 
        The relation \eqref{relation_interfaces} can be given at the interfaces $\mathcal{X}_1$, 
        $\mathcal{X}_2$ and $\mathcal{X}_3$ on the figure. 
    }
\end{figure}
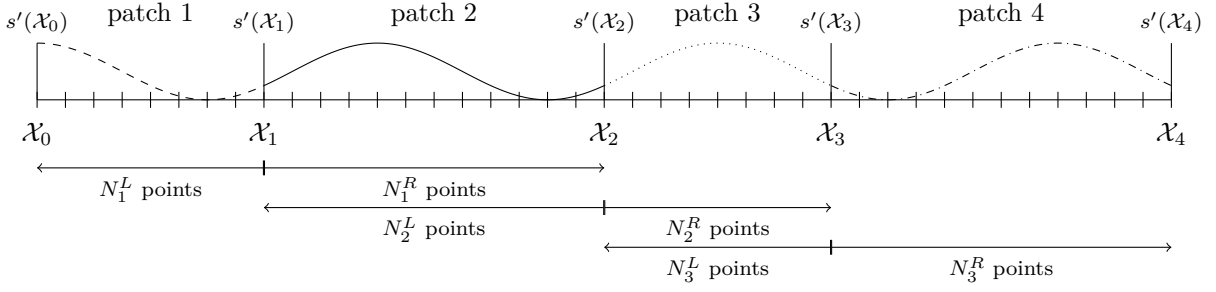

Let's split the global domain into $N_p$ patches and write the interface positions $\{\mathcal{X}_I\}_{I=0, ..., N_p}$
(see Fig. \ref{multi-patch_split_4}). 
The relation \eqref{relation_NM} is established at the interpolation point $x_i$ where an 
interface $\mathcal{X}_I$ is located.
$N^R$ corresponds to the number of cells between the interface $\mathcal{X}_{I}$ and $\mathcal{X}_{I+1}$. 
It is the number of cells on the patch to the right of the interface $\mathcal{X}_{I}$.
We set $N^L$ equal to the number of cells between the interface $\mathcal{X}_{I-1}$ and $\mathcal{X}_{I}$.
It is the number of cells on the patch to the left of the interface $\mathcal{X}_{I}$.
The relation \eqref{relation_NM} becomes a relation linking derivatives at three consecutive interfaces.
\begin{equation}
    s'(\mathcal{X}_I) 
        = c^I_{N^L_I,N^R_I} 
        + a^I_{N^L_I,N^R_I} s'(\mathcal{X}_{I+1}) 
        + b^I_{N^L_I, N^R_I} s'(\mathcal{X}_{I-1}).
    \label{relation_interfaces}
\end{equation}

We then have several relations for every three consecutive interfaces, 
\begin{equation*}
    \left\{
    \begin{aligned}
            s'(\mathcal{X}_1) 
                & = c^1_{N^L_1,N^R_1} 
                + a^1_{N^L_1,N^R_1} s'(\mathcal{X}_{2}) 
                + b^1_{N^L_1,N^R_1} s'(\mathcal{X}_{0}) \\
            s'(\mathcal{X}_2) 
                & = c^2_{N^L_2,N^R_2} 
                + a^2_{N^L_2,N^R_2} s'(\mathcal{X}_{3}) 
                + b^2_{N^L_2,N^R_2} s'(\mathcal{X}_{1}) \\
        & \vdots \\
            s'(\mathcal{X}_{N_p-1}) & = c^{N_p-1}_{N^L_{N_p-1},N^R_{N_p-1}} 
            + a^{N_p-1}_{N^L_{N_p-1},N^R_{N_p-1}} s'(\mathcal{X}_{{N_p}}) 
            + b^{N_p-1}_{N^L_{N_p-1},N^R_{N_p-1}} s'(\mathcal{X}_{{N_p-2}}). \\
    \end{aligned}
    \right.
\end{equation*}
We arrange them in the following matrix system, 
\begin{equation*}
    \begin{aligned}
        \mathbf{S} = 
        \begin{bmatrix}
            s'(\mathcal{X}_1) \\
            \vdots \\
            s'(\mathcal{X}_{I}) \\
            \vdots \\
            s'(\mathcal{X}_{N_p -1})
        \end{bmatrix}
        = \begin{bmatrix}
            1 & -b^{1}_{N^L_1,N^R_1} &  & \\
            -a^{2}_{N^L_2,N^R_2} & 1& -b^{2}_{N^L_2,N^R_2} & \\
            & \ddots& \ddots& \ddots\\
            & & -a^{N_p-1}_{N^L_{N_p-1},N^R_{N_p-1}}& 1\\
        \end{bmatrix}^{-1}
        \begin{bmatrix}
            c^{1}_{N^L_1,N^R_1} + a^{1}_{N^L_1,N^R_1} s'(\mathcal{X}_{0}) \\
            \vdots \\
            c^{I}_{N^L_I,N^R_I} \\
            \vdots \\
            c^{N_p-1}_{N^R_{N_p-1},N^L_{N_p-1}} + b^{N_p-1}_{N^R_{N_p-1},N^L_{N_p-1}} s'(\mathcal{X}_{N_p}) \\
        \end{bmatrix}
    \end{aligned}
\end{equation*}

\begin{equation}
    \mathbf{S} = (\mathbf{I} - \mathbf{M})^{-1} \mathbf{C}
    \label{Matrix_Hermite}
\end{equation}
with $s'(\mathcal{X}_{0})$ and $s'(\mathcal{X}_{N_p})$ the boundary conditions, supposed to be known. 

\paragraph{Remark.}
The size of the matrix is only $(N_p-1)\times (N_p-1)$. It corresponds to the number of inner interfaces. 
The matrix is small, so inverting it is not costly. We will also see in Section 
\ref{subsubsection_Approximation_of_the_relation}
that it is a dominant diagonal matrix.

\subsubsection{Independence of coefficients from the function}
\label{subsubsection_independence_of_coef_from_the_function}
\paragraph{}
As mentioned previously the coefficients $\gamma_i$ in Equation \eqref{relation_3_interpol_points} are 
a weighted sum of function values. The only coefficient in Equation \eqref{relation_interfaces} 
depending on the coefficients $\gamma_i$
is the coefficient $c^I_{N^L_I, N^R_I}$. It is also a weighted sum of function values. 
We can write it as follows
(with the index $i$ referring to the interpolation points $x_i$ located on the interface $\mathcal{X}_I$),
\begin{equation}
    c^I_{N^L_I, N^R_I} = \sum_{k= -N^L}^{N^R} \omega_{k,N^L_I, N^R_I} f_{i+k}.
    \label{c_weighted_sum}
\end{equation}

The weights $\{\omega_{k,N^L_I, N^R_I}\}_k$ only depend on the interpolation points. 
They can be determined at the initialization for each interface. 
In addition, the coefficients $a^I_{N^L_I, N^R_I}$ and $b^I_{N^L_I, N^R_I}$ in the relation 
\eqref{relation_interfaces} only depend
on $\alpha_i$ and $\beta_i$ from the relation \eqref{relation_3_interpol_points}. 
They then depend only on the proportionality of the distances between the interpolation points, 
(i.e. rescaling the domain will not affect the values of $a^I_{N^L_I, N^R_I}$ and $b^I_{N^L_I, N^R_I}$).

The computation of the coefficients $a^I_{N^L_I, N^R_I}$ and $b^I_{N^L_I, N^R_I}$ and 
the weights $\{\omega_k\}_k$ needs only to be done at the initialization for each interface. 
It is also done pseudo-locally: only the information from two patches is needed to determine them.

\subsubsection{Different boundary conditions}
\paragraph{}
To establish the matrix system \eqref{Matrix_Hermite}, we assumed that the derivatives on the 
external boundaries (i.e. the boundaries of the global domain) were known. 
In practice, we do not always know these derivatives. 
The two main examples are the case where 
we have defined enough interpolation points to not require additional closure conditions on the global domain, 
and the case where the global domain is periodic. 

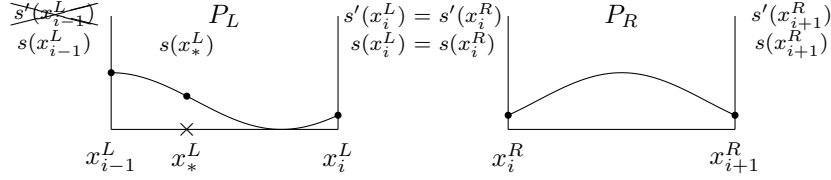
\begin{figure}[h]
    \centering
    \begin{tikzpicture}[xscale = 1.5, yscale = 1.5]
        \draw (0,0) -- (2,0); 
        \foreach \x in {0,2}{
            \draw (\x, 0) -- (\x, 1); 
        }
        \def\coef{120}; 
        \draw[domain=0:2, smooth, variable=\x, samples=50]
                plot ({\x}, {0.25+0.25*cos(\x*\coef)});

        \foreach \x in {0, 0.66667, 2}{
            \draw[fill=black] ({\x}, {0.25+0.25*cos(\x*\coef)}) circle (0.75pt); 
        }

        \node at (0.66667, 0) {$\times$}; 

        \node at (1,1){$P_L$};
        \node at (0,-0.25){ $x^L_{i-1}$}; 
        \node at (2,-0.25){ $x^L_{i}$};
        \node at (0.66667,-0.25){ $x^L_{*}$};
        \node at (-0.5, 0.75){\footnotesize $s(x^L_{i-1})$}; 
        \node at (-0.5, 1.00){\footnotesize $\xcancel{s'(x^L_{i-1})}$}; 
        \node at (0.66667, 0.75){\footnotesize ${s(x^L_{*})}$}; 

        \node at (2.75,0.75){\footnotesize $s(x^L_{i})$
            = $s(x^R_{i})$}; 
        \node at (2.75,1.00){\footnotesize $s'(x^L_{i})$
            = $s'(x^R_{i})$};  

        \begin{scope}[xshift=1.5cm]
            \draw (2,0) -- (4,0); 
            \foreach \x in {2,4}{
                \draw (\x, 0) -- (\x, 1); 
            }
            \draw[domain=2:4, smooth, variable=\x, samples=50]
                plot ({\x}, {0.25+0.25*cos(\x*\coef)});	
            \foreach \x in {2, 4}{
                \draw[fill=black] ({\x}, {0.25+0.25*cos(\x*\coef)}) circle (0.75pt); 
            }

            \node at (3,1){$P_R$}; 
            \node at (2,-0.25){ $x^R_{i}$}; 
            \node at (4,-0.25){ $x^R_{i+1}$}; 
            \node at (4.5,0.75){\footnotesize $s(x^R_{i+1})$}; 
            \node at (4.5,1.00){\footnotesize $s'(x^R_{i+1})$}; 
        \end{scope}

    \end{tikzpicture}
    \caption{\label{2_Hermite_polynomials_Greville} 
        Example of the case of interpolation points as closure conditions with
        two Hermite polynomials: $P_L$ defined on $[x^L_{i-1}, x^L_{i}]$
        and $P_R$ defined on $[x^R_{i}, x^R_{i+1}]$ with functions connected across patches 
        between $x_i^L$ and $x_i^R$.
        The left cell is a boundary cell of a global domain. We do not know the derivative at $x^L_{i-1}$.
        But, we know the function value at $x^L_{*}$ to provide closure for the system of equations. 
    }
\end{figure}

\paragraph{Interpolation points as closure conditions.}
If we define $N_c +3$ interpolation points, we do not need the derivatives, and we do not know them. 
All the degrees of freedom can be fixed. 
If we use Greville points as interpolation points, 
the previous figure, Fig. \ref{2_Hermite_polynomials} becomes the Fig. \ref{2_Hermite_polynomials_Greville} 
for the left boundary.
To determine the cubic Hermite polynomials $P_L$, we use the value of the function at this additional 
interpolation point $x^L_*$ in Fig. \ref{2_Hermite_polynomials_Greville}.

The Hermite polynomial $P_L$ has to interpolate the function at the additional interpolation point $x^L_{*}$, 
so we have the following relation
\begin{equation*}
    \begin{aligned}
        s(x^L_{*}) 
            =& s(x^L_{i-1}) H_0^{*,L}
            + s(x^L_{i})    H_1^{*,L}
            + s'(x^L_{i-1}) \Delta x^L_{i} K_0^{*,L} 
            + s'(x^L_{i})   \Delta x^L_{i} K_1^{*,L} \\
    \end{aligned}
\end{equation*}
with $H_0^{*,L}$, $H_1^{*,L}$, $K_0^{*,L}$ and $K_1^{*,L}$ the evaluations of 
$H_0$, $H_1$, $K_0$ and $K_1$ \eqref{Hermite_basis_functions} at 
$\frac{x^L_{*} - x^L_{i-1}}{\Delta x^L_{i}}$.

As $s'(x^L_{i-1})$ is unknown, we define it from the known values,  
and inject it into the relation \ref{relation_3_interpol_points} 
established before. This gives, 
\begin{equation*}
    s'(x^L_{i}) = 
    \frac{1}{1 + \beta_i \frac{K_1^{*,L}}{K_0^{*,L}}}
    \left[
        \gamma_i 
        + \frac{\beta_i }{\Delta x^L_{i} K_0^{*,L}} s(x^L_{*})
        - \frac{\beta_i H_0^{*,L}}{\Delta x^L_{i} K_0^{*,L}} s(x^L_{i-1}) 
        - \frac{\beta_i H_1^{*,L}}{\Delta x^L_{i} K_0^{*,L}} s(x^L_{i})  
    \right]
    +  \frac{\alpha_i}{1 + \beta_i \frac{K_1^{*,L}}{K_0^{*,L}}}s'(x^R_{i+1}).
\end{equation*}

So the relation on $[x^L_{i-1}, x^L_{i}]$ is 
\begin{equation}
    s'(x^L_{i}) = \gamma_i^{*,L} + \alpha_i^{*,L} s'(x^R_{i+1}).
    \label{relation_3_interpol_points_Greville_L}
\end{equation}
with 
\begin{equation*}
    \left\{
    \begin{aligned}
        & \gamma_i^{*,L} = \frac{1}{1 + \beta_i \frac{K_1^{*,L}}{K_0^{*,L}}}
            \left[
                \gamma_i 
                + \beta_i \frac{1}{\Delta x^L_{i} K_0^{*,L}} s(x^L_{*})
                - \beta_i \frac{H_0^{*,L}}{\Delta x^L_{i} K_0^{*,L}} s(x^L_{i-1}) 
                - \beta_i \frac{H_1^{*,L}}{\Delta x^L_{i} K_0^{*,L}} s(x^L_{i})  
            \right],
        \\
        & \alpha_i^{*,L} = \frac{\alpha_i}{1 + \beta_i \frac{K_1^{*,L}}{K_0^{*,L}}}, \\
        & \beta_i^{*,L} = 0.
    \end{aligned}
    \right. 
\end{equation*}

If the additional Greville point is on $[x^R_{i}, x^R_{i+1}]$, the relation
is similar
\begin{equation}
    s'(x^R_{i}) = \gamma_i^{*,R} + \beta_i^{*,R} s'(x^L_{i-1}).
    \label{relation_3_interpol_points_Greville_R}
\end{equation}
with 
\begin{equation*}
    \left\{
    \begin{aligned}
        & \gamma_i^{*,R} = \frac{1}{1 + \alpha_i \frac{K_0^{*,R}}{K_1^{*,R}}}
            \left[
                \gamma_i 
                + \alpha_i \frac{1}{\Delta x^R_{i} K_1^{*,R}} s(x^R_{*})
                - \alpha_i \frac{H_0^{*,R}}{\Delta x^R_{i} K_1^{*,R}} s(x^R_{i}) 
                - \alpha_i \frac{H_1^{*,R}}{\Delta x^R_{i} K_1^{*,R}} s(x^R_{i+1})  
            \right],
        \\
        & \alpha_i^{*,R} = 0, \\
        & \beta_i^{*,R} = \frac{\beta_i}{1 + \alpha_i \frac{K_0^{*,R}}{K_1^{*,R}}}.
    \end{aligned}
    \right. 
\end{equation*}
and 
$H_0^{*,R}$, $H_1^{*,R}$, $K_0^{*,R}$ and $K_1^{*,R}$ the evaluations of 
$H_0$, $H_1$, $K_0$ and $K_1$ \eqref{Hermite_basis_functions} at 
$\frac{x^R_{*} - x^R_{i}}{\Delta x^R_{i}}$.

In the recursions \eqref{coeff_forward_n} and \eqref{coeff_backward_m}, 
we use equations \eqref{relation_3_interpol_points_Greville_L} 
and/or \eqref{relation_3_interpol_points_Greville_R} instead of \eqref{relation_3_interpol_points}
for the last cells. 
This change affects only the first and the last interfaces for 
an equivalent global spline with interpolation points as closure conditions. 

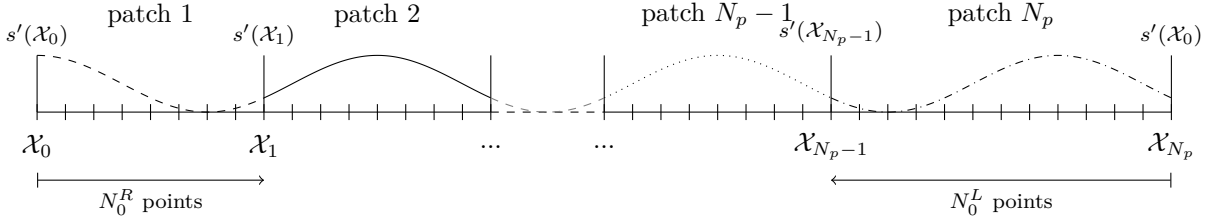
\begin{figure}[h]
    \centering
    \begin{tikzpicture}[xscale=1.5, yscale=1.5]
        \draw (0,0) -- (4,0); 
        \draw[dashed] (4,0) -- (5,0);
        \draw (5,0) -- (10,0);
        \foreach \x in {0, 0.25, ..., 4, 5, 5.25, ..., 10}{
            \draw (\x, -2pt) -- (\x, 2pt); 
        }

        \foreach \x in {0, 2, 4, 5, 7, 10}{
            \draw (\x, 0) -- (\x, 0.5); 
        }
        \node at (1   , 0.85){patch 1};
        \node at (3   , 0.85){patch 2};
        \node at (6   , 0.85){patch $N_p-1$};
        \node at (8.5 , 0.85){patch $N_p$};

        \def\coef{120}; 
        \draw[dashed, domain=0:2, smooth, variable=\x, samples=50]
                plot ({\x}, {0.25+0.25*cos(\x*\coef)});
        \draw[domain=2:4, smooth, variable=\x, samples=50]
                plot ({\x}, {0.25+0.25*cos(\x*\coef)});
        \draw[gray, dashed, domain=4:5, smooth, variable=\x, samples=50]
                plot ({\x}, {0.25+0.25*cos(\x*\coef)});
        \draw[dotted, domain=5:7, smooth, variable=\x, samples=50]
                plot ({\x}, {0.25+0.25*cos(\x*\coef)});
        \draw[dashdotted, domain=7:10, smooth, variable=\x, samples=50]
                plot ({\x}, {0.25+0.25*cos(\x*\coef)});

        \node at (0 , -0.3){$\mathcal{X}_{0}$};
        \node at (2 , -0.3){$\mathcal{X}_{1}$};
        \node at (4 , -0.3){$...$};
        \node at (5 , -0.3){$...$};
        \node at (7 , -0.3){$\mathcal{X}_{N_p-1}$};
        \node at (10, -0.3){$\mathcal{X}_{N_p}$};

        \draw[|->] (0 , -0.6) -- (2 , -0.6); 
        \node[below] at (1 , -0.6) {\footnotesize $N^R_0$ points}; 
        \draw[|->] (10 , -0.6) -- (7 , -0.6); 
        \node[below] at (8.5 , -0.6) {\footnotesize $N^L_0$ points};  

        \node[above] at (0 , 0.5) {\footnotesize $s'(\mathcal{X}_0)$};
        \node[above] at (10, 0.5) {\footnotesize $s'(\mathcal{X}_0)$};
        \node[above] at (2, 0.5) {\footnotesize $s'(\mathcal{X}_1)$};
        \node[above] at (7, 0.5) {\footnotesize $s'(\mathcal{X}_{N_p-1})$};
    \end{tikzpicture}
    \caption{\label{multi-patch_split_4_periodic} 
        Domain split into $N_p$ patches. 
        In the periodic case, the derivatives at $\mathcal{X}_0$, $\mathcal{X}_1$, $\mathcal{X}_{N_p}$ 
        and $\mathcal{X}_{N_p-1}$ depend on each other. 
        Moreover, the interfaces at $\mathcal{X}_0$ and $\mathcal{X}_{N_p-1}$ are physically overlapping. 
    }
\end{figure}

\paragraph{Interpolation with periodic boundary conditions.}
If the boundary conditions are periodic, we do not know the derivatives at the boundaries, but 
we can establish additional relations between the derivatives at $\mathcal{X}_0$, $\mathcal{X}_1$, 
$\mathcal{X}_{N_p}$ and $\mathcal{X}_{N_p-1}$ 
(see Fig. \ref{multi-patch_split_4_periodic}), 

\begin{equation*}
    \left\{
    \begin{aligned}
            s'(\mathcal{X}_0) & = s'(\mathcal{X}_{N_p}), \\
            s'(\mathcal{X}_0) 
                & = c^0_{N^L_0,N^R_0} 
                + a^0_{N^L_0,N^R_0} s'(\mathcal{X}_{1}) 
                + b^0_{N^L_0,N^R_0} s'(\mathcal{X}_{N_p-1}). \\
    \end{aligned}
    \right.
\end{equation*}

We add these relations to the matrix system \eqref{Matrix_Hermite}.  
It adds another line, so the matrix is now of size $N_p\times N_p$. 

\begin{equation*}
    \begin{aligned}
        \mathbf{S} = 
        \begin{bmatrix}
            s'(\mathcal{X}_0) \\
            s'(\mathcal{X}_1) \\
            \vdots \\
            s'(\mathcal{X}_{I}) \\
            \vdots \\
            s'(\mathcal{X}_{N_p -1})
        \end{bmatrix}
        = \begin{bmatrix}
            1 & -b^{0}_{N^L_0,N^R_0} & & & -a^{0}_{N^L_0,N^R_0} \\
            -a^{1}_{N^L_1,N^R_1} & 1 & -b^{1}_{N^L_1,N^R_1} & & \\
            & -a^{2}_{N^L_2,N^R_2} & 1& -b^{2}_{N^L_2,N^R_2} & \\
            & & \ddots& \ddots& \ddots\\
            -b^{N_p-1}_{N^L_{N_p-1},N^R_{N_p-1}}& & & -a^{N_p-1}_{N^L_{N_p-1},N^R_{N_p-1}}& 1\\
        \end{bmatrix}^{-1}
        \begin{bmatrix}
            c^{0}_{N^L_0,N^R_0} \\
            c^{1}_{N^L_1,N^R_1} \\
            \vdots \\
            c^{I}_{N^L_I,N^R_I} \\
            \vdots \\
            c^{N_p-1}_{N^L_{N_p-1},N^R_{N_p-1}} \\
        \end{bmatrix}.
    \end{aligned}
\end{equation*}

\subsubsection{Explicit formulation for uniform cells per patch}
\paragraph{}
As mentioned before, for optimization reasons, it is more interesting to deal with 
uniform patches. 
To locally refine the global domain, some patches may also not have the 
same cell length.  
In this situation, we can also optimize the computation by giving an explicit expression for 
the coefficients $a^I_{N^L_I, N^R_I}$, $b^I_{N^L_I, N^R_I}$ and $c^I_{N^L_I, N^R_I}$, 
instead of the recursive relations \eqref{relation_interfaces}, \eqref{Matrix_Hermite} given 
for the general case. 
The coefficients in 
\eqref{relation_3_interpol_points} are simplified and are almost constant in the uniform per patch case.  
Let's take two uniform patches surrounding the interface $\mathcal{X}_I$ at the interpolation point $x_i$.
We introduce the following sequence $u$ such that 
\begin{equation*}
    u_{k} = (2+\sqrt{3})^k - (2-\sqrt{3})^k, 
    \qquad \forall k \in \mathbb{N}.
\end{equation*}

We can then prove by recursion that the coefficients in \eqref{relation_NM}
can be explicitly expressed as follows with $\{u_{k}\}_{k \in\mathbb{N}}$,
for $n,m\in\mathbb{N}^*$, 
(with $n$ representing $N^R_I$ and  $m$ representing $N^L_I$ to simplify the notation), 
\begin{equation}
\left\{
\begin{aligned}
    & a^I_{m,n} 
        = \frac{(-1)^{n-1}u_{1}a^I_{1,1}u_{m}}
               {u_{m}u_{n} + u_{m}u_{n-1}a^I_{1,1} + u_{n} u_{m-1}b^I_{1,1}},
    \\
    & b^I_{m,n} 
        = \frac{(-1)^{m-1} u_{1} b^I_{1,1} u_{n}}
               {u_{m}u_{n} + u_{m}u_{n-1}a^I_{1,1} + u_{n} u_{m-1}b^I_{1,1}},
    \\
    & c^I_{m,n} = \sum_{k=-m}^{n} \omega^I_{k,m,n} f_{i+k},
\end{aligned}
\right.
\label{explicit_relation}
\end{equation}
with the following weights, 
\begin{equation}
\left\{
\begin{aligned}
    & \omega^I_{k,m,n} = 3(-1)^{k}
        \frac{\frac{a^I_{1,1}}{\Delta x^R} u_{m}u_{1}}
        {u_{m}u_{n} + u_{m}u_{n-1}a^I_{1,1} + u_{n} u_{m-1}b^I_{1,1}}
        , &k = n, \\
    & \omega^I_{k,m,n} = 3(-1)^{k}
        \frac{\frac{a^I_{1,1}}{\Delta x^R} u_{m}(u_{n-k+1} - u_{n-k-1})}
        {u_{m}u_{n} + u_{m}u_{n-1}a^I_{1,1} + u_{n} u_{m-1}b^I_{1,1}}
        , &k = 1, ..., n-1, \\
    & \omega^I_{k,m,n} = 3(-1)^k
        \frac{\frac{a^I_{1,1}}{\Delta x^R}u_{m} (u_{n} - u_{n-1}) 
            - \frac{b^I_{1,1}}{\Delta x^L}u_{n}(u_{m} - u_{m-1})}
        {u_{m}u_{n} + u_{m}u_{n-1}a_{1,1} + u_{n} u_{m-1}b^I_{1,1}}
        , &k = 0, \\
    & \omega^I_{k,m,n} = 3(-1)^{k+1}
        \frac{\frac{b^I_{1,1}}{\Delta x^L} u_{n}(u_{m+k+1} - u_{m+k-1})}
        {u_{m}u_{n} + u_{m}u_{n-1}a^I_{1,1} + u_{n} u_{m-1}b^I_{1,1}}
        , &k = -(m-1), ..., -1, \\
    & \omega^I_{k,m,n} = 3(-1)^{k+1}
        \frac{\frac{b^I_{1,1}}{\Delta x^L} u_{n}u_{1}}
        {u_{m}u_{n} + u_{m}u_{n-1}a^I_{1,1} + u_{n} u_{m-1}b^I_{1,1}}
        , &k = -m, \\
\end{aligned}
\right.
\label{explicit_weights}
\end{equation}
and $a^I_{1,1} = -\frac{1}{2} \frac{\Delta x^L}{\Delta x^R +\Delta x^L}$ 
and $b^I_{1,1} =  -\frac{1}{2} \frac{\Delta x^R}{\Delta x^R +\Delta x^L}$.
Some hints of the proof are given in \ref{Appendix_proof_explicit_formulation}

\subsubsection{Approximation of the relation}
\label{subsubsection_Approximation_of_the_relation}
\paragraph{}
Equation \eqref{Matrix_Hermite} gives derivatives which are exactly equal to the derivatives of an 
equivalent global spline. To compute them, we need the information at every interface.
However, when the number of cells increases, 
we expect that the derivative at one interface will depend less on the derivatives 
at the other interfaces. 
We see here how the coefficients of Equation \eqref{Matrix_Hermite} behave when the number of 
cells becomes large. If they become negligible, the dependency then also becomes
negligible. 

In this section, we work with the numbers of cells $N^L_I$ and $N^R_I$ tending to infinity.
Thanks to the symmetry of the problem, we can fix $N_I = N^R_I = N^L_I$
and make $N_I$ tend to infinity.

\begin{figure}[h]
    \centering
    \includegraphics[width=\textwidth]{"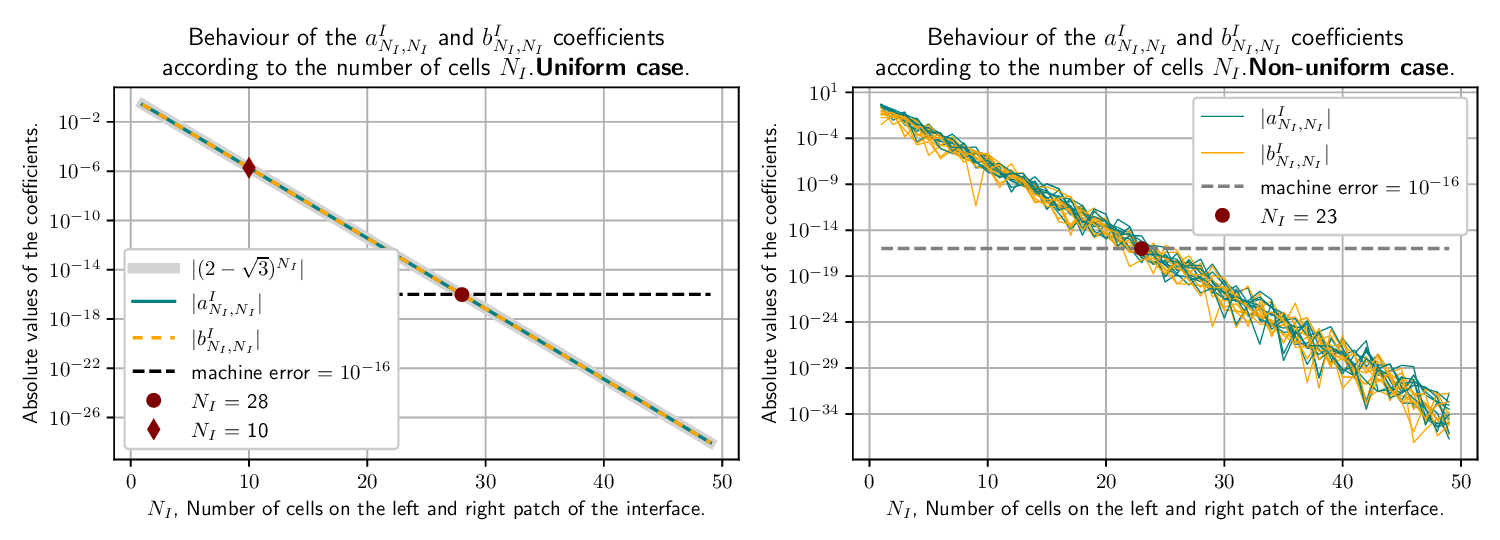"}
    \caption{\label{a_and_b_coeff_num} 
    Absolute values of the coefficients $a^I_{N_I, N_I}$ and $b^I_{N_I, N_I}$ for the same number 
    of cells $N_I$ in the two patches surrounding the interface $I$ with $n=N_I$
    varying between 1 and 50. On the left, it is for uniform meshes. 
    On the right several examples for non-uniform meshes generated with random functions. 
    We can observe that for $N_I \simeq 28$, the values of these coefficients reach machine 
    precision. The dependency on the derivatives at the other interfaces becomes negligible. 
    On the left, we compare the coefficients to the equivalence relation \eqref{equivalence_ab}.
    We can also note that for $N_I = 10$ in the uniform case, we retrieve the precision of 
    $10^{-6}$ mentioned in \cite{Crouseilles2009} for Equation \eqref{Crouseilles_formula}.
    }
\end{figure}

\paragraph{Behavior of the coefficients $a^I_{N_I,N_I}$ and $b^I_{N_I,N_I}$.}
These coefficients describe the dependency on the derivatives at the previous and next interfaces.
Their absolute values for an increasing number of interpolation points
are shown in Fig. \ref{a_and_b_coeff_num}. 

In the uniform per patch case, the expressions in \eqref{explicit_relation} give us an equivalence. 
Indeed, it is easy to see that $u_{k} = (2+\sqrt{3})^k + o((2+\sqrt{3})^k)$ when $k$ tends to infinity. 
So, injecting it in the equations, we find, for $n=N_I$ tending to infinity, 
\begin{equation}
\left\{
\begin{aligned}
    & a^I_{n,n} \simeq (-(2-\sqrt{3}))^n (-4a^I_{11}) \\
    & b^I_{n,n} \simeq (-(2-\sqrt{3}))^n (-4b^I_{11}) \\
\end{aligned}
\right.
\label{equivalence_ab}
\end{equation}
These results agree with the tests shown in Fig. \ref{a_and_b_coeff_num}. 

\paragraph{Behavior of the coefficient $c^I_{N_I,N_I}$.}
As mentioned before in \eqref{c_weighted_sum}, the coefficient $c^I_{N_I,N_I}$ can be expressed 
as a weighted sum of the function values: $c^I_{N_I, N_I} = \sum_{k= -N_I}^{N_I} \omega_{k,N_I,N_I} f_{i+k}$.
It is more interesting to study the behavior of these weights than the behavior of the sum.

\begin{figure}[h]
    \centering
    \includegraphics[width=\textwidth]{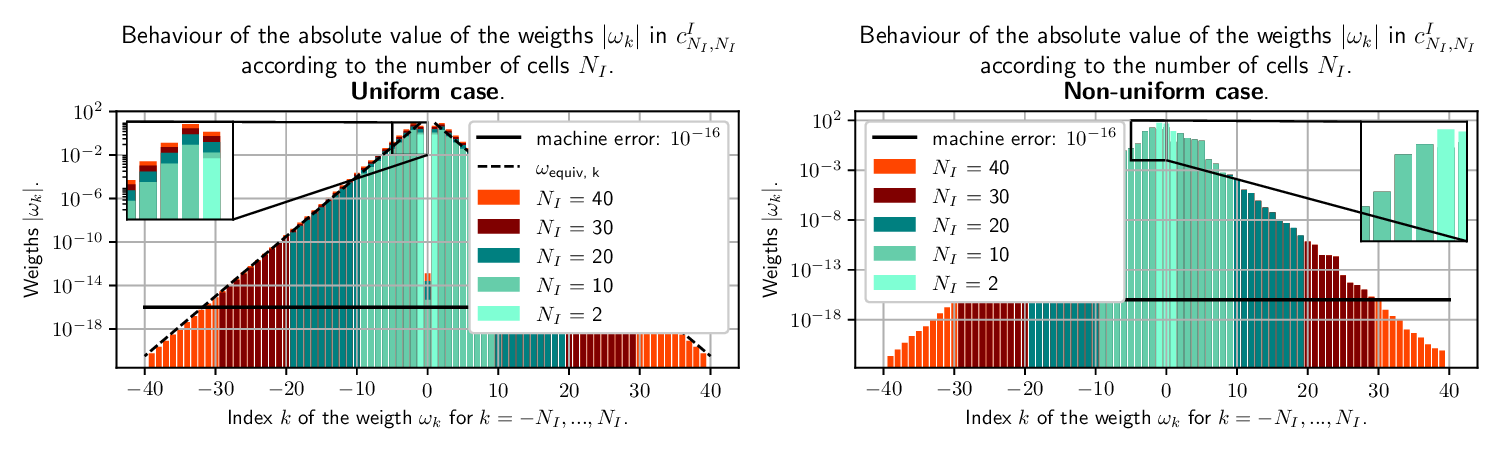}
    \caption{\label{c_coeff_num} 
    Absolute values of the weights $\omega_k$ in \eqref{c_weighted_sum} for the same number 
    of cells $N_I$ in the two patches surrounding the interface $I$ with $N_I$
    varying between 1, 10, 20, 30 and 40. 
    The weights for different $N_I$ are overlaid in the figures. 
    The values of weights close to the interface also change with $N_I$.
    On the left, it is for uniform meshes. On the right, an example with
    non-uniform meshes generated with random functions. 
    We can observe that for $N_I \simeq 30$, the values of these 
    weights reach machine precision. The dependency on more distant function values becomes negligible. 
    On the left, the weights for $N_I = 40$ are compared to 
    $\omega_{k, N_I, N_I} \simeq \omega_{\text{equiv}, k, N_I, N_I} 
    = 3 N_I \left[(2-\sqrt{3})^{|k|} + (2-\sqrt{3})^{2N_I - |k|}\right]$ ,
    following the equivalence relation \eqref{equivalence_w}.
    }
\end{figure}

The absolute values of the weights for an increasing number of interpolation points 
are shown in the Fig. \ref{c_coeff_num}. The support of one patch is of length 1.
 
In the uniform per patch case, the expressions in \eqref{explicit_weights} give us the following equivalences, 
for $n=N_I$ tending to infinity.
\begin{equation}
\left\{
\begin{aligned}
    & \omega_{k,n,n} 
        \simeq
        (-1)^{k} u_1^2\frac{a_{1,1}}{\Delta x^R} \left((2-\sqrt{3})^{k} + (2-\sqrt{3})^{2n-k}\right)
        ,\quad &k = 1, ..., n, \\
    & \omega_{k,n,n} 
        \simeq
        u_{1}(1 + u_{1})
        \left(\frac{a^I_{1,1}}{\Delta x^R} - \frac{b^I_{1,1}}{\Delta x^L}\right) 
        ,\quad &k = 0, \\
    & \omega_{k,n,n} 
        \simeq
        -(-1)^{k} u_1^2\frac{b_{1,1}}{\Delta x^L} \left((2-\sqrt{3})^{-k} + (2-\sqrt{3})^{2n+k}\right)
        ,\quad &k = -n, ..., -1, \\
\end{aligned}
\right. 
\label{equivalence_w}
\end{equation}
These results agree with the tests shown in Fig. \ref{c_coeff_num}. 

\paragraph{Consequences in the matrix system.}
If we now look at the matrix system \eqref{Matrix_Hermite}, we can see that a lot of the terms are negligible
if the number of cells becomes large enough. 

\begin{equation*}
    \begin{aligned}
        \mathbf{S} = 
        \begin{bmatrix}
            s'(\mathcal{X}_1) \\
            \vdots \\
            s'(\mathcal{X}_{I}) \\
            \vdots \\
            s'(\mathcal{X}_{N_p -1})
        \end{bmatrix}
        = \begin{bmatrix}
            1 & \xcancel{-b^{1}_{N^L_1,N^R_1}} & & \\
            \xcancel{-a^{2}_{N^L_2,N^R_2}} & 1& \xcancel{-b^{2}_{N^L_2,N^R_2}} & \\
            & \ddots& \ddots& \ddots\\
            & & \xcancel{-a^{N_p-1}_{N^L_{N_p-1},N^R_{N_p-1}}}& 1\\
        \end{bmatrix}^{-1}
        \begin{bmatrix}
            c^{1}_{N^L_1,N^R_1} 
                + \xcancel{a^{1}_{N^L_1,N^R_1} s'(\mathcal{X}_{0})} \\
            \vdots \\
            c^{I}_{N^L_I,N^R_I} \\
            \vdots \\
            c^{N_p-1}_{N^L_{N_p-1},N^R_{N_p-1}} 
                + \xcancel{b^{N_p-1}_{N^L_{N_p-1},N^R_{N_p-1}} s'(\mathcal{X}_{N_p})} \\
        \end{bmatrix}
    \end{aligned}
\end{equation*}

\begin{equation}
    \begin{aligned}
        \mathbf{S} \rightarrow \mathbf{C_{trunc}}, \qquad \text{ when } N^L_I,N^R_I \rightarrow \infty,     
            \forall I
    \end{aligned}
    \label{Matrix_Hermite_approx}
\end{equation}
with $\mathbf{C_{trunc}}$ the vector of sums on truncated grids,
\begin{equation*}
    \begin{aligned}
        s'(\mathcal{X}_I) \rightarrow 
            \sum_{k= -N^L_{trunc}}^{N^R_{trunc}} \omega_{k, N^L_{trunc}, N^R_{trunc}} f_{i+k}.
    \end{aligned}
\end{equation*}

We observe the same behavior for the different boundary conditions mentioned earlier.
We can choose different truncation values depending on the desired precision. 
We refer to Table \ref{table_equivalence_diff_n} to select the precision. 
\begin{table}[h]
    \centering
    \begin{tabular}{|c|c|c|c|c|c|c|}
        \hline 
        $N_{I, trunc}$ & $5 $
                & $10 $
                & $15 $
                & $20 $
                & $25 $
                & $30 $
                \\
        \hline
        $(2-\sqrt{3})^{N_I}$ 
            & $1.38 \cdot 10^{-3}$ 
            & $1.91 \cdot 10^{-6}$ 
            & $2.64 \cdot 10^{-9}$ 
            & $3.64 \cdot 10^{-12}$ 
            & $5.03 \cdot 10^{-15}$ 
            & $6.94 \cdot 10^{-18}$ 
            \\
        $|a^I_{N_I, N_I}|/|a^I_{1, 1}|$ 
            & $5.52 \cdot 10^{-3}$ 
            & $7.63 \cdot 10^{-6}$ 
            & $1.05 \cdot 10^{-8}$ 
            & $1.46 \cdot 10^{-11}$ 
            & $2.01 \cdot 10^{-14}$ 
            & $2.78 \cdot 10^{-17}$ 
            \\
        $|b^I_{N_I, N_I}|/|b^I_{1, 1}|$ 
            & $5.52 \cdot 10^{-3}$ 
            & $7.63 \cdot 10^{-6}$ 
            & $1.05 \cdot 10^{-8}$ 
            & $1.46 \cdot 10^{-11}$ 
            & $2.01 \cdot 10^{-14}$ 
            & $2.78 \cdot 10^{-17}$ 
            \\
        $|\omega_{N_I, N_I, N_I}||\Delta x^R|/|a^I_{1, 1}|$ 
            & $1.66 \cdot 10^{-2}$ 
            & $2.29 \cdot 10^{-5}$ 
            & $3.16 \cdot 10^{-8}$ 
            & $4.37 \cdot 10^{-11}$ 
            & $6.03 \cdot 10^{-14}$ 
            & $8.33 \cdot 10^{-17}$ 
            \\
        \hline
    \end{tabular}
    \caption{
        \label{table_equivalence_diff_n}
        Different precisions according to the number of cells we take into account. 
        The values in this table are expressed to three significant figures. 
        If we use the approximation, it is equivalent to using the exact formula on subgrids of $N_I$
        cells of the patches and removing the dependency on the other derivatives. 
        We will have an error of order  $|a^I_{N_I, N_I}||s'(i+N_I)| + |b^I_{N_I, N_I}||s'(i-N_I)|$.
    }
\end{table}

\subsection{Semi-Lagrangian advection on a multi-patch domain}
\label{subsection_BSL_multi_patch_1D}
\paragraph{Application of a multi-patch interpolation in BSL.}
Now that we know how to interpolate on a multi-patch domain, let's detail the backward semi-Lagrangian
algorithm on a multi-patch domain. 

\begin{enumerate}
    \item We first solve the equations of the characteristics to get the feet of the characteristics. 
    The feet may not land on the same patch as the starting point (i.e. their initial condition in 
    the equation of the characteristics).

    \item We then identify on which patch the feet are physically located.
    Different algorithms exist. If the mappings $\mathcal{F}$ from each of the logical domains
    to the physical domain are analytically invertible, we can identify on which subdomain a physical 
    coordinate can be found.
    For special geometries such as the concentric rings implemented in our simulations, we can also 
    identify on which patch a coordinate is physically located by computing the distance 
    to the O-point.

    \item The function values have been computed at the previous time step. We can now compute the 
    derivatives with the method described in the previous section, \eqref{Matrix_Hermite} from the 
    function values. 

    \item The function values and the derivatives at the interfaces are then known. We have all the necessary 
    information to build a local spline representation of the function at $t-\Delta t$ on each patch. 
    We interpolate the corresponding local spline at the foot of each characteristic. 
    These values are the values of the function at $t$ at the interpolation points. 
\end{enumerate}

\section{Semi-Lagrangian scheme on a 2D multi-patch domain}
\label{section_Semi-Lagrangian_2D_multi-patch}
\paragraph{}
Now, we know how to handle 1D multi-patch domains. We may also want to handle 
2D multi-patch domains. 
Indeed, in Fig. \ref{JET_SOLEDGE} we have a 2D multi-patch case. 

In this case, the global domain is now split into multiple patches along two dimensions.
In the guiding-center model \eqref{guiding_center_eq}, it is the $r$ direction and the $\theta$ direction. 
Local refinement in a 2D multi-patch domain can be split into two cases: a conforming global 
mesh or a non-conforming global mesh. A mesh is called non-conforming when it is not the result
of a tensor product of two 1D meshes. In multi-patch geometries, this means that there are break points 
defined on one patch with no counterpart in its neighboring patch. 
Moreover, the 2D multi-patch case leads to another case which must be studied: T-joints. 
T-joints designate patch connections involving more than two patches such that only part of an edge 
of a patch is connected to the edge of another patch.
In this section, we present the three cases. 

\subsection{Conforming mesh: a direct application of the 1D case}
\paragraph{}
In the conforming case, the 2D B-spline basis of an equivalent global domain is a 
tensor product of two 1D B-spline bases of equivalent global 1D domains. The computation of the
interface derivatives is then a direct application of the 1D formulae, \eqref{relation_interfaces}, 
\eqref{Matrix_Hermite}. 
In the 2D multi-patch case, we apply Hermite boundary conditions in both directions. So, we need 
to define the derivatives along $r$ and $\theta$ on the edges and the cross-derivatives on the corners. 
Let's note $N_r$ the number of cells in the $r$-direction 
and $N_\theta$ the number of cells in the $\theta$-direction.
\begin{equation*}
    \left\{
    \begin{aligned}
        & \partial_r s(r_i, \theta_j) 
            = \sum_{k=0}^{N_r+2}\sum_{l=0}^{N_\theta+2} c_k c_l b'_k(r_i) b_l(\theta_j), 
            & \qquad i = 0, N_r, \quad j = 0, ..., N_\theta, \\
        & \partial_\theta s(r_i, \theta_j) 
            = \sum_{k=0}^{N_r+2}\sum_{l=0}^{N_\theta+2} c_k c_l b_k(r_i) b'_l(\theta_j), 
            & \qquad i = 0, ..., N_r, \quad j = 0, N_\theta, \\
        & \partial_{r\theta} s(r_i, \theta_j) 
            = \sum_{k=0}^{N_r+2}\sum_{l=0}^{N_\theta+2} c_k c_l b'_k(r_i) b'_l(\theta_j),
            & \qquad i = 0, N_r, \quad j = 0, N_\theta. \\
    \end{aligned}
    \right.
\end{equation*}

For a fixed $j\in \llbracket 0 ,  N_\theta\rrbracket $, $\forall i \in \llbracket 0 ,  N_r \rrbracket$, 
we can write the 2D-spline $(r,\theta)\mapsto s(r,\theta)$ as a 1D-spline $r\mapsto s_{\theta_j}(r)$ 
and its derivatives as, 
\begin{equation*}
    \begin{aligned}
        \partial_r s(r_i, \theta_j) 
            = \left(\sum_{k=0}^{N_r+2} c_k b'_k(r_i)\right) 
            \left(\sum_{l=0}^{N_\theta+2} c_l b_l(\theta_j)\right)
            = s_{\theta_j}'(r_i)
            = \gamma_{\theta_j, i} + \alpha_i s_{\theta_j}'(r_{i+1}) + \beta_i s_{\theta_j}'(r_{i-1}), 
    \end{aligned}
\end{equation*}
by applying the relation \eqref{relation_3_interpol_points} established in the previous section 
to the function $r\mapsto s_{\theta_j}(r)$. The subscript ${\theta_j, i}$ is used to highlight that 
here $\gamma_{\theta_j, i}$ is a weighted sum of $s_{\theta_j}(r_{i-1})$, $s_{\theta_j}(r_{i})$ and 
$s_{\theta_j}(r_{i+1})$.
Similarly, with the other direction, we then have the 
following relation, 
\begin{equation*}
    \begin{aligned}
        \partial_\theta s(r_i, \theta_j)
            = \gamma_{r_i, j} + \alpha_j s_{r_i}'(\theta_{j+1}) + \beta_j s_{r_i}'(\theta_{j-1}).
    \end{aligned}
\end{equation*}

By fixing $\theta_j$ or $r_i$, we can proceed similarly for the cross-derivatives. 
Let's take for example $\theta_j$ fixed:
\begin{equation*}
    \begin{aligned}
        \partial_{r\theta} s(r_i, \theta_j) 
            = \left(\sum_{k=0}^{N_r+2} c_k b'_k(r_i)\right) 
            \left(\sum_{l=0}^{N_\theta+2} c_l b'_l(\theta_j)\right)
            = s_{\theta_j'}'(r_i)
            = \gamma_{\theta_j', i} + \alpha_i s_{\theta_j'}'(r_{i+1}) + \beta_i s_{\theta_j'}'(r_{i-1}), 
    \end{aligned}
\end{equation*}
with $s_{\theta_j'}: r \mapsto \partial_\theta s(r, \theta_j)$.

We then have equivalent formulae for the derivatives along $r$ and $\theta$ and the cross-derivatives. 
With the method applied in the previous section \ref{section_Semi-Lagrangian_1D_multi-patch}, 
we can obtain equations linking the derivatives at the interfaces. 
Let's denote by $\{R_I\}_I$ the interfaces in the $r$-direction and $\{\Theta_J\}_J$ the interfaces 
in the $\theta$-direction. We note $N_p^R$ the number of patches in the $r$ direction and 
$N_p^\Theta$ the number of patches in the $\theta$ direction. 
The equations linking the derivatives are: 
\begin{equation}
    \left\{
    \begin{aligned}
        & \partial_r s(R_I, \theta_j) 
            = c^{I}_{N^L_I, N^R_I} 
            + a^{I}_{N^L_I, N^R_I} \partial_r s(R_{I+1}, \theta_j) 
            + b^{I}_{N^L_I, N^R_I} \partial_r s(R_{I-1}, \theta_j) , 
            & \forall j\in \llbracket 0 ,  N_\theta\rrbracket, 
            I \in \llbracket 1 ,  N_p^R \rrbracket, \\
        & \partial_\theta s(r_i, \Theta_J) 
            = c^{J}_{N^L_J, N^R_J} 
            + a^{J}_{N^L_J, N^R_J} \partial_\theta s(r_i, \Theta_{J+1}) 
            + b^{J}_{N^L_J, N^R_J} \partial_\theta s(r_i, \Theta_{J-1}) , 
            &  \forall i\in \llbracket 0 ,  N_r \rrbracket, 
            J \in \llbracket 1 ,  N_p^\Theta \rrbracket, \\
        & \partial_{r\theta} s(R_I, \Theta_J) 
            = c^{I,*}_{N^L_I, N^R_I} 
            + a^{I}_{N^L_I, N^R_I} \partial_{r\theta} s(R_{I+1}, \Theta_J) 
            + b^{I}_{N^L_I, N^R_I} \partial_{r\theta} s(R_{I-1}, \Theta_J) , 
            & I \in \llbracket 1 ,  N_p^R \rrbracket, 
            J \in \llbracket 1 ,  N_p^\Theta \rrbracket, \\
    \end{aligned}
    \right.
    \label{relations_derivatives_2D}
\end{equation}
with the coefficients $a^{I}_{N^L_I, N^R_I}$, $b^{I}_{N^L_I, N^R_I}$ and $c^{I}_{N^L_I, N^R_I}$
calculated with the same formula as in \eqref{relation_interfaces}
for $N^R_I$ cells along the $r$-direction in the right patch and 
$N^L_I$ cells along the $r$-direction in the left patch. 
We proceed similarly for the $a^{J}_{N^L_J, N^R_J}$, $b^{J}_{N^L_J, N^R_J}$ and $c^{J}_{N^L_J, N^R_J}$
coefficients in the $\theta$-direction. 
The $c^{I}_{N^L_I, N^R_I}$ coefficient is a weighted sum of the function values 
$\{f_{i+k,j}\}_{k=-N^L_I,...,N^R_I}$ in the $r$-direction
for a fixed $j$ and $i$ corresponding to the interpolation 
point located on the interface $R_I$. 
In the equation describing the cross-derivatives, $c^{I,*}_{N^L_I, N^R_I}$ is a weighted sum 
of the function derivatives $\{\partial_\theta f_{i+k,j}\}_{k=-N^L_I,...,N^R_I}$ along the 
$\theta$-direction in the $r$-direction for a fixed $j$ and $i$ corresponding to the interpolation 
point located on the interface $R_I$. 

\paragraph{Remark.}
A remark concerning the given relation between the cross-derivatives is that we can also establish
it in the other direction, 
\begin{equation*}
    \partial_{r\theta} s(R_I, \Theta_J) 
            = c^{J,*}_{N^L_J, N^R_J} 
            + a^{J}_{N^L_J, N^R_J} \partial_{r\theta} s(R_I, \Theta_{J+1}) 
            + b^{J}_{N^L_J, N^R_J} \partial_{r\theta} s(R_I, \Theta_{J-1}) , 
            \qquad  
            I \in \llbracket 1 ,  N_p^R \rrbracket, 
            J \in \llbracket 1 ,  N_p^\Theta \rrbracket,
\end{equation*}
where $c^{J,*}_{n_J, m_J}$ is the weighted sum 
of the function derivatives $\{\partial_r f_{i,j+k}\}_{k=-N^L_J, ..., N^R_J}$ along the
$r$-direction in the $\theta$-direction for a fixed $i$ and $j$ corresponding to the interpolation 
point located on the interface $\Theta_J$. Both formulae are equivalent and the cross-derivatives
are equal to the cross-derivatives of an equivalent 2D global spline.

\paragraph{Algorithm.}
As the equations describing the cross-derivatives depend on the derivatives of the function at the 
interfaces, to compute all the derivatives in the conforming case, the algorithm needs only two steps: 
\begin{enumerate}
    \item From the function values at the interpolation points, compute the derivatives 
    along $r$ and $\theta$ with \eqref{relations_derivatives_2D}. 
    \item From the computed derivatives at the interfaces, compute the cross-derivatives 
    at the corners by choosing one of the two directions with \eqref{relations_derivatives_2D}. 
\end{enumerate}

\subsection{Non-Conforming mesh}
\paragraph{}
Non-conforming multi-patch meshes refer to meshes where break points 
defined on one patch have no counterpart on its neighboring patch. 
The global mesh is not the result of a tensor product of two 1D global meshes.
So, there is no equivalent global spline. 
The mesh local to a given patch remains conforming and the result of a tensor product of two 1D local meshes.
So we can still define local splines. 
We suppose here that the non-conforming mesh is not too degenerated. 
By this, we mean that for an interface between two patches, 
the one-dimensional mesh of one of the two patches is included in the one-dimensional mesh of the other patch. 

We cannot apply formulae 
\eqref{relations_derivatives_2D} along all the lines of each direction because they may end at 
some interfaces. The idea is then to progressively compute the derivatives and the cross-derivatives. 
As mentioned before, the cross-derivatives can be computed from the derivatives along one direction. 
From the derivatives and the cross-derivatives, we can build spline representations of the derivatives. 
If we want to impose a $\mathcal{C}^1$ regularity between the patches as mentioned previously, we 
build a spline on a coarser patch of an interface. We interpolate it at the interpolation points 
on the refined patch which have 
no counterpart on the coarser patch. The interpolated values are then assigned to 
the refined patch to enforce the $\mathcal{C}^1$ regularity.

In some cases, we can define an algorithm to compute all the derivatives and cross-derivatives. 
For other cases, we may be limited and will not use all the information of the patches. 
For these cases, it is recommended to apply an approximation of the exact formulae \eqref{Matrix_Hermite_approx},
if the number of cells per patch is big enough.

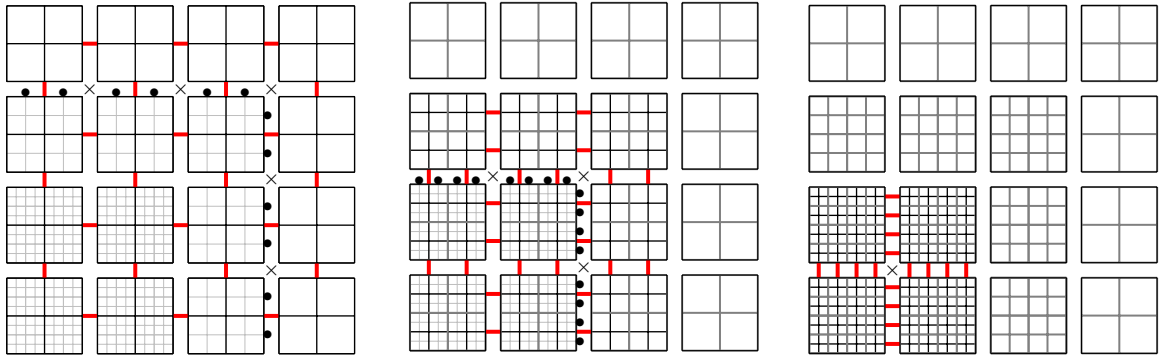
\begin{figure}[H]
    \centering
    \begin{subfigure}[t]{0.3\textwidth}
        \centering
        \begin{tikzpicture}[xscale = 1, yscale = 1]
            \foreach \x in {0, 1}{
                \foreach \y in {0, 1}{
                    \begin{scope}[xshift=\x*1.2cm, yshift=\y*1.2cm]
                        \draw[lightgray, line width= 0.25] (0,0) grid [step = 0.125] (1,1);
                    \end{scope}
                }
            }
            \foreach \x in {0, 1, 2}{
                \foreach \y in {0, 1, 2}{
                    \begin{scope}[xshift=\x*1.2cm, yshift=\y*1.2cm]
                        \draw[lightgray, line width = 0.25] (0,0) grid [step = 0.25] (1,1);
                    \end{scope}
                }
            }
            \foreach \x in {0, 1, 2, 3}{
                \foreach \y in {0, 1, 2, 3}{
                    \begin{scope}[xshift=\x*1.2cm, yshift=\y*1.2cm]
                        \draw[black, line width = 0.5] (0,0) grid [step = 0.5] (1,1);
                    \end{scope}
                }
            }    
                 
            \foreach \x in {0, 1, 2}{
                \foreach \y in {0, 1, 2, 3}{
                    \begin{scope}[xshift=\x*1.2cm, yshift=\y*1.2cm]
                        \draw[red, line width = 1.5] (1, 0.5) -- (1.2, 0.5);
                    \end{scope}
                }
            }            
            \foreach \x in {0, 1, 2, 3}{
                \foreach \y in {0, 1, 2}{
                    \begin{scope}[xshift=\x*1.2cm, yshift=\y*1.2cm]
                        \draw[red, line width = 1.5] (0.5, 1) -- (0.5, 1.2);
                    \end{scope}
                }
            }   
                
            \foreach \x in {0, 1, 2, 3}{
                \foreach \y in {0, 1, 2, 3}{
                    \begin{scope}[xshift=\x*1.2cm, yshift=\y*1.2cm]
                        \draw[black] (0,0) grid[step = 1] (1,1);
                    \end{scope}
                }
            }

            \foreach \x in {1, 2, 3}{
                \begin{scope}[xshift=\x*1.2cm, yshift= 3*1.2cm]
                    \node at (-0.1, -0.1) {\scriptsize$\times$}; 
                \end{scope}
            }
            \foreach \y in {1, 2}{
                \begin{scope}[xshift=3*1.2cm, yshift= \y*1.2cm]
                    \node at (-0.1, -0.1) {\scriptsize$\times$}; 
                \end{scope}
            }

            \foreach \x in {0, 1, 2}{
                \begin{scope}[xshift=\x*1.2cm, yshift= 3*1.2cm]
                    \node at (0.25, -0.15) {\scriptsize$\bullet$}; 
                    \node at (0.75, -0.15) {\scriptsize$\bullet$}; 
                \end{scope}
            }
            \foreach \y in {0, 1, 2}{
                \begin{scope}[xshift=3*1.2cm, yshift= \y*1.2cm]
                    \node at (-0.15, 0.25) {\scriptsize$\bullet$}; 
                    \node at (-0.15, 0.75) {\scriptsize$\bullet$}; 
                \end{scope}
            }
                
        \end{tikzpicture}
        \subcaption{
            Compute the derivatives (red lines) along the lines crossing 
            the whole domain in $r$ and $\theta$ directions, the "conforming" lines 
            with \eqref{relations_derivatives_2D}.
            From these derivatives, we compute the cross-derivatives ($\times$) 
            on the coarser patches. With the derivatives along $r$ and the cross-derivatives, 
            we can build spline representations of the $r$-derivatives on the coarse patches. 
            We do the same for the $\theta$-derivatives.
            We interpolate the local splines at the interpolation points of the refined patches with no 
            counterpart in the coarse patches ($\bullet$). 
        }
    \end{subfigure}
    \hfill
    \begin{subfigure}[t]{0.3\textwidth}
        \centering
        \begin{tikzpicture}[xscale = 1, yscale = 1]
            \foreach \x in {0, 1}{
                \foreach \y in {0, 1}{
                    \begin{scope}[xshift=\x*1.2cm, yshift=\y*1.2cm]
                        \draw[lightgray, line width= 0.25] (0,0) grid [step = 0.125] (1,1);
                    \end{scope}
                }
            }
            \foreach \x in {0, 1, 2}{
                \foreach \y in {0, 1, 2}{
                    \begin{scope}[xshift=\x*1.2cm, yshift=\y*1.2cm]
                        \draw[black, line width = 0.5] (0,0) grid [step = 0.25] (1,1);
                    \end{scope}
                }
            }
            \foreach \x in {0, 1, 2 ,3}{
                \foreach \y in {0, 1, 2, 3}{
                    \begin{scope}[xshift=\x*1.2cm, yshift=\y*1.2cm]
                        \draw[gray, line width = 0.75] (0,0) grid [step = 0.5] (1,1);
                    \end{scope}
                }
            }

            \foreach \x in {0, 1}{
                \foreach \y in {0, 1, 2}{
                    \begin{scope}[xshift=\x*1.2cm, yshift=\y*1.2cm]
                        \draw[red, line width = 1.5] (1, 0.25) -- (1.2, 0.25);
                        \draw[red, line width = 1.5] (1, 0.75) -- (1.2, 0.75);
                    \end{scope}
                }
            }            
            \foreach \x in {0, 1, 2}{
                \foreach \y in {0, 1}{
                    \begin{scope}[xshift=\x*1.2cm, yshift=\y*1.2cm]
                        \draw[red, line width = 1.5] (0.25, 1) -- (0.25, 1.2);
                        \draw[red, line width = 1.5] (0.75, 1) -- (0.75, 1.2);
                    \end{scope}
                }
            }

            \foreach \x in {0, 1, 2 ,3}{
                \foreach \y in {0, 1, 2, 3}{
                    \begin{scope}[xshift=\x*1.2cm, yshift=\y*1.2cm]
                        \draw[black] (0,0) grid [step = 1] (1,1);
                    \end{scope}
                }
            }

            \foreach \x in {1, 2}{
                \begin{scope}[xshift=\x*1.2cm, yshift= 2*1.2cm]
                    \node at (-0.1, -0.1) {\scriptsize$\times$}; 
                \end{scope}
            }
            \foreach \y in {1}{
                \begin{scope}[xshift=2*1.2cm, yshift= \y*1.2cm]
                    \node at (-0.1, -0.1) {\scriptsize$\times$}; 
                \end{scope}
            }

            \foreach \x in {0, 1}{
                \begin{scope}[xshift=\x*1.2cm, yshift= 2*1.2cm]
                    \node at (0.125, -0.15) {\scriptsize$\bullet$}; 
                    \node at (0.375, -0.15) {\scriptsize$\bullet$}; 
                    \node at (0.625, -0.15) {\scriptsize$\bullet$}; 
                    \node at (0.875, -0.15) {\scriptsize$\bullet$}; 
                \end{scope}
            }
            \foreach \y in {0, 1}{
                \begin{scope}[xshift=2*1.2cm, yshift= \y*1.2cm]
                    \node at (-0.15, 0.125) {\scriptsize$\bullet$}; 
                    \node at (-0.15, 0.375) {\scriptsize$\bullet$}; 
                    \node at (-0.15, 0.625) {\scriptsize$\bullet$}; 
                    \node at (-0.15, 0.875) {\scriptsize$\bullet$}; 
                \end{scope}
            }

        \end{tikzpicture}
        \subcaption{      
            With the derivatives computed at the previous step at the "non-conforming" 
            interpolation points ($\bullet$) and the boundary conditions, 
            we can compute the derivatives along the 
            red lines with \eqref{relations_derivatives_2D}. 
            From these derivatives, we compute the cross-derivatives at $\times$.
            As previously, we build a spline representation of the derivatives 
            along $r$ and $\theta$ using the newly computed derivatives and 
            cross-derivatives. The splines have to be interpolated at the $\bullet$
            points. 
        }
    \end{subfigure}
    \hfill
    \begin{subfigure}[t]{0.3\textwidth}
        \centering
        \begin{tikzpicture}[xscale = 1, yscale = 1]
            \foreach \x in {0, 1}{
                \foreach \y in {0, 1}{
                    \begin{scope}[xshift=\x*1.2cm, yshift=\y*1.2cm]
                        \draw[black, line width= 0.5] (0,0) grid [step = 0.125] (1,1);
                    \end{scope}
                }
            }
            \foreach \x in {0, 1, 2}{
                \foreach \y in {0, 1, 2}{
                    \begin{scope}[xshift=\x*1.2cm, yshift=\y*1.2cm]
                        \draw[gray, line width = 0.75] (0,0) grid [step = 0.25] (1,1);
                    \end{scope}
                }
            }
            \foreach \x in {0, 1, 2 ,3}{
                \foreach \y in {0, 1, 2, 3}{
                    \begin{scope}[xshift=\x*1.2cm, yshift=\y*1.2cm]
                        \draw[gray, line width = 0.75] (0,0) grid [step = 0.5] (1,1);
                    \end{scope}
                }
            }
            
            \foreach \x in {0}{
                \foreach \y in {0, 1}{
                    \begin{scope}[xshift=\x*1.2cm, yshift=\y*1.2cm]
                        \draw[red, line width = 1.5] (1, 0.125) -- (1.2, 0.125);
                        \draw[red, line width = 1.5] (1, 0.375) -- (1.2, 0.375);
                        \draw[red, line width = 1.5] (1, 0.625) -- (1.2, 0.625);
                        \draw[red, line width = 1.5] (1, 0.875) -- (1.2, 0.875);
                    \end{scope}
                }
            }      
            \foreach \x in {0, 1}{
                \foreach \y in {0}{
                    \begin{scope}[xshift=\x*1.2cm, yshift=\y*1.2cm]
                        \draw[red, line width = 1.5] (0.125, 1) -- (0.125, 1.2);
                        \draw[red, line width = 1.5] (0.375, 1) -- (0.375, 1.2);
                        \draw[red, line width = 1.5] (0.625, 1) -- (0.625, 1.2);
                        \draw[red, line width = 1.5] (0.875, 1) -- (0.875, 1.2);
                    \end{scope}
                }
            }   

            \foreach \x in {0, 1, 2 ,3}{
                \foreach \y in {0, 1, 2, 3}{
                    \begin{scope}[xshift=\x*1.2cm, yshift=\y*1.2cm]
                        \draw[black] (0,0) grid [step = 1] (1,1);
                    \end{scope}
                }
            }

            \foreach \x in {1}{
                \begin{scope}[xshift=\x*1.2cm, yshift= 1*1.2cm]
                    \node at (-0.1, -0.1) {\scriptsize$\times$}; 
                \end{scope}
            }  

        \end{tikzpicture}
        \subcaption{      
            With the derivatives computed at the previous step and the boundary conditions, 
            we can compute the derivatives along the 
            red lines using \eqref{relations_derivatives_2D} again. 
            From them, we compute the cross-derivatives at $\times$.
            The function values are known at every interpolation point and all the 
            derivatives are now determined. We have all the necessary information to 
            build local spline representations of the function on the whole domain. 
        }
    \end{subfigure}
    \hfill
    \vspace*{-0.2cm}

    \caption{
        \label{example_1_non_conforming}
        Example of an algorithm for non-conforming patches where we can successively compute all 
        the derivatives and cross-derivatives. 
    }
\end{figure}

To illustrate this point, here are two examples of non-conforming patches. 
Fig. \ref{example_1_non_conforming} shows an example of a situation where we can compute all the derivatives and 
cross-derivatives using all the available information.
Fig. \ref{example_2_non_conforming} shows an example of a situation where we cannot use all 
the available information if we want to apply the exact formulae \eqref{relations_derivatives_2D}.

\begin{figure}[H]
    \centering
    \begin{subfigure}[t]{0.3\textwidth}
        \centering
        \begin{tikzpicture}[xscale = 1, yscale = 1]
            \begin{scope}[xshift=2*1.2cm, yshift=3*1.2cm]
                \draw[lightgray, line width=0.25] (0,0) grid [step = 0.125] (1,1);
            \end{scope}
            \begin{scope}[xshift=3*1.2cm, yshift=2*1.2cm]
                \draw[lightgray, line width= 0.25] (0,0) grid [step = 0.125] (1,1);
            \end{scope}
            \begin{scope}[xshift=1*1.2cm, yshift=0*1.2cm]
                \draw[lightgray, line width= 0.25] (0,0) grid [step = 0.125] (1,1);
            \end{scope}
            \begin{scope}[xshift=0*1.2cm, yshift=1*1.2cm]
                \draw[lightgray, line width= 0.25] (0,0) grid [step = 0.125] (1,1);
            \end{scope}
    
            \begin{scope}[xshift=2*1.2cm, yshift=3*1.2cm]
                \draw[lightgray, line width= 0.25] (0,0) grid [step = 0.25] (1,1);
            \end{scope}
            \begin{scope}[xshift=3*1.2cm, yshift=2*1.2cm]
                \draw[lightgray, line width= 0.25] (0,0) grid [step = 0.25] (1,1);
            \end{scope}
            \begin{scope}[xshift=1*1.2cm, yshift=0*1.2cm]
                \draw[lightgray, line width= 0.25] (0,0) grid [step = 0.25] (1,1);
            \end{scope}
            \begin{scope}[xshift=0*1.2cm, yshift=1*1.2cm]
                \draw[lightgray, line width= 0.25] (0,0) grid [step = 0.25] (1,1);
            \end{scope}
    
            \foreach \x in {1, 2}{
                \foreach \y in {1, 2}{
                    \begin{scope}[xshift=\x*1.2cm, yshift=\y*1.2cm]
                        \draw[lightgray, line width= 0.25] (0,0) grid [step = 0.25] (1,1);
                    \end{scope}
                }
            }
            \foreach \x in {0, 1, 2 ,3}{
                \foreach \y in {0, 1, 2, 3}{
                    \begin{scope}[xshift=\x*1.2cm, yshift=\y*1.2cm]
                        \draw[black, line width= 0.5] (0,0) grid [step = 0.5] (1,1);
                    \end{scope}
                }
            }
            \foreach \x in {0, 1, 2}{
                \foreach \y in {0, 1, 2, 3}{
                    \begin{scope}[xshift=\x*1.2cm, yshift=\y*1.2cm]
                        \draw[red, line width = 1.5] (1, 0.5) -- (1.2, 0.5);
                    \end{scope}
                }
            }            
            \foreach \x in {0, 1, 2, 3}{
                \foreach \y in {0, 1, 2}{
                    \begin{scope}[xshift=\x*1.2cm, yshift=\y*1.2cm]
                        \draw[red, line width = 1.5] (0.5, 1) -- (0.5, 1.2);
                    \end{scope}
                }
            }   
                
            \foreach \x in {0, 1, 2 ,3}{
                \foreach \y in {0, 1, 2, 3}{
                    \begin{scope}[xshift=\x*1.2cm, yshift=\y*1.2cm]
                        \draw[black] (0,0) grid [step = 1] (1,1);
                    \end{scope}
                }
            }

            \foreach \x in {1, 2, 3}{
                \begin{scope}[xshift=\x*1.2cm, yshift= 3*1.2cm]
                    \node at (-0.1, -0.1) {\scriptsize$\times$}; 
                \end{scope}
            }

            \begin{scope}[xshift=3*1.2cm, yshift=3*1.2cm]
                \node[fill=white] at (0.5,0.5) {16}; 
            \end{scope}
            \begin{scope}[xshift=2*1.2cm, yshift=3*1.2cm]
                \node[fill=white] at (0.5,0.5) {15}; 
            \end{scope}
            \begin{scope}[xshift=2*1.2cm, yshift=2*1.2cm]
                \node[fill=white] at (0.5,0.5) {11}; 
            \end{scope}
                
        \end{tikzpicture}
        \subcaption{
            As in the first example in Fig. \ref{example_1_non_conforming}, we start by computing
            the derivatives along the "conforming" lines in red on the figure. 
            Then, we can compute the cross-derivatives at the $\times$.
            However, by doing that we do not use the fact that patches 15 and 11 are refined. 
            If we build a spline representation of the derivative on the interface 15/11, 
            it will be as coarse as if we build it in the patch 16 (for example). 
            So, we lose information from the patches 15 and 11. 
        }
    \end{subfigure}
    \hfill
    \begin{subfigure}[t]{0.3\textwidth}
        \centering
        \begin{tikzpicture}[xscale = 1, yscale = 1]
            \begin{scope}[xshift=2*1.2cm, yshift=3*1.2cm]
                \draw[lightgray, line width=0.25] (0,0) grid [step = 0.125] (1,1);
            \end{scope}
            \begin{scope}[xshift=3*1.2cm, yshift=2*1.2cm]
                \draw[lightgray, line width= 0.25] (0,0) grid [step = 0.125] (1,1);
            \end{scope}
            \begin{scope}[xshift=1*1.2cm, yshift=0*1.2cm]
                \draw[lightgray, line width= 0.25] (0,0) grid [step = 0.125] (1,1);
            \end{scope}
            \begin{scope}[xshift=0*1.2cm, yshift=1*1.2cm]
                \draw[lightgray, line width= 0.25] (0,0) grid [step = 0.125] (1,1);
            \end{scope}

            \begin{scope}[xshift=2*1.2cm, yshift=3*1.2cm]
                \draw[lightgray, line width= 0.25] (0,0) grid [step = 0.25] (1,1);
            \end{scope}
            \begin{scope}[xshift=3*1.2cm, yshift=2*1.2cm]
                \draw[lightgray, line width= 0.25] (0,0) grid [step = 0.25] (1,1);
            \end{scope}
            \begin{scope}[xshift=1*1.2cm, yshift=0*1.2cm]
                \draw[lightgray, line width= 0.25] (0,0) grid [step = 0.25] (1,1);
            \end{scope}
            \begin{scope}[xshift=0*1.2cm, yshift=1*1.2cm]
                \draw[lightgray, line width= 0.25] (0,0) grid [step = 0.25] (1,1);
            \end{scope}

            \foreach \x in {1, 2}{
                \foreach \y in {1, 2}{
                    \begin{scope}[xshift=\x*1.2cm, yshift=\y*1.2cm]
                        \draw[lightgray, line width= 0.25] (0,0) grid [step = 0.25] (1,1);
                    \end{scope}
                }
            }
            \foreach \x in {0, 1, 2 ,3}{
                \foreach \y in {0, 1, 2, 3}{
                    \begin{scope}[xshift=\x*1.2cm, yshift=\y*1.2cm]
                        \draw[black, line width= 0.5] (0,0) grid [step = 0.5] (1,1);
                    \end{scope}
                }
            }
            \foreach \x in {0, 1, 2}{
                \foreach \y in {0, 1, 2, 3}{
                    \begin{scope}[xshift=\x*1.2cm, yshift=\y*1.2cm]
                        \draw[red, line width = 1.5] (1, 0.5) -- (1.2, 0.5);
                    \end{scope}
                }
            }            
            \foreach \x in {0, 1, 2, 3}{
                \foreach \y in {0, 1, 2}{
                    \begin{scope}[xshift=\x*1.2cm, yshift=\y*1.2cm]
                        \draw[red, line width = 1.5] (0.5, 1) -- (0.5, 1.2);
                    \end{scope}
                }
            }   
            
            \foreach \x in {0, 1, 2 ,3}{
                \foreach \y in {0, 1, 2, 3}{
                    \begin{scope}[xshift=\x*1.2cm, yshift=\y*1.2cm]
                        \draw[black] (0,0) grid [step = 1] (1,1);
                    \end{scope}
                }
            }

            \foreach \x in {1, 2, 3}{
                \begin{scope}[xshift=\x*1.2cm, yshift= 1*1.2cm]
                    \node at (-0.1, -0.1) {\scriptsize$\times$}; 
                \end{scope}
            }


            \begin{scope}[xshift=1*1.2cm, yshift=0*1.2cm]
                \node[fill=white] at (0.5,0.5) {2}; 
            \end{scope}
            \begin{scope}[xshift=1*1.2cm, yshift=1*1.2cm]
                \node[fill=white] at (0.5,0.5) {6}; 
            \end{scope}
            \begin{scope}[xshift=2*1.2cm, yshift=0*1.2cm]
                \node[fill=white] at (0.5,0.5) {3}; 
            \end{scope}
            \begin{scope}[xshift=2*1.2cm, yshift=1*1.2cm]
                \node[fill=white] at (0.5,0.5) {7}; 
            \end{scope}
                
        \end{tikzpicture}
        \subcaption{
            If the option in the sub-figure (a) is a bad idea, we can then think that
            we can start by computing the cross-derivatives at the $\times$ 
            on this sub-figure. 
            However, we will encounter the same problem at the interface 2/6. 
            Because of the dependency of the derivatives at all the interfaces 
            in \eqref{relations_derivatives_2D}, it will affect the cross-derivatives
            used to define a spline at the interface 3/7. 
        }
    \end{subfigure}
    \hfill
    \begin{subfigure}[t]{0.3\textwidth}
        \centering
        \begin{tikzpicture}[xscale = 1, yscale = 1]
            \begin{scope}[xshift=2*1.2cm, yshift=3*1.2cm]
                \draw[lightgray, line width=0.25] (0,0) grid [step = 0.125] (1,1);
            \end{scope}
            \begin{scope}[xshift=3*1.2cm, yshift=2*1.2cm]
                \draw[lightgray, line width= 0.25] (0,0) grid [step = 0.125] (1,1);
            \end{scope}
            \begin{scope}[xshift=1*1.2cm, yshift=0*1.2cm]
                \draw[lightgray, line width= 0.25] (0,0) grid [step = 0.125] (1,1);
            \end{scope}
            \begin{scope}[xshift=0*1.2cm, yshift=1*1.2cm]
                \draw[lightgray, line width= 0.25] (0,0) grid [step = 0.125] (1,1);
            \end{scope}

            \begin{scope}[xshift=2*1.2cm, yshift=3*1.2cm]
                \draw[lightgray, line width= 0.25] (0,0) grid [step = 0.25] (1,1);
            \end{scope}
            \begin{scope}[xshift=3*1.2cm, yshift=2*1.2cm]
                \draw[lightgray, line width= 0.25] (0,0) grid [step = 0.25] (1,1);
            \end{scope}
            \begin{scope}[xshift=1*1.2cm, yshift=0*1.2cm]
                \draw[lightgray, line width= 0.25] (0,0) grid [step = 0.25] (1,1);
            \end{scope}
            \begin{scope}[xshift=0*1.2cm, yshift=1*1.2cm]
                \draw[lightgray, line width= 0.25] (0,0) grid [step = 0.25] (1,1);
            \end{scope}

            \foreach \x in {1, 2}{
                \foreach \y in {1, 2}{
                    \begin{scope}[xshift=\x*1.2cm, yshift=\y*1.2cm]
                        \draw[lightgray, line width= 0.25] (0,0) grid [step = 0.25] (1,1);
                    \end{scope}
                }
            }
            \foreach \x in {0, 1, 2 ,3}{
                \foreach \y in {0, 1, 2, 3}{
                    \begin{scope}[xshift=\x*1.2cm, yshift=\y*1.2cm]
                        \draw[black, line width= 0.5] (0,0) grid [step = 0.5] (1,1);
                    \end{scope}
                }
            }
            
            \foreach \x in {0, 1, 2 ,3}{
                \foreach \y in {0, 1, 2, 3}{
                    \begin{scope}[xshift=\x*1.2cm, yshift=\y*1.2cm]
                        \draw[black] (0,0) grid [step = 1] (1,1);
                    \end{scope}
                }
            }

            \foreach \x in {2}{
                \begin{scope}[xshift=\x*1.2cm, yshift= 3*1.2cm]
                    \node at (0.25, -0.15) {\scriptsize$\bullet$}; 
                    \node at (0.75, -0.15) {\scriptsize$\bullet$}; 
                \end{scope}
            }

            \begin{scope}[xshift=1.5*1.2cm, yshift=1.5*1.2cm]
                \node[fill=white] at (0.5,0.5) { \Huge \textbf{?} }; 
            \end{scope}
                
        \end{tikzpicture}
        \subcaption{
            The case of the two first sub-figures can also be encountered if we compute the 
            cross-derivatives in the $\theta$-direction. 
            This illustrates that if we apply the exact formulae in \eqref{relations_derivatives_2D}, 
            we may not always use all the available information. 
            A solution is to add enough cells in the problematic patches in order 
            to use an approximation of the formulae (see \eqref{Matrix_Hermite_approx}). 
            Instead of computing the derivatives at $\bullet$ by interpolating a spline 
            representation of the derivative, we approximate the derivatives with 
            \eqref{Matrix_Hermite_approx}.
        }
    \end{subfigure}
    \hfill
    \vspace*{-0.2cm}

    \caption{
        \label{example_2_non_conforming}
        Example of non-conforming patches where the application of the exact formulae 
        \eqref{relations_derivatives_2D}
        may not use all the available information. In this case, the application of the approximated 
        formulae \eqref{Matrix_Hermite_approx} seems to be a better solution if the number of the 
        cells in the patches is big enough. 
    }
\end{figure}

In practice, the first case of non-conforming meshes may appear more often.
It is for example the case in Fig. \ref{JET_SOLEDGE}.
The second case was described in this article mainly as a warning.

\paragraph{Remark.}
As mentioned previously, we impose a $\mathcal{C}^1$ regularity at the interfaces. 
It is imposed by building spline representations of the derivatives and evaluating 
them at the interpolation points with no counterpart in the neighboring patches.
Thus, contrary to the conforming case, we lose the 
$\mathcal{C}^2$ regularity that we had with the Hermite polynomials in Section 
\ref{subsubsection_Relation_between_derivatives_at_three_consecutive_points}.

\subsection{Geometries with T-joints}
\subsubsection{Definition of a T-joint}
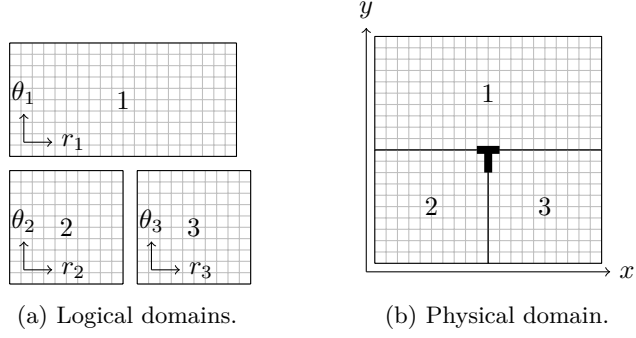
\begin{figure}[h]
    \centering
    \begin{subfigure}[t]{0.3\textwidth}
        \centering
        \begin{tikzpicture}[xscale=0.75, yscale=0.75]
            \draw[lightgray] (0,0) grid[step=0.2] (4,2); 
            \draw (0,0) grid[xstep=4, ystep=2] (4,2); 
            \node at (2,1) {1}; 
            \draw[->] (0.25, 0.25) -- (0.75, 0.25) node[right] {$r_1$}; 
            \draw[->] (0.25, 0.25) -- (0.25, 0.75) node[above] {$\theta_1$}; 

            \begin{scope}[xshift=0cm, yshift=-2.25cm]
                \draw[lightgray] (0,0) grid[step=0.2] (2,2); 
                \draw (0,0) grid[step=2] (2,2); 
                \node at (1,1) {2}; 
                \draw[->] (0.25, 0.25) -- (0.75, 0.25) node[right] {$r_2$}; 
                \draw[->] (0.25, 0.25) -- (0.25, 0.75) node[above] {$\theta_2$}; 
            \end{scope}
            
            \begin{scope}[xshift=2.25cm, yshift=-2.25cm]
                \draw[lightgray] (0,0) grid[step=0.2] (2,2); 
                \draw (0,0) grid[step=2] (2,2); 
                \node at (1,1) {3}; 
                \draw[->] (0.25, 0.25) -- (0.75, 0.25) node[right] {$r_3$}; 
                \draw[->] (0.25, 0.25) -- (0.25, 0.75) node[above] {$\theta_3$}; 
            \end{scope}

        \end{tikzpicture}
        \caption{Logical domains.}
    \end{subfigure}
    \begin{subfigure}[t]{0.3\textwidth}
        \centering
        \begin{tikzpicture}[xscale=0.75, yscale=0.75]
            \draw[->] (-.15, -2.15) -- (4.15, -2.15) node[right] {$x$}; 
            \draw[->] (-.15, -2.15) -- (-.15,  2.15) node[above] {$y$}; 

            \draw[lightgray] (0,0) grid[step=0.2] (4,2); 
            \draw (0,0) grid[xstep=4, ystep=2] (4,2); 
            \node at (2,1) {1}; 

            \begin{scope}[xshift=0cm, yshift=-2cm]
                \draw[lightgray] (0,0) grid[step=0.2] (2,2); 
                \draw (0,0) grid[step=2] (2,2); 
                \node at (1,1) {2}; 
            \end{scope}
            
            \begin{scope}[xshift=2cm, yshift=-2cm]
                \draw[lightgray] (0,0) grid[step=0.2] (2,2); 
                \draw (0,0) grid[step=2] (2,2); 
                \node at (1,1) {3}; 
            \end{scope}

            \draw[black, line width = 3] (1.8, 0) -- (2.2, 0); 
            \draw[black, line width = 3] (2, 0) -- (2, -0.4); 

        \end{tikzpicture}
        \caption{Physical domain.}
    \end{subfigure}
    \caption{\label{Tjoint_example}
        Example of a T-joint between three patches. The mapping between the logical domains
        and the physical domain is the identity.  
        We can see that part of an edge of patch 1 is connected to an edge of patch 2 on a subdomain
        and to an edge of patch 3 on a complementary subdomain. This type of connection is called a 
        \textit{T-joint}.
    }
\end{figure}

\paragraph{}
We define a \textit{T-joint} as a connection between at least three patches, 
where an edge of one patch is not fully connected to the edge of another patch along its entire length.
An illustration between three patches is given in Fig. \ref{Tjoint_example}. 
The introduction of T-joints may be necessary for some complex geometries involving an
O-point and an X-point. This can be observed in Fig. \ref{JET_SOLEDGE}.
We will analyze a case with T-joints on a domain with an O-point.

\subsubsection{Application of the formula to the T-joint cases}
\paragraph{}
The computation of the derivatives on geometries with T-joints stays pretty similar
to the cases mentioned before. For conforming meshes, we remark that the cross-derivatives 
can be computed along one of the directions without a difference
in the result. On geometries with T-joints, the choice of the direction may be influenced 
by the geometry. For instance, in the example of Fig. \ref{Tjoint_example}, the cross-derivatives
will be computed along the $r$-direction. Some $r$-derivatives on patch 1, needed to apply equation 
\eqref{relations_derivatives_2D} along the $\theta$-direction, are not available. 

\begin{figure}[h]
    \begin{subfigure}[t]{0.5\textwidth}
        \centering
        \begin{tikzpicture}[xscale=0.75, yscale=0.75]
            \draw[lightgray] (0,0) grid[step=0.2] (2,6); 
            \draw (0,0) grid[xstep=2, ystep=6] (2,6); 
            \node at (1,3) {1}; 
            \draw[->] (0.25, 0.25) -- (0.75, 0.25) node[right] {$r_1$}; 
            \draw[->] (0.25, 0.25) -- (0.25, 0.75) node[above] {$\theta_1$}; 

            \begin{scope}[xshift=2.25cm, yshift=0]
                \draw[lightgray] (0,0) grid[step=0.2] (2,2); 
                \draw (0,0) grid[step=2] (2,2); 
                \node at (1,1) {2}; 
                \draw[->] (0.25, 0.25) -- (0.75, 0.25) node[right] {$r_2$}; 
                \draw[->] (0.25, 0.25) -- (0.25, 0.75) node[above] {$\theta_2$}; 
            \end{scope}

            \begin{scope}[xshift=2.25cm, yshift=2.25cm]
                \draw[lightgray] (0,0) grid[step=0.2] (2,2); 
                \draw (0,0) grid[step=2] (2,2); 
                \node at (1,1) {3}; 
                \draw[->] (0.25, 0.25) -- (0.75, 0.25) node[right] {$r_3$}; 
                \draw[->] (0.25, 0.25) -- (0.25, 0.75) node[above] {$\theta_3$}; 
            \end{scope}

            \begin{scope}[xshift=2.25cm, yshift=4.5cm]
                \draw[lightgray] (0,0) grid[step=0.2] (2,2); 
                \draw (0,0) grid[step=2] (2,2); 
                \node at (1,1) {4}; 
                \draw[->] (0.25, 0.25) -- (0.75, 0.25) node[right] {$r_4$}; 
                \draw[->] (0.25, 0.25) -- (0.25, 0.75) node[above] {$\theta_4$}; 
            \end{scope}

            \begin{scope}[xshift=4.5cm, yshift=0]
                \draw[lightgray] (0,0) grid[step=0.2] (2,6); 
                \draw (0,0) grid[xstep=2, ystep=6] (2,6); 
                \node at (1,3) {5};
                \draw[->] (0.25, 0.25) -- (0.75, 0.25) node[right] {$r_5$}; 
                \draw[->] (0.25, 0.25) -- (0.25, 0.75) node[above] {$\theta_5$};  
            \end{scope}
        \end{tikzpicture}
        \caption{Logical domains.}
    \end{subfigure}
    \begin{subfigure}[t]{0.5\textwidth}
        \centering
        \begin{tikzpicture}[xscale=0.70, yscale=0.70]
            \draw[->] (-3.25, -3.25) -- ( 3.25, -3.25) node[right] {$x$}; 
            \draw[->] (-3.25, -3.25) -- (-3.25,  3.25) node[above] {$y$}; 

            \foreach \r in {0, 0.1, ..., 3}{
                \draw[lightgray, line width=0.5] (0,0) circle (\r); 
            }
            \foreach \theta in {0,6, ..., 360}{
                \draw[lightgray, line width=0.5] (\theta:0) -- (\theta:3);
            }
            \draw[line width = 1] (0,0) circle (1); 
            \draw[line width = 1] (0,0) circle (2); 
            \draw[line width = 1] (0,0) circle (3); 
            \draw[line width = 1] (0:1) -- (0:2);
            \draw[line width = 1] (120:1) -- (120:2);
            \draw[line width = 1] (240:1) -- (240:2);

            \node at (0,0)     {1}; 
            \node at (60:1.5)  {2}; 
            \node at (180:1.5) {3}; 
            \node at (300:1.5) {4}; 
            \node at (90:2.5)  {5}; 
        \end{tikzpicture}
        \caption{Physical domain.}
    \end{subfigure}
    \caption{\label{Tjoint_study_case}
        Example of the geometry with T-joints and an O-point studied for our simulations 
        in the Section \ref{section_Numerical_results}. 
        On the left, the patches in the logical domains. On the right, the physical domain. 
        The global mesh is a conforming mesh. 
    }
\end{figure}
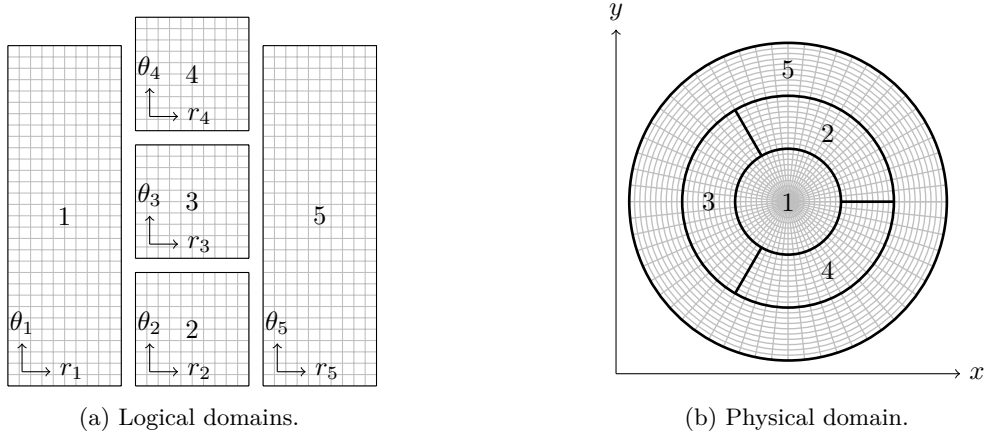

\paragraph{Example applied later in our simulation.}
\label{Example_T_joints_geometry}
Let's study the geometry in Fig. \ref{Tjoint_study_case}, which will be an application case for the 
simulations presented in the Section \ref{section_Numerical_results} of this article. 
The global mesh is conforming, so we can build an equivalent global spline on this domain. 
To compute all the derivatives, we start by computing the derivatives along $r$
on each patch using \eqref{relations_derivatives_2D}. We then have all the necessary information 
to build spline representations on patches 1 and 5, but not on patches 2, 3 and 4. 
The derivatives still need to be computed along $\theta$ on patches 2, 3 and 4. 
From the derivatives along $r$, we compute the cross-derivatives at the T-joints. 
All the required derivatives are then computed, 
and we can build spline representations on patches 2, 3 and 4.

\section{Numerical results}
\label{section_Numerical_results}
\paragraph{}
Here, we present some numerical results of the method presented in the previous sections.
The implementation can be found in the Gyselalib++ library \cite{Gyselalibxx}.
This section is organized as follows.
The first results are spline interpolation results on three concentric patches. 
We check that the approximation formula returns the expected precision. 
Secondly, we test the interpolation in an advection simulation with a backward semi-Lagrangian method. 
Finally, results corresponding to more complete simulation on a guiding-center model, 
\eqref{guiding_center_eq}, involving a coupling with a Poisson equation are shown. 

\subsection{Spline interpolation tests}
\paragraph{}
We compare the results of an interpolation on local splines built with the derivatives 
computed using the approximation relations and with the exact formula to an equivalent
global spline. 
We test it with the following function 
\begin{equation*}
    f: \widehat{\Omega} \rightarrow \mathbb{R}: 
    (r, \theta) \mapsto \cos(2\pi x(r, \theta)) \sin(2\pi y(r,\theta)), 
\end{equation*}
with $x$ and $y$ given by the Czarny mapping \cite{Czarny2008}, 
\begin{equation}
    \mathcal{F}: \widehat{\Omega} \rightarrow \Omega: 
    (r, \theta) \mapsto 
    \left\{
    \begin{aligned}
        & x(r, \theta) 
            = \frac{1}{\varepsilon}(1 - \sqrt{1 + \varepsilon(\varepsilon + 2r\cos(\theta))}),  \\
        & y(r, \theta) 
            = \frac{e\xi r\sin(\theta)}{2 - \sqrt{1 + \varepsilon(\varepsilon + 2r\cos(\theta))}}, \\
    \end{aligned}
    \right. 
    \label{Czarny_mapping}
\end{equation}
with $\xi = \frac{1}{\sqrt{1 - \frac{1}{4}\varepsilon^2}}$, and $\varepsilon$ and $e$ given constants, 
that we will fix to $\varepsilon = 0.3$ and $e = 1.4$.
The interpolation is carried out on the logical domain. The physical domain is not involved. 
We introduced a mapping in the test function to have curves which do not follow the mesh lines, 
as will be the case in the next simulations.

The global domain is split into three patches. The splitting of the global domain is only along the $r$-direction. 
It is a 1D multi-patch case as only the $r$-grid is split. 
For the patch connections, we only need to compute the r-derivatives.
We test two cases: a uniform global mesh and a non-uniform global mesh. 
In the non-uniform global mesh case, the patches are uniform, but the cell length differs from one patch 
to another. 
In the following tests, the number of cells in each patch stays the same in both cases but the position 
of the interfaces changes.
The grid parameters of the patches and the global domain are given in Table \ref{Interpolation_table_grid}.

\begin{table}[h]
    \centering
    \footnotesize
    \begin{tabular}{ c|c|c|c|c|c } 
        \hline
        Test & Patch & Number of cells & Domain along $r$ & $dr$  & Domain along $\theta$ \\ 
        \hline\hline
        1.1 Uniform 
            & Patch1 & $[42\times 256]$ & $\left[0, 42/128\right]$      & $7.8125\cdot10^{-3}$ & $\left[0, 2\pi\right]$ \\
            & Patch2 & $[44\times 256]$ & $\left[42/128,86/128\right]$  & $7.8125\cdot10^{-3}$ & $\left[0, 2\pi\right]$ \\
            & Patch3 & $[42\times 256]$ & $\left[86/128, 1\right]$      & $7.8125\cdot10^{-3}$ & $\left[0, 2\pi\right]$ \\
            \arrayrulecolor{lightgray} \cline{2-5}  \arrayrulecolor{black}
            & Global & $[128\times 256]$ & $\left[0, 1\right]$          & $7.8125\cdot10^{-3}$ & $\left[0, 2\pi\right]$ \\
        \hline 
        1.2 Non-uniform 
            & Patch1 & $[42\times 256]$ & $\left[0, 42/128\right]$      & $7.8125\cdot10^{-3}$ & $\left[0, 2\pi\right]$ \\
            & Patch2 & $[44\times 256]$ & $\left[42/128, 60/128\right]$ & $3.1960\cdot10^{-3}$ & $\left[0, 2\pi\right]$ \\
            & Patch3 & $[42\times 256]$ & $\left[60/128, 1\right]$      & $1.2649\cdot10^{-2}$ & $\left[0, 2\pi\right]$ \\
            \arrayrulecolor{lightgray} \cline{2-5}  \arrayrulecolor{black}
            & Global & $[128\times 256]$ & $\left[0, 1\right]$          &                      & $\left[0, 2\pi\right]$ \\
        \hline
    \end{tabular}
    \caption{ \label{Interpolation_table_grid}
        Number of cells and domains of the patches and equivalent global domain 
        for the spline interpolation tests. 
        }
\end{table}

The results for different numbers of cells for the approximation \eqref{Matrix_Hermite_approx}
and the exact formula \eqref{relation_interfaces} are shown in Table \ref{Interpolation_tables}.
We compare the values of the interface derivatives computed with the formulae to the 
derivatives of an equivalent global spline. We also compare the values of the local 
and global spline representations on the interpolation points. The given errors are the 
$L_\infty$-norm of the difference, 
\begin{equation*}
    err_{\text{funct}} 
        = || s_{local} -  s_{global}||_{L_\infty\left(\widehat{\Omega}\right)}, 
    \qquad
    err_{\text{deriv}} 
        = || \partial_r s_{local} -  \partial_r  s_{global}||_{L_\infty\left(\partial\widehat{\Omega}^k\right)}. 
\end{equation*}
We observe that the results in Table \ref{Interpolation_tables} agree with the estimations in 
Table \ref{table_equivalence_diff_n}. Moreover, as expected, for thirty cells, 
the computed interface derivative values reach the same precision as if we use all the cells 
with the exact formula.

\begin{table}[h]
    \begin{subtable}[b]{0.5\textwidth}
        \centering
        \footnotesize
        \begin{tabular}{ |c|c|c| } 
            \hline
            Number & Errors on local  & Errors on interfaces  \\ 
            of cells  & splines & derivatives \\ 
            \hline\hline
            5 & 5.15e-05 & 4.16e-02 \\
            \hline
            10 & 6.37e-08 & 5.15e-05 \\
            \hline
            20 & 7.86e-14 & 6.35e-11 \\
            \hline
            30 & 8.44e-15 & 2.37e-13 \\
            \hline
            40 & 8.44e-15 & 2.36e-13 \\
            \hline
            Exact & 8.44e-15 & 2.82e-13 \\
            \hline
        \end{tabular}
        \caption{Uniform mesh between the three patches.}
    \end{subtable}
    \begin{subtable}[b]{0.5\textwidth}
        \centering
        \footnotesize
        \begin{tabular}{ |c|c|c| } 
            \hline
            Number & Errors on local & Errors on interfaces \\ 
            of cells  & splines & derivatives \\ 
            \hline\hline
            5 & 4.61e-05 & 2.57e-02 \\
            \hline
            10 & 6.08e-08 & 3.23e-05 \\
            \hline
            20 & 9.55e-14 & 4.77e-11 \\
            \hline
            30 & 8.44e-15 & 2.59e-13 \\
            \hline
            40 & 8.44e-15 & 3.84e-13 \\
            \hline
            Exact & 8.44e-15 & 3.36e-13 \\
            \hline
        \end{tabular}
        \caption{Non-uniform mesh between the three patches.}
    \end{subtable}
    \caption{\label{Interpolation_tables}
        Interpolation results.
        The results are given for the uniform mesh on the left and the non-uniform mesh 
        on the right with the configurations from Table \ref{Interpolation_table_grid}. 
        The given errors are the difference of the results on local splines 
        defined on patch 1, patch 2 and patch 3 to an equivalent global spline. 
        The indicated number of cells corresponds to the number of cells taken 
        in the approximation of the derivatives. 
        We see that indeed for $30$, it reaches the same precision as if we 
        take all the cells with the exact formula. 
        The results agree with the analytical ones given in Table \ref{table_equivalence_diff_n}
        until saturation at $10^{-13}$ (for a space step around $10^{-2}$). 
        }
\end{table}

\begin{figure}[H]
    \centering
    \begin{subfigure}[h]{\textwidth}
        \includegraphics[width=\textwidth,trim={2cm 11.8cm 2cm 0.5cm}, clip=true]{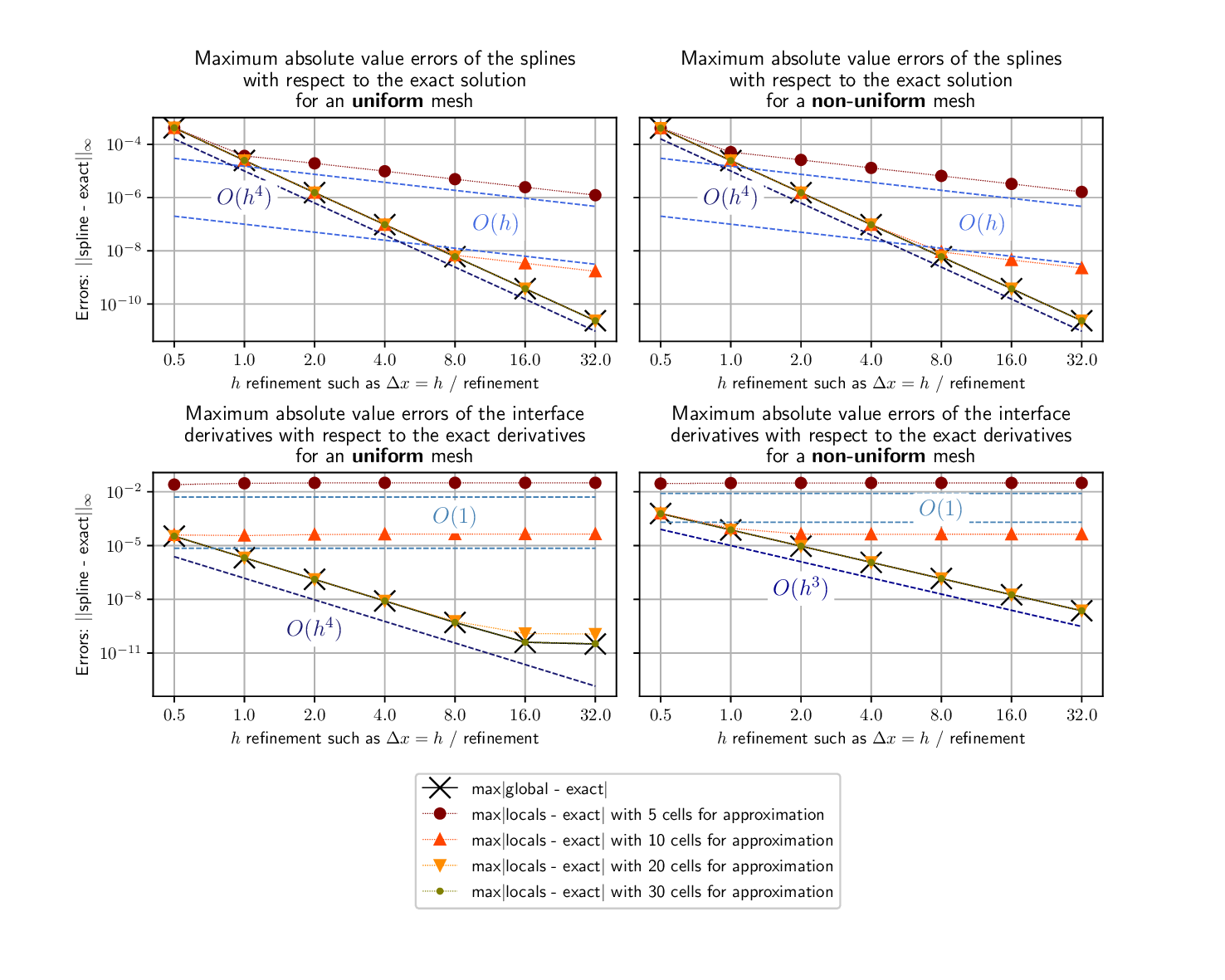}
        \caption{
            \label{interpolation_convergence_order_a}
            Convergence orders on the interpolation of the local splines and the equivalent global spline. 
        }
    \end{subfigure}
    \begin{subfigure}[h]{\textwidth}
        \includegraphics[width=\textwidth,trim={2cm 1.0cm 2cm 8.0cm}, clip=true]{Fig14.eps}
        \caption{
            \label{interpolation_convergence_order_b}
            Convergence order of the computed interface derivatives and the derivatives 
            at the interfaces of the global spline. 
        }
    \end{subfigure}
    \vspace*{-0.25cm}

    \caption{
        \label{interpolation_convergence_order}
        Convergence orders. 
        The results are given for the uniform mesh on the left and the non-uniform mesh on the right
        with the configurations from Table \ref{Interpolation_table_grid}. 
        The legend in subfigure (b) is the legend for both subfigures.
        The $\Delta x$ given for the abscissa indicates the refinement in both dimensions. 
        It corresponds to $dr$ and $d\theta$ with their associated $h_r = 1 / 128$ 
        and $h_\theta = 2\pi / 256$ in the uniform case.
    }
\end{figure}

The convergence orders of the local splines compared to the equivalent global spline are indicated in 
Fig. \ref{interpolation_convergence_order_a}.
For an approximation of the interface derivatives using 20 cells or 30 cells, 
the local splines converge similarly to the global spline with an order of 4. 
For an approximation of the interface derivatives using 10 cells, 
we observe a saturation starting with a refinement of 8 of 
the initial grid given in Table \ref{Interpolation_table_grid}, 
i.e., for patches on $[336\times2048]$, $[352\times2048]$ and $[336\times2048]$ cells.
With an approximation of the interface derivatives using only 5 cells, 
this saturation already starts for a refinement of 1, 
i.e., the grid given in Table \ref{Interpolation_table_grid}.

To understand the saturation, we need to understand where the errors come from. 
Let's write $\varepsilon_{tot}$ the error of the local splines with respect to the exact solution, 
$\varepsilon_{global}$ the error of the equivalent global spline with respect to the exact solution, 
and
$\varepsilon_{local}$ the error of the local splines with respect to the equivalent global spline. 
We have the following relation, 
\begin{equation*}
    |\varepsilon_{tot}| \leq |\varepsilon_{global}| + |\varepsilon_{local}|.
\end{equation*}

Fig. \ref{interpolation_convergence_order_a} confirms that $|\varepsilon_{global}|$ is order 4. 
This is the expected convergence order for a cubic spline. 

When we compute the interface derivatives with the approximation of $N_a$ cells instead of 
the exact formula, the error $\varepsilon_{deriv \ local}$ we make is 
\begin{equation*}
    |\varepsilon_{deriv \ local}|
        = |a^I_{N_a,N_a} s'(x_{i+N_a})
        + b^I_{N_a,N_a} s'(x_{i-N_a})|
        < (|a^I_{N_a,N_a}| + |b^I_{N_a,N_a}|) \  || s'||_{L_\infty}.
\end{equation*}

This can be seen as the exact formula applied to subgrids on the patches of $N_a$ cells,
where we removed the dependency on the other derivatives. 

As mentioned earlier in section \ref{subsubsection_independence_of_coef_from_the_function}, 
$a^I_{N_a,N_a}$ and $b^I_{N_a,N_a}$ do not depend directly on the cell lengths 
but on the proportionality between the cell lengths. While refining, they stay constant. 
If we refer to Table \ref{table_equivalence_diff_n}, we see that for a global uniform mesh, 
$|a^I_{5, 5}| = |b^I_{5, 5}| = 4\times5.52 \cdot 10^{-3}$ so $|a^I_{5, 5}| + |b^I_{5, 5}| = 4.42 \cdot 10^{-2}$. 
Similarly, $|a^I_{10, 10}| + |b^I_{10, 10}| = 6.10\cdot 10^{-5}$. 

Similar values can be found in Table \ref{Interpolation_tables}, 
which gives the exact errors on the interface derivatives for the chosen function. 
It indicates that with a 5-cell approximation in the uniform case, 
$|\varepsilon_{deriv \ local}| = 4.16\times10^{-2}$. 
Fig. \ref{interpolation_convergence_order_b} shows that this error is bigger than 
$|\varepsilon_{deriv \ global}|$ 
(in black on the figure). This explains why the saturation is already present from a refinement of 0.5. 
We have similar results for an approximation with 10 cells. 
With an approximation of 20 cells, Table \ref{Interpolation_tables} indicates an error 
of $|\varepsilon_{deriv \ local}| = 6.35\times10^{-11}$ (in the uniform case). This value is reached by 
$|\varepsilon_{deriv \ global}|$ after a refinement of $h / 16$. 
We can then see the saturation appearing in Fig. \ref{interpolation_convergence_order_b}.

For a convergence order $p$ given on derivatives, the associated spline has a convergence order of $p+1$. 
Here, the convergence order of the derivatives is zero, so the spline has an order one of convergence 
when the saturation is reached. 

This constant error in the interface derivatives causes the saturation in the local spline convergence order
that we observe. For enough cells in the approximation of the interface derivatives, 
the behaviors of the local splines and the equivalent global spline are the same. 

In the uniform case, the convergence of the derivatives of a cubic spline is order four 
\cite{MehrenbergerHDR, Besse2008}. 
This explains the difference that can be seen in Fig. \ref{interpolation_convergence_order_b}
between the uniform and non-uniform cases.

\subsection{Advection tests}
\paragraph{}
We test the method on an advection with a test case similar to the one detailed in \cite{Zoni2019}.
The exact formula \eqref{Matrix_Hermite} is applied, so we expect the same results between the 
local splines and the equivalent global spline.
The initial condition on the global domain is given by, 
\begin{equation}
    \left\{
    \begin{aligned}
        & \rho(x,y) = \frac{1}{2}\left[G(r_1(x,y)) + G(r_2(x,y))\right], 
            \qquad (x,y)\in\Omega^2\\ 
        & r_1(x,y) = \sqrt{(x - x_0)^2 + 8(y - y_0)^2}, \\
        & r_2(x,y) = \sqrt{8(x - x_0)^2 + (y - y_0)^2}, \\
        & G(r) = \cos\left(\frac{\pi r}{2 a}\right)^4  \mathbbm{1}_{r \leq a},
    \end{aligned}
    \right.
    \label{function_advection_test}
\end{equation}
with $(x_0, y_0) = \mathcal{F}(0.5,0)$ and $a = 0.3$. The physical domain is defined by 
the Czarny mapping $\mathcal{F}$ \eqref{Czarny_mapping}, \cite{Czarny2008}. 
This function is advected according to the following advection field $\mathbf{A}$,
\begin{equation}
    \mathbf{A}
    = 
    \mathbf{J}_{\mathcal{F}}
    \begin{bmatrix}
        0 \\
        \omega \\
    \end{bmatrix},
    \label{advection_field_test}
\end{equation}
with $\mathbf{J}_{\mathcal{F}}$ the Jacobian matrix of the mapping $\mathcal{F}$ \eqref{Czarny_mapping}, 
and $\omega = 2\pi$.

The trajectory is a rotation following the mesh lines.
The results in \cite{Zoni2019} are given for a uniform global domain defined on 
$N_r\times N_\theta = 128\times256$.
The chosen time step is $\Delta t = 0.01$. 
To find the feet of the characteristics, we apply a third order Runge-Kutta method.

In our tests, we split the domain into three (concentric rings) for (Test 2.1) or five patches 
(see Fig. \ref{Tjoint_study_case}) for (Test 2.2). 
We apply the exact formulae \eqref{Matrix_Hermite} in the first test and \eqref{relations_derivatives_2D}
in the second. 
The domain discretizations of these two tests are detailed in Table \ref{Advection_table_grid}. 

\begin{table}[h]
    \centering
    \footnotesize
    \begin{tabular}{ c|c|c|c|c } 
        \hline
        Test & Patch & Number of cells & Domain along $r$ & Domain along $\theta$ \\ 
        \hline\hline
        2.1 Non-conforming (1D) 
            & Patch1 & $[21\times 128]$ & $\left[0, 42/128\right]$ & $\left[0, 2\pi\right]$ \\
            & Patch2 & $[76\times 256]$ & $\left[42/128, 80/128\right]$ & $\left[0, 2\pi\right]$ \\
            & Patch3 & $[24\times 128]$ & $\left[80/128, 1\right]$ & $\left[0, 2\pi\right]$ \\
            \arrayrulecolor{lightgray} \cline{2-5}  \arrayrulecolor{black}
            & Refined global & $[128\times 256]$ & $\left[0, 1\right]$ & $\left[0, 2\pi\right]$ \\
            & Coarse global  & $[128\times 128]$ & $\left[0, 1\right]$ & $\left[0, 2\pi\right]$ \\
        \hline
        2.2 T-joints (2D) 
            & Patch1 & $[42\times 255]$ & $\left[0, 42/128\right]$ & $\left[0, 2\pi\right]$ \\
            & Patch2 & $[38\times  85]$ & $\left[42/128, 80/128\right]$ & $\left[0, 2\pi/3\right]$ \\
            & Patch3 & $[38\times  85]$ & $\left[42/128, 80/128\right]$ & $\left[2\pi/3,4\pi/3\right]$ \\
            & Patch4 & $[38\times  85]$ & $\left[42/128, 80/128\right]$ & $\left[4\pi/3, 2\pi\right]$ \\
            & Patch5 & $[48\times 255]$ & $\left[80/128, 1\right]$ & $\left[0, 2\pi\right]$ \\
            \arrayrulecolor{lightgray} \cline{2-5}  \arrayrulecolor{black}
            & Global & $[128\times 255]$ & $\left[0, 1\right]$ & $\left[0, 2\pi\right]$ \\
        \hline
    \end{tabular}
    \caption{ \label{Advection_table_grid}
        Number of cells and domains of the patches and equivalent global domains for the 
        advection tests. 
        In parentheses, we indicate the dimensions of the multi-patch mesh. If the multi-patch mesh is 1D, 
        we only need to compute the derivatives along $r$. If it is 2D, we compute the 
        derivatives along $r$, along $\theta$ and the cross-derivatives. 
        }
\end{table}

\begin{figure}[H]
    \centering
    \includegraphics[width=\textwidth]{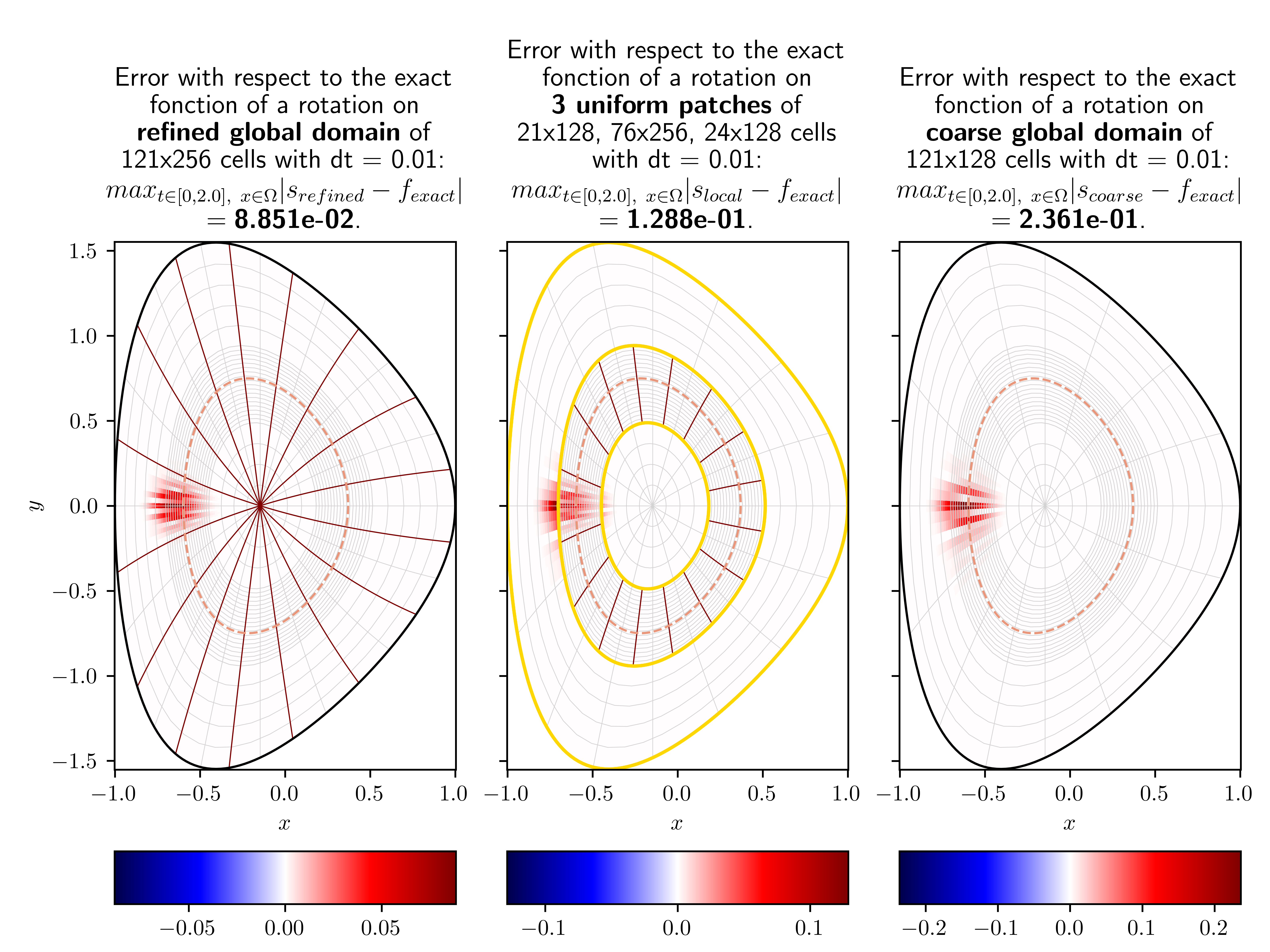}

    \vspace*{-2.5mm}
    \caption{ \label{test_2.2}
        Advection test 2.1: non-conforming patches.
        There is no comparison to an equivalent global spline as it does not exist, 
        but we compare here the errors to the exact solution with the errors of a refined spline 
        and a coarse spline.
        Errors of the spline representations to the exact solution at $t = 2.0$. 
        The dashed line shows the trajectory of the maximum density point.
    }
\end{figure}

\begin{figure}[H]
    \centering
    \includegraphics[width=\textwidth]{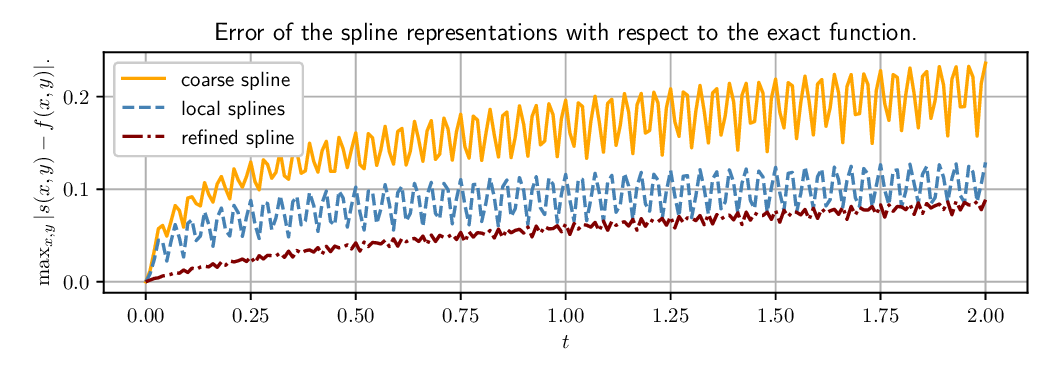}

    \vspace*{-2.5mm}

    \caption{ \label{test_2.2_error}
        Advection test 2.1: non-conforming patches.
        Time variation of the errors of the spline representations with respect to the exact solution
        at each time steps. The errors of the local splines are bounded by the errors of the 
        refined spline and the errors of the coarse spline as expected.  
        We also observe that the local splines on the non-conforming mesh tend to be as good 
        as the spline on the refined mesh. 
    }
\end{figure}

\paragraph{Test 2.1 Non-conforming patches.}
Patch 1 and Patch 3 are made coarser along the $\theta$ direction. 
The global mesh is then no longer conforming. After computing the derivatives along
$r$ which are on the conforming lines, we build $\theta$-spline representations of the $r$-derivatives
to compute the $r$-derivatives on the non-conforming lines at the boundaries of Patch 2. 
Fig. \ref{test_2.2} and Fig. \ref{test_2.2_error} show the results. 
Here, there is no equivalent global spline. 
So, we cannot expect a machine precision agreement to a global spline. However, we can build a refined
global spline and a coarse global spline and compare the three spline representations to the exact solution. 
We expect the error with respect to the exact solution of the local splines to be bounded between 
the error of the refined global spline and the coarse global spline. 
This is what we observe in Fig. \ref{test_2.2_error}. 
The results of this figure also highlight the efficiency of the non-conforming mesh: despite using fewer cells, 
the local splines tend to be as precise as the refined spline. 
The non-conforming mesh is, of course, adapted to the solution trajectory presented in Fig. \ref{test_2.2}. 
The errors we obtain in Fig. \ref{test_2.2_error} for $t = 0.8$ agree with the error obtained in 
\cite{Zoni2019}, see Table 2.


\begin{figure}[H]
    \centering
    \includegraphics[width=\textwidth]{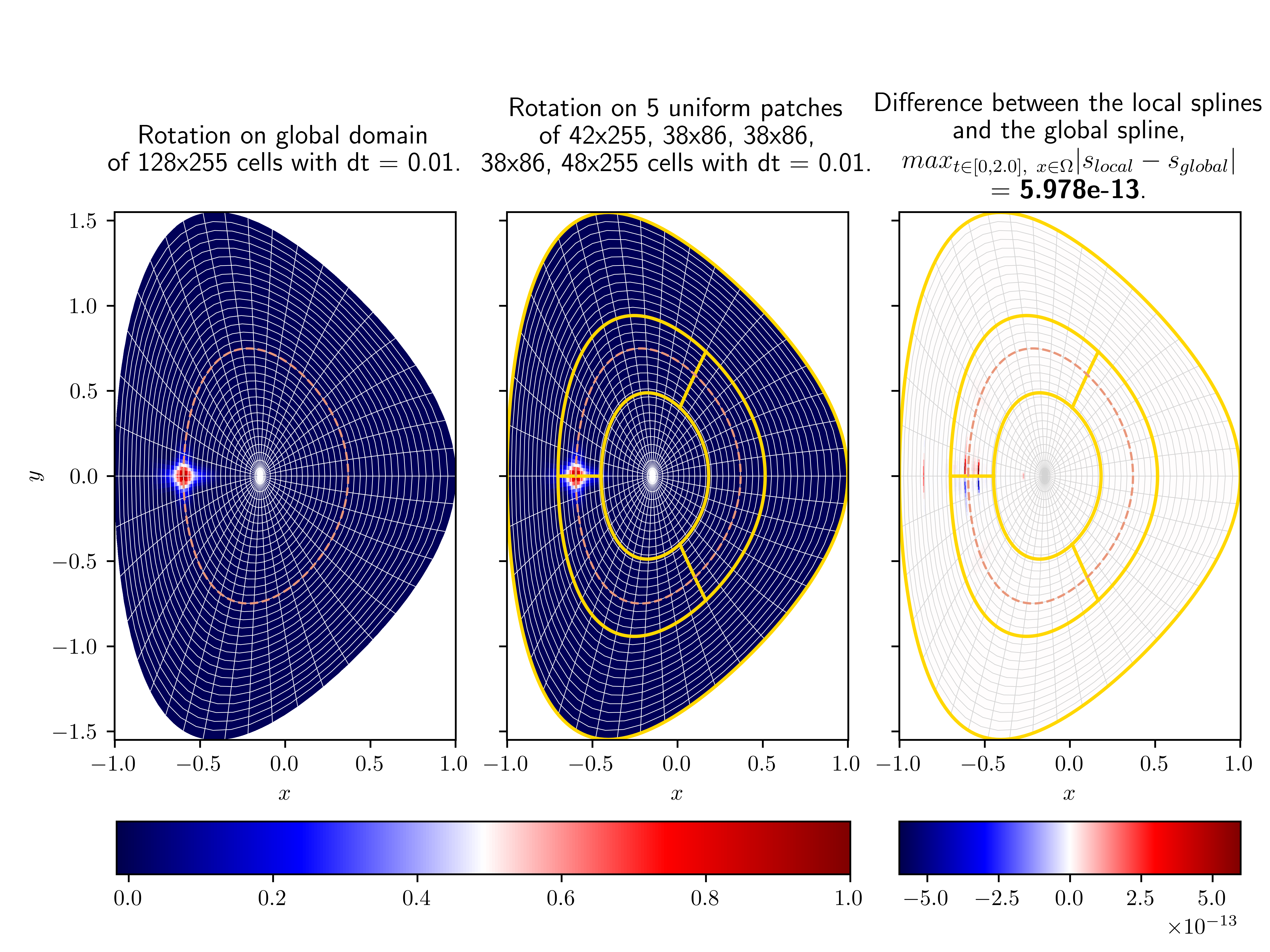}

    \vspace*{-2.5mm}
    
    \caption{ \label{test_2.3}
        Advection test 2.3: T-joints.
        Difference between the results on multiple patches and the results on an 
        equivalent global spline. 
        The dashed line shows the trajectory of the maximum density point. 
        We can see that at $t = 2.0$ (200 iterations), 
        $||s_{\text{local}} - s_{\text{global}}||_{L_\infty} = 5.978\cdot 10^{-13}$. 
        The two spline representations agree even when the trajectory crosses the interfaces. 
    }
\end{figure}

\paragraph{Test 2.2 T-joints.}
We apply the method to a case with T-joints. 
Here the multi-patch splitting is 2D, so 
we need to compute the $r$-derivatives, the $\theta$-derivatives and the cross-derivatives
as explained in paragraph \ref{Example_T_joints_geometry} for the studied geometry described 
in Fig. \ref{Tjoint_study_case}. 
We choose a uniform global mesh. 
So, we expect a difference between the final local spline representations
and the final equivalent global spline representation at machine precision. 
Fig. \ref{test_2.3} shows that the spline representations agree. 
We obtain an error of $6\cdot 10^{-13}$ due to machine precision. 

\subsection{Guiding-center model - diocotron instability simulation}
\paragraph{}
The diocotron instability simulation is a reference numerical test for the guiding-center equations 
\eqref{guiding_center_eq}. 
Physically, it simulates a non-neutral plasma in a cylinder, 
whose particles are confined radially by a uniform axial magnetic field 
with a cylindrical conducting wall located at the outer boundary.
The simulation is detailed in 
\cite{Davidson1990, Petri2007, Petri2009, Madaule2013, Zoni2019, ZoniPhD}.

The electric potential $\phi$ is null at the boundaries. 
The initial density profile is  given by 
\begin{equation}
    \begin{aligned}
        \rho(t=0, r,\theta) 
        = \rho_0(r) 
        + \rho_1(t=0, r,\theta) 
        = \left(1 + \varepsilon\cos(m\theta)\right)
            \exp\left[\left(-\frac{r - \bar{r}}{d}\right)^p\right] 
            \mathbbm{1}(r^{-}\leq r \leq r^{+}),
    \end{aligned}
    \label{diocotron_initial_condition}
\end{equation}
with $\rho_0$ the equilibrium and $\rho_1$, the perturbation of amplitude $\varepsilon$
and mode $m$. The initial density profile is smoothed by an exponential to avoid discontinuities.

The functions are defined on $\widehat{\Omega} = [0, 1]\times[0, 2\pi]$. 
The physical domain is the arrival space of the following circular mapping, 
\begin{equation*}
    \mathcal{F} : \widehat{\Omega} \rightarrow \Omega, \  
    (r, \theta) \mapsto 
    \left\{
    \begin{aligned}
        & x = r\cos(\theta), \\
        & y = r\sin(\theta). \\
    \end{aligned}
    \right.
\end{equation*}

\begin{figure}[h]
    \centering
    \includegraphics[width=\textwidth]{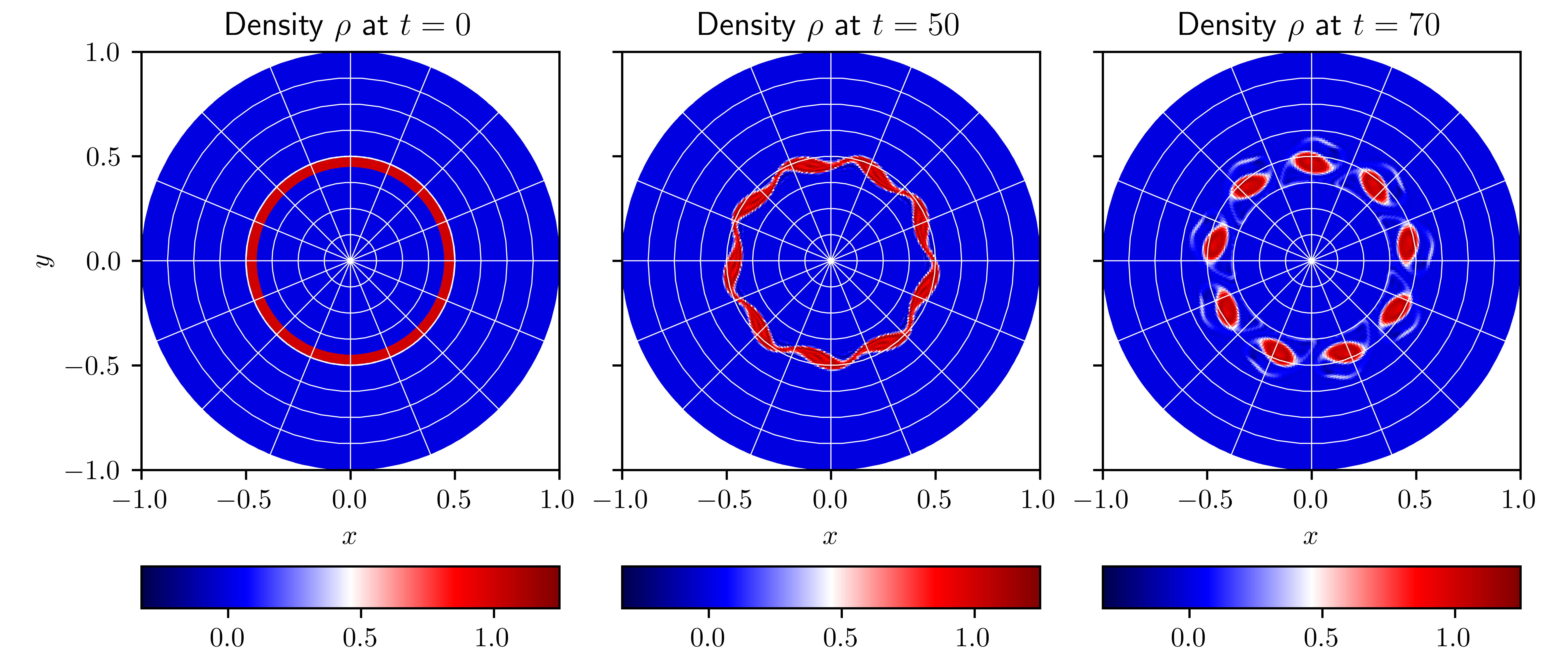}

    \vspace*{-0.25cm}
    
    \caption{\label{diocotron_steps}
        Density profile for different time steps. 
        On the right, the initial condition ($t = 0.0$). 
        In the middle, end of the linear regime ($ = 50.0$). 
        On the left, after the linear regime ($t = 70.0$), $m = 9$ vortexes are formed. 
        This is also the final time chosen in the simulation in the article \cite{Zoni2019}.
    }
\end{figure}

The reference simulation in \cite{Zoni2019}, 
is run on a grid of $N_r\times N_\theta = 128\times 256$ cells, with a time 
step of $dt = 0.1$. 
The initial condition is defined with the amplitude and mode perturbation, 
$\varepsilon = 10^{-4}$ and $m = 9$. It is a ring of diameters $r^- = 0.45$ 
and $r^+ = 0.50$. The parameters of the smoothing function are: 
$\bar{r} = \frac{r^- + r^+}{2}$, $d= \frac{r^+ - r^- }{2}$ and $p = 50$. 
The coupling of the two equations is made with an explicit predictor-corrector method 
described in \cite{Zoni2019}. Fig. \ref{diocotron_steps} shows the 
evolution of the density profile for different time steps during the linear regime.
The perturbation introduced in the initial condition grows linearly and creates $m$
vortices. 

In the article \cite{Zoni2019}, the simulations are run on a global spline. We will compare 
the results of this simulation to a simulation run with local splines. 
To compute the derivatives at the interface, we will use the exact formulae
developed in \eqref{relation_interfaces} and \eqref{relations_derivatives_2D}.
The difference between a simulation on a global spline and a simulation on 
local splines is computed by taking the infinity norm of the difference of the 
density profile, 
\begin{equation*}
    \max_{t\geq 0} 
    || \rho_{\text{global}} - \rho_{\text{local}} ||
    _{L_\infty}
    = 
    \max_{t\geq 0}
    \left(
        \max_{(r_i,\theta_j)}
            (|\rho_{\text{global}}(r_i,\theta_j) 
            - \rho_{\text{local}}(r_i,\theta_j)|)
    \right).
\end{equation*}

Two tests have been implemented in Gyselalib++ and are presented below: 
(Test 3.1) an introduction of patches on the same simulation as in the article \cite{Zoni2019}
in \ref{diocotron_uniform_section}; 
(Test 3.2) a modification of the simulation mesh to add a local refinement
with an interface in the middle of the solution
and a cell depopulation of the central patch
in a 2D multi-patch geometry with T-joints
in \ref{diocotron_Tjoints_refined_section}.

\paragraph{Remark.}
The study in this article focuses on the advection on a multi-patch grid. The Poisson 
equation of the guiding-center equations \eqref{guiding_center_eq} is solved 
on an equivalent global spline of the local splines used for the advection. The method used is
exactly the same as described in the reference article \cite{Zoni2019}.
It has been implemented in Gyselalib++ \cite{BournePhD}. 
For this reason, the results presented below will be only on conforming global meshes.

\paragraph{List of tests.}
The diocotron instability simulation has been run on different multi-patch grids. 
Table \ref{Diocotron_table_grid} contains a list of the different parameters we chose for the grids. 

\begin{table}[h]
    \centering
    \footnotesize
    \begin{tabular}{ c|c|c|c|c } 
        \hline
        Test & Patch & Number of cells & Domain along $r$ & Domain along $\theta$ \\ 
        \hline\hline
        3.1 Uniform (1D) 
            & Patch1 & $[42\times 256]$ & $\left[0, 42/128\right]$ & $\left[0, 2\pi\right]$ \\
            & Patch2 & $[44\times 256]$ & $\left[42/128, 86/128\right]$ & $\left[0, 2\pi\right]$ \\
            & Patch3 & $[42\times 256]$ & $\left[86/128, 1\right]$ & $\left[0, 2\pi\right]$ \\
            \arrayrulecolor{lightgray} \cline{2-5}  \arrayrulecolor{black}
            & Global & $[128\times 256]$ & $\left[0, 1\right]$ & $\left[0, 2\pi\right]$ \\
        \hline 
        3.2 T-joints with refinement (2D) 
            & Patch1 & $[4\times 255]$ & $\left[0, 42/128\right]$ & $\left[0, 2\pi\right]$ \\
            & Patch2 & $[72\times  85]$ & $\left[42/128, 60/128\right]$ & $\left[0, 2\pi/3\right]$ \\
            & Patch3 & $[72\times  85]$ & $\left[42/128, 60/128\right]$ & $\left[2\pi/3, 4\pi/3\right]$ \\
            & Patch4 & $[72\times  85]$ & $\left[42/128, 60/128\right]$ & $\left[4\pi/3, 2\pi\right]$ \\
            & Patch5 & $[48\times 255]$ & $\left[60/128, 1\right]$ & $\left[0, 2\pi\right]$ \\
            \arrayrulecolor{lightgray} \cline{2-5}  \arrayrulecolor{black}
            & Global & $[124\times 255]$ & $\left[0, 1\right]$ & $\left[0, 2\pi\right]$ \\
        \hline
    \end{tabular}
    \caption{ \label{Diocotron_table_grid}
        Diocotron instability simulation. 
        Number of cells and domains of the patches and equivalent global domains for the 
        advection tests. 
        In parentheses, we indicate the type of splitting. If the multi-patch splitting is 1D, 
        we only need to compute the derivatives along $r$. If the splitting is 2D, we compute the 
        derivatives along $r$, along $\theta$ and the cross-derivatives. 
        }
\end{table}

\paragraph{Test 3.1 Uniform global mesh}
\label{diocotron_uniform_section}
\paragraph{}
This first simulation is run on the same mesh as in \cite{Zoni2019}. 
The mesh with $128\times 256$ cells is split into three patches. 
The number of total cells and their shape are the same as in the article \cite{Zoni2019}.

\begin{figure}[H]
    \centering
    \includegraphics[width=\textwidth]{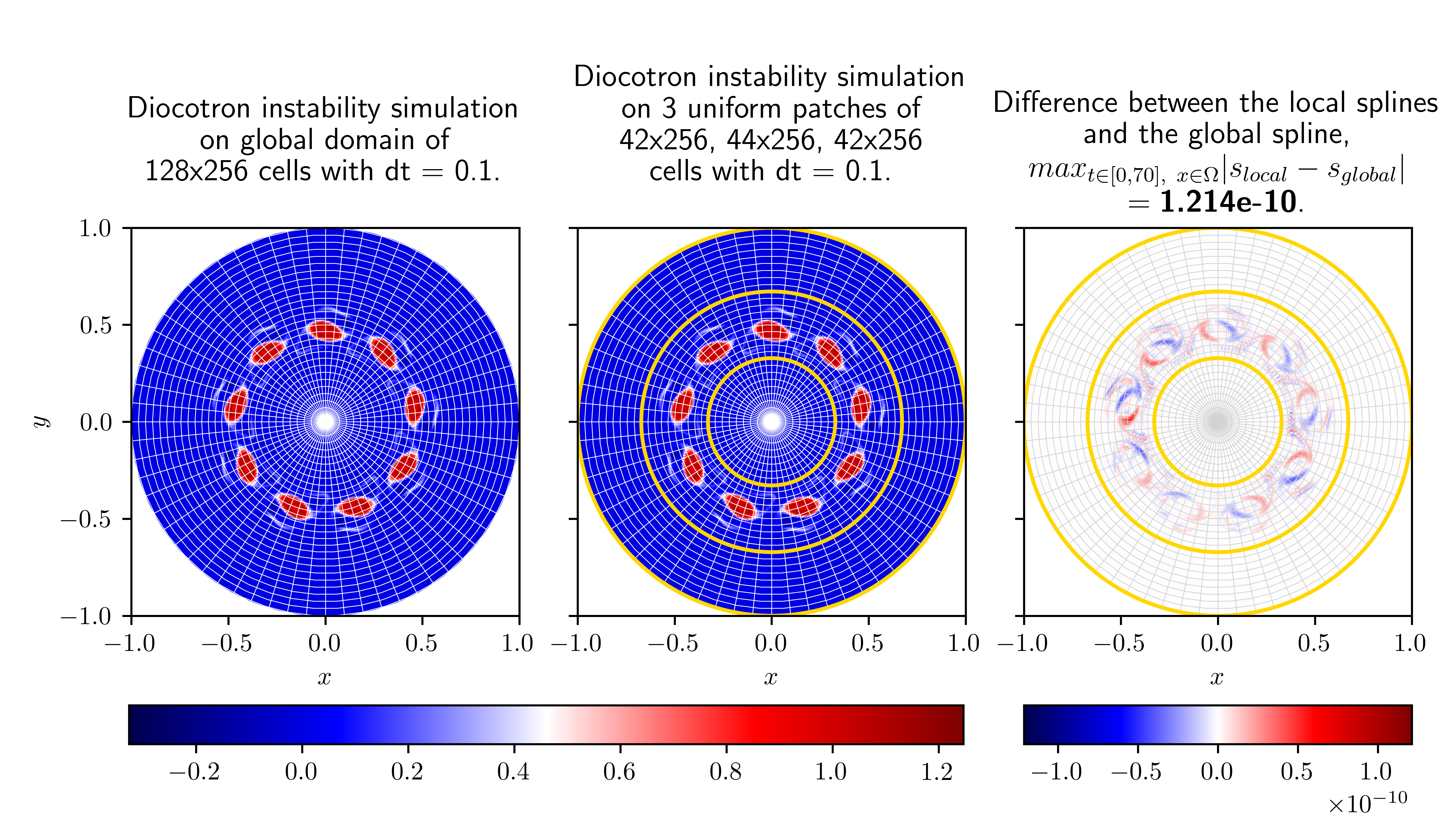}

    \vspace*{-0.25cm}
    
    \caption{\label{diocotron_uniform_fig}
        Diocotron instability simulation.
        On the left, results with advection on a global spline at $t = 70.0$. 
        This result agrees with the result in article \cite{Zoni2019}.
        In the middle, results with advection on local splines at $t = 70.0$. 
        The three yellow circles represent the interfaces between the patches. 
        On the right, the difference $(\rho_{\text{global}} - \rho_{\text{local}})$ at
        the same time. 
    }
\end{figure}

The results are presented in Fig. \ref{diocotron_uniform_fig}. 
We expect an error of $10^{-14}$ between the two curves at each interpolation 
(see Table \ref{Interpolation_tables}). 
The time step is $\Delta t = 0.1$ and the final time, $T = 70.0$. There are $700$ iterations. 
In the predictor-corrector method applied to solve the guiding-center equations, the density
is interpolated twice and the advection field once. So, obtaining an error of $10^{-10}$ in Fig.
\ref{diocotron_uniform_fig} is consistent with the estimations.  The results are not as good as for the 
simple advection. This is due to the errors on the advection field, which is also interpolated 
using the same method. 

\begin{figure}[H]
    \centering
    \includegraphics[width=0.49\textwidth]{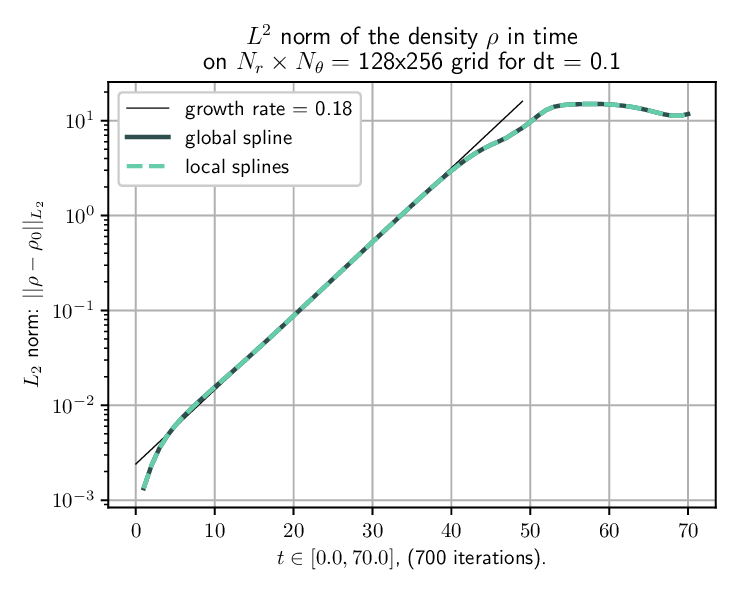}
    \includegraphics[width=0.49\textwidth]{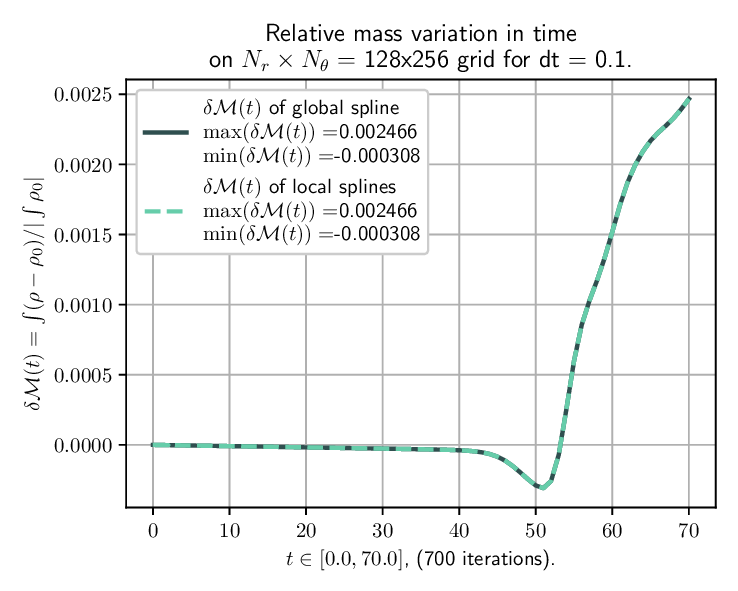}

    \vspace*{-0.25cm}

    \caption{\label{diocotron_norms_mass_fig}
        Diocotron instability simulation.
        Verification of the results: $L^2$-norm of the perturbed density and relative errors 
        on the total mass.
        As expected, the growth rate of the $L^2$-norm follows the analytical one given in 
        Fig. 6 of \cite{Zoni2019}. The relative errors on the mass are small. 
    }
\end{figure}

Fig. \ref{diocotron_norms_mass_fig} confirms that the $L^2$-norm of the perturbation in the 
multi-patch case gives rise to the expected slope of $0.18$ (for $m = 9$) as calculated in \cite{Zoni2019}. 
It does not affect the mass variation. The variation stays below $0.25\%$ of the initial mass. 

\paragraph{Test 3.2. 2D multi-patch case with T-joints and local refined}
\label{diocotron_Tjoints_refined_section}
We now test the simulation on a domain with a 2D multi-patch case. The patches also form T-joints
with the interfaces. The simulation uses the same layout as in Fig. \ref{Tjoint_study_case}, 
in paragraph \ref{Example_T_joints_geometry}.

One advantage of the multi-patch approach is that we can deal with local refinement
while working on a uniform B-spline basis, which can be more optimized than the non-uniform 
ones. Here, the most interesting patch is Patch 2. The density is quasi-null in the 
other patches. So, we refine Patch 2 and remove some cells in the patches 1 and 3. 
The global domain is no longer uniform but still conforming.

We depopulated the central patch 1 to the minimum number of cells necessary to define a cubic spline 
(i.e. four cells). 

In addition, we tested shifting one interface to the middle of the 
solution so as not to be influenced by the fact that the density was quasi-null at the interfaces 
in the previous simulations.
The interface between patch 2 and patch 3 has been shifted to the middle of the
ring of the initial density profile. 

\begin{figure}[H]
    \centering
    \includegraphics[width=\textwidth]{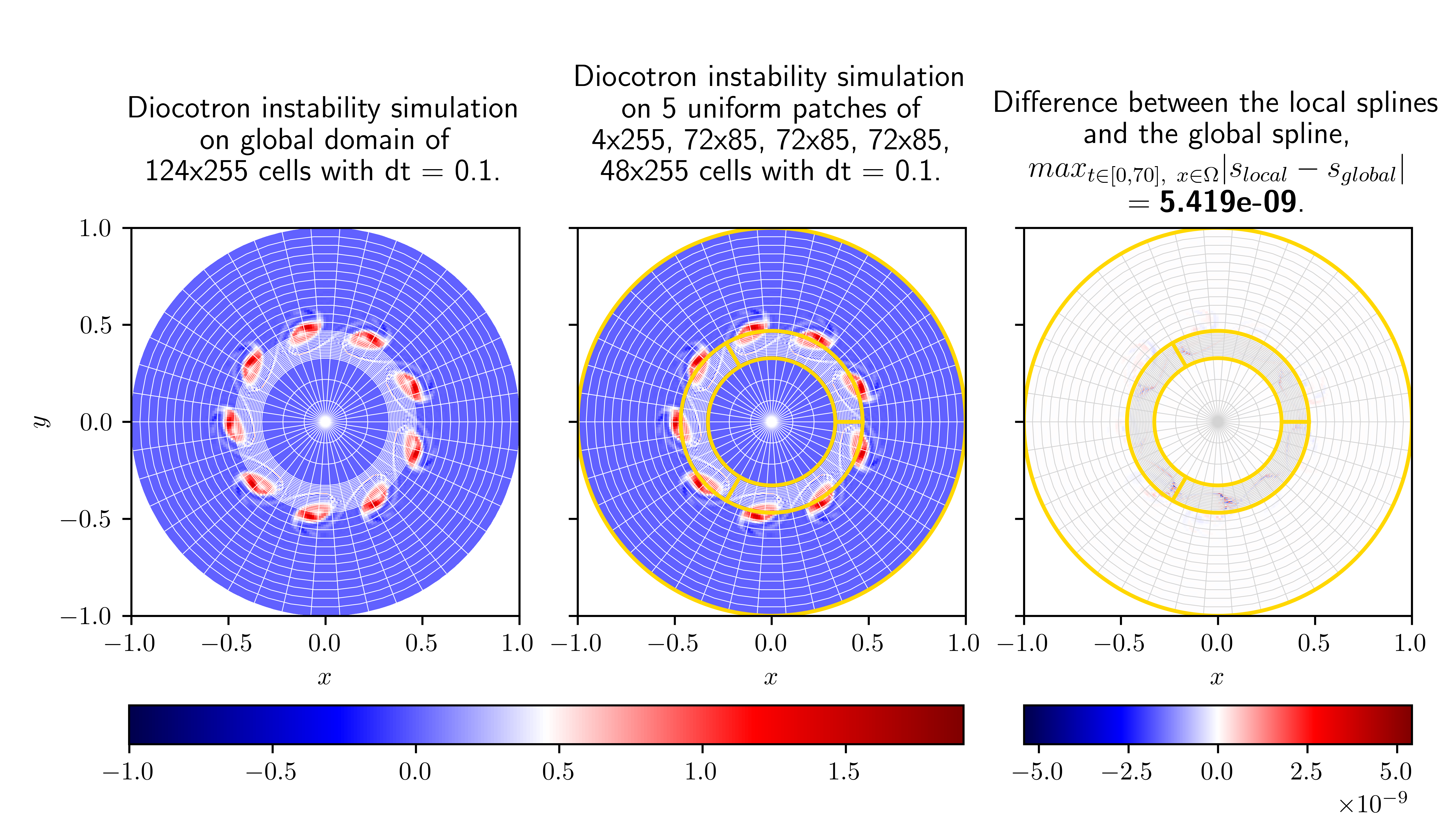}

    \vspace*{-0.25cm}

    \caption{\label{diocotron_Tjoint_refined_fig}
        Diocotron instability simulation with a 2D multi-patch case and T-joints. 
        We observe results similar to the previous simulations despite the use of a
        non-optimal mesh for this simulation (ill-placed interface and depopulation 
        of the central patch).
        The local splines and the equivalent global spline still agree.  
    }
\end{figure}

The patches are defined in Table \ref{Diocotron_table_grid}. 
The results are presented in Fig. \ref{diocotron_Tjoint_refined_fig}.
The error here is higher than in the previous simulation. This is due to the less adapted mesh. 
The local spline on patch 1 is not precise due to the depopulation. 
The location of the interfaces with patch 5 is not optimized.
 
\begin{figure}[H]
    \centering
    \includegraphics[width=0.49\textwidth]{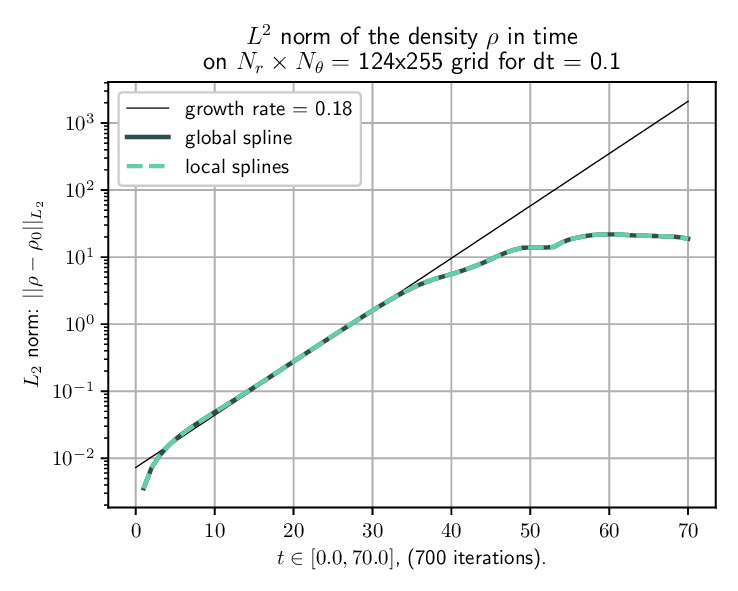}
    \includegraphics[width=0.49\textwidth]{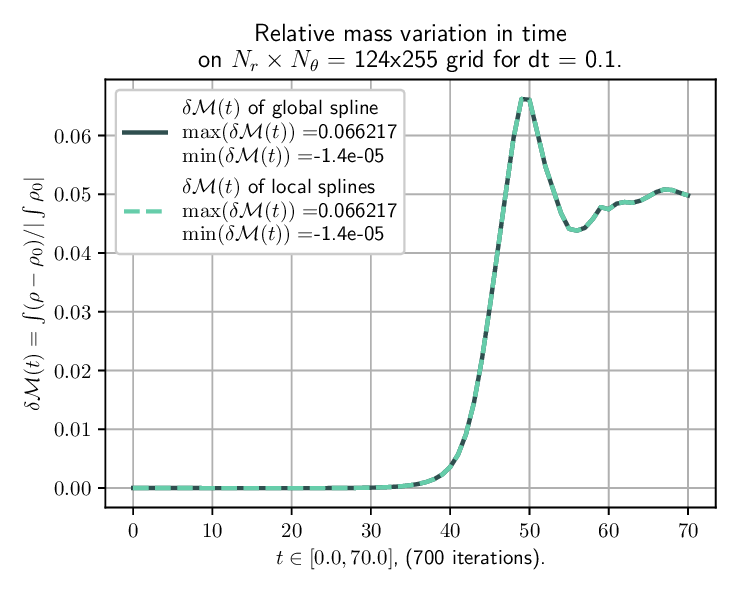}

    \vspace*{-0.25cm}

    \caption{\label{diocotron_norms_mass_Tjoint_refined_fig}
        Diocotron instability simulation with a 2D multi-patch case and T-joints. 
        Verification of the results: $L^2$-norm of the perturbed density and relative errors 
        on the total mass. We check that the suggested simulation on T-joints with refinement 
        is an admissible simulation. 
        The growth rate of the $L^2$-norm follows as expected the analytical one given in 
        \cite{Zoni2019}, Fig. 6. The relative errors on the mass stay small. 
    }
\end{figure}

Fig. \ref{diocotron_norms_mass_Tjoint_refined_fig} shows the admissibility of this simulation
despite the depopulation in patch 1 and the ill-placed interfaces at $r = \frac{60}{128}$.
The analytically predicted slope of the $L^2$-norm of the perturbation (\cite{Zoni2019}, Fig. 6) is retrieved.
The mass variation reaches $7\%$ of the initial mass. This choice of mesh is not recommended 
for this simulation. However, despite this mesh being poorly adapted, the local splines follow the 
equivalent global spline, and the results show a good agreement between them. 

\section{Conclusions and perspectives}
\paragraph{}
We presented a method to solve a classical advection equation on a multi-patch domain with 
the backward semi-Lagrangian method. In the backward semi-Lagrangian method, we mainly focus
on interpolating the function at the previous step.
The interpolation is made on cubic B-splines. 
In a multi-patch geometry, the conditions used to enforce patch connections are important to 
avoid instabilities.
In this paper, we showed that a $\mathcal{C}^0$ regularity at the patch interfaces may be 
not sufficient. 
Instead, we impose a $\mathcal{C}^1$ regularity. We then work with Hermite boundary conditions
and determine the derivatives at interfaces every time we want to interpolate a function. 
Starting from a formula developed in \cite{Crouseilles2009}, 
we demonstrate here a new generalization to compute the derivatives. 
This allows us to work with 
a variety of meshes (uniform and non-uniform with any cell length, conforming and some non-conforming). 
The generalized formula is exact but requires to solve a small global linear system. 
However, as the dependency on points far away from the interface decreases with distance, 
more local formulas can be derived at any precision.
In the uniform per patch case, we also give an explicit formulation
and study its asymptotic behavior when the number of cells per patch increases. 
The ideal case may be an approximation using thirty cells per patch to reconcile locality of the
computations and machine precision. 
Section \ref{section_Semi-Lagrangian_2D_multi-patch} also describes the application 
on 2D multi-patch domains, especially the non-conforming and T-joints cases. 
The numerical results of Section \ref{section_Numerical_results} show good agreement between the 
local splines defined on the patches and an equivalent global spline (if it exists). 

In the future, we will consider developing new tools and methods to deal with X-point geometries
as described in the Introduction \ref{section_introduction}. The X-point raises new difficulties
that cannot be handled using the same methods as for the O-point. 
The study in this article with T-joints should allow us to connect O-point geometries 
to X-point geometries in a next step. 
So, once interpolations on X-point geometry are handled, we will be able to advect on a 
complete tokamak-like poloidal cross-section as represented in Fig. \ref{JET_SOLEDGE}. 
 
\section{Declarations}
\subsection{Funding}
\paragraph{}
The project has received funding from the EoCoE-III project
(Energy-oriented Centre of Excellence for Exascale HPC applications).
This work has been carried out within the framework of the EUROfusion Consortium, 
funded by the European Union via the Euratom Research and Training Program 
(Grant Agreement No. 101052200 — EUROfusion). 
Emily Bourne's salary was paid for by the EUROfusion Advanced Computing Hub (ACH).

\subsection{Employment}
\paragraph{}
Emily Bourne is employed by École polytechnique fédérale de Lausanne (EPFL) in Switzerland. 
Virginie Grandgirard is employed by Commissariat à l’énergie atomique et aux énergies alternatives 
(CEA) in France. 
Michel Mehrenberger is employed by Institut de Mathématiques de Marseille (I2M) 
and Aix-Marseille University (AMU) in France. 
Eric Sonnendrücker is employed by Max-Planck-Institut für Plasmaphysik (IPP)
and Technische Universität München (TUM) in Germany. 
Pauline Vidal is employed by Max-Planck-Institut für Plasmaphysik (IPP) in Germany. 

\subsection{Competing interests}
\paragraph{}
The authors have no competing interests to declare that are relevant to the content of this article.

\subsection{Data Availability Statements}
\paragraph{}
This work was carried out using the Gyselalib++ library.
The library is open source and available on GitHub \cite{Gyselalibxx}.

\bibliographystyle{abbrv} 
\bibliography{biblio}

\appendix

\begin{appendices}
\renewcommand{\thesection}{Appendix \Alph{section}.}

\section[]{Knot sequences to define B-splines} 
\addcontentsline{toc}{section}{Appendix A. Knot sequences to define B-splines}
\label{Appendix_knots_and_Bsplines}
\normalsize
\paragraph{}
For $N_c$ cells, the 1D spline $s$ of degree $d = 3$ is given by 	
\begin{equation}
    s(x) = \sum_{k=0}^{N_c+d-1} c_k b_k(x), \qquad x\in\Omega
    \label{spline_representation}
\end{equation}
with $\{b_k\}_k$ the B-splines, and $\{c_k\}_k$ coefficients to determine. 
There are $N_c+1$ break points $\{x_{b,i}\}_{i=0}^{N_c}$.
We add $N_c+2d$ knots $\{k_i\}_{i=0}^{N_c+2d-1}$ to build the B-splines. 
The knot sequence is generally given as follows (especially in the uniform case),
\begin{equation*}
    \left\{
    \begin{aligned}
        & k_i = x_{b,0} - \frac{x_{b,N_c}-x_{b,0}}{N_c}(d-i),     
            & \qquad & 0 \leq i < d,      \\
        & k_i = x_{b,i-d},                                    
            & \qquad & d \leq i < N_c+d,    \\
        & k_i = x_{b,N_c} + \frac{x_{b,N_c}-x_{b,0}}{N_c}(i-N_c+d+1), 
            & \qquad & N_c+d \leq i < N_c+2d. \\
    \end{aligned}
    \right.
\end{equation*}

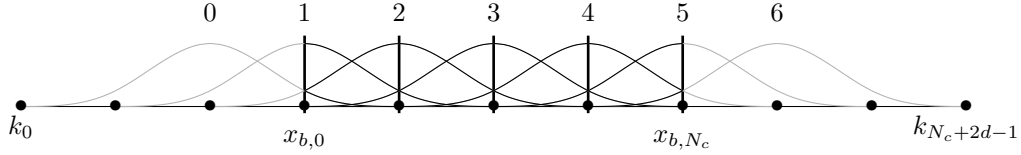
\begin{figure}[h]
    \centering
    \begin{tikzpicture}[xscale=1.25, yscale=1.25]
        \draw (0,0) -- (4,0); 
        \foreach \x in {0, 1, 2, 3, 4}{
            \draw[line width =1] (\x,-2pt) -- (\x,0.75); 
        }
        \node[below] at (0,-0.15) {$x_{b,0}$}; 
        \node[below] at (4,-0.15) {$x_{b,N_c}$}; 

        \draw[black] (-3,0) -- (7,0); 
        \foreach \xshft in {0, 1, 2, 3, 4, 5, 6}{
            \begin{scope}[xshift=(\xshft-3)*1cm]
                \node at (2, 1){\xshft}; 

                \begin{scope}[xshift=3cm]
                    \draw[black!30, domain=0:1, smooth, variable=\x, samples=10]
                    plot ({\x}, {(1-\x)^3 /6});
                \end{scope}

                \begin{scope}[xshift=2cm]
                    \draw[black!30, domain=0:1, smooth, variable=\x, samples=10]
                        plot ({\x}, {(1-\x)^2*(\x+2)/6 
                                    + (2-\x)*(1-\x)*(\x+1) /6
                                    + (2-\x)^2*(\x-0) /6});
                \end{scope}

                \begin{scope}[xshift=1cm]
                    \draw[black!30, domain=0:1, smooth, variable=\x, samples=10]
                        plot ({\x}, {(1-\x)*(\x+1)^2 /6
                                    + (2-\x)*(\x+1)*(\x-0) /6
                                    + (3-\x)*(\x-0)^2 /6});
                \end{scope}

                \draw[black!30, domain=0:1, smooth, variable=\x, samples=10]
                    plot ({\x}, {(\x-0)^3 /6});
            \end{scope}
        }

        \foreach \xshft in {0, 1, 2, 3}{
            \begin{scope}[xshift=(\xshft-3)*1cm]
                \begin{scope}[xshift=3cm]
                    \draw[black, domain=0:1, smooth, variable=\x, samples=10]
                    plot ({\x}, {(1-\x)^3 /6});
                \end{scope}
            \end{scope}
        }

        \foreach \xshft in {1, 2, 3, 4}{
            \begin{scope}[xshift=(\xshft-3)*1cm]
                \begin{scope}[xshift=2cm]
                    \draw[black, domain=0:1, smooth, variable=\x, samples=10]
                        plot ({\x}, {(1-\x)^2*(\x+2)/6 
                                    + (2-\x)*(1-\x)*(\x+1) /6
                                    + (2-\x)^2*(\x-0) /6});
                \end{scope}
            \end{scope}
        }

        \foreach \xshft in {2, 3, 4, 5}{
            \begin{scope}[xshift=(\xshft-3)*1cm]
                \begin{scope}[xshift=1cm]
                    \draw[black, domain=0:1, smooth, variable=\x, samples=10]
                        plot ({\x}, {(1-\x)*(\x+1)^2 /6
                                    + (2-\x)*(\x+1)*(\x-0) /6
                                    + (3-\x)*(\x-0)^2 /6});
                \end{scope}
            \end{scope}
        }

        \foreach \xshft in {3, 4, 5, 6}{
            \begin{scope}[xshift=(\xshft-3)*1cm]
                \draw[black, domain=0:1, smooth, variable=\x, samples=10]
                    plot ({\x}, {(\x-0)^3 /6});
            \end{scope}
        }

        \foreach \x in {-3, -2, -1, 0, 1, 2, 3, 4, 5, 6, 7}{
            \node at (\x,0) {$\bullet$}; 
        }
        \node[below] at (-3,0) { $k_0$}; 
        \node[below] at (7,0) { $k_{N_c+2d-1}$}; 

    \end{tikzpicture}
    \caption{\label{knots_seq}
        Cubic B-splines defined on $N_c = 4$ uniform cells with a knot 
        sequence with uniform outside knots ($\bullet$). 
        }
\end{figure}

The Fig. \ref{knots_seq} shows B-splines we obtain with this type of knot sequence. 

In the $2\pi$-periodic case, they can be given by 
\begin{equation*}
    \left\{
    \begin{aligned}
        & k_i = x_{b,N_c-1-i} - 2\pi, & \qquad & 0 \leq i < d,      \\
        & k_i = x_{b,i-d},          & \qquad & d \leq i < N_c+d,    \\
        & k_i = x_{b,i-N_c+d} + 2\pi, & \qquad & N_c+d \leq i < N_c+2d, \\
    \end{aligned}
    \right.
\end{equation*}
we obtained B-splines basis is similar to the one shown in Fig. \ref{knots_seq}.

We call \textit{open knots sequence}, if they are defined with repetitions of the first and the 
last break points as follows
\begin{equation*}
    \left\{
    \begin{aligned}
        & k_i = x_{b,0},   & \qquad & 0 \leq i < d,      \\
        & k_i = x_{b,i-d}, & \qquad & d \leq i < N_c+d,    \\
        & k_i = x_{b,N_c},   & \qquad & N_c+d \leq i < N_c+2d. \\
    \end{aligned}
    \right.
\end{equation*}

\begin{figure}[h]
    \centering
    \begin{tikzpicture}[xscale=1.25, yscale=1.25]
        \draw (0,0) -- (4,0); 
        \foreach \x in {0, 1, 2, 3, 4}{
            \draw[line width =1] (\x,-2pt) -- (\x, 1); 
        }
        \node[below] at (0,-0.5) {$x_{b,0}$}; 
        \node[below] at (4,-0.5) {$x_{b,N_c}$}; 

        \foreach \xshft in {3}{
            \begin{scope}[xshift=(\xshft-3)*1cm]
                \node at (2, 1.25){{\xshft}}; 

                \begin{scope}[xshift=3cm]
                    \draw[black, domain=0:1, smooth, variable=\x, samples=10]
                    plot ({\x}, {(1-\x)^3 /6});
                \end{scope}

                \begin{scope}[xshift=2cm]
                    \draw[black, domain=0:1, smooth, variable=\x, samples=10]
                        plot ({\x}, {(1-\x)^2*(\x+2)/6 
                                    + (2-\x)*(1-\x)*(\x+1) /6
                                    + (2-\x)^2*(\x-0) /6});
                \end{scope}

                \begin{scope}[xshift=1cm]
                    \draw[black, domain=0:1, smooth, variable=\x, samples=10]
                        plot ({\x}, {(1-\x)*(\x+1)^2 /6
                                    + (2-\x)*(\x+1)*(\x-0) /6
                                    + (3-\x)*(\x-0)^2 /6});
                \end{scope}

                \draw[black, domain=0:1, smooth, variable=\x, samples=10]
                    plot ({\x}, {(\x-0)^3 /6});
            \end{scope}
        }

        \draw[black, domain=0:1, smooth, variable=\x, samples=10]
            plot ({\x}, {(1 - \x)^3});

        \begin{scope}[xshift=4cm]
            \draw[black, domain=0:1, smooth, variable=\x, samples=10]
                plot ({-\x}, {(1 - \x)^3});
        \end{scope}

        \def\alpha{0.3333}; 
        \def\beta{-1.3333};
        \draw[black, domain=0:1, smooth, variable=\x, samples=10]
            plot ({\x}, {(-2*\alpha+1+\beta)*\x^3
                        +(-\beta - 2 +3*\alpha) * \x^2
                        + \x});
        \begin{scope}[xshift=1cm]
            \draw[black, domain=0:1, smooth, variable=\x, samples=10]
                plot ({\x}, {(\beta + 2*\alpha)*\x^3
                            +(-2*\beta - 3*\alpha)*\x^2 
                            + \beta*\x
                            + \alpha});
        \end{scope}	

        \begin{scope}[xshift=4cm]
            \draw[black, domain=0:1, smooth, variable=\x, samples=10]
                plot ({-\x}, {(-2*\alpha+1+\beta)*\x^3
                            +(-\beta - 2 +3*\alpha) * \x^2
                            + \x});
        \end{scope}	
        \begin{scope}[xshift=3cm]
            \draw[black, domain=0:1, smooth, variable=\x, samples=10]
                plot ({-\x}, {(\beta + 2*\alpha)*\x^3
                            +(-2*\beta - 3*\alpha)*\x^2 
                            + \beta*\x
                            + \alpha});
        \end{scope}

        \def\d{0.5}; 
        \def\c{0.83333333333};
        \draw[black, domain=0:1, smooth, variable=\x, samples=10]
            plot ({\x}, {\x^2*(\x+2) /6});   
        \begin{scope}[xshift=1cm]
            \draw[black, domain=0:1, smooth, variable=\x, samples=10]
                plot ({\x}, {((-5/6+\c+2*\d)*\x^3 
                            + (1-2*\c-3*\d)*\x^2
                            + \c*\x + \d)});
        \end{scope}	  
        \begin{scope}[xshift=2cm]
            \draw[black, domain=0:1, smooth, variable=\x, samples=10]
            plot ({\x}, {(1-\x)^3 /6});
        \end{scope}	

        \begin{scope}[xshift=4cm]
        \draw[black, domain=0:1, smooth, variable=\x, samples=10]
            plot ({-\x}, {\x^2*(\x+2) /6});  
        \end{scope}	 
        \begin{scope}[xshift=3cm]
            \draw[black, domain=0:1, smooth, variable=\x, samples=10]
                plot ({-\x}, {((-5/6+\c+2*\d)*\x^3 
                            + (1-2*\c-3*\d)*\x^2
                            + \c*\x + \d)});
        \end{scope}	  
        \begin{scope}[xshift=2cm]
            \draw[black, domain=0:1, smooth, variable=\x, samples=10]
            plot ({-\x}, {(1-\x)^3 /6});
        \end{scope}

        \node at (0.00, 1.25){0}; 
        \node at (0.66, 1.25){1}; 
        \node at (1.33, 1.25){2}; 
        \node at (2.66, 1.25){4}; 
        \node at (3.33, 1.25){5}; 
        \node at (4.00, 1.25){6}; 

        \foreach \x in {0, 1, 2, 3, 4}{
            \node at (\x,0) {$\bullet$}; 
        }
        \node[below left ] at (0,-0.1) {$k_0, k_1, k_2$}; 
        \node[below right] at (4,-0.1) {$k_{N_c+d-2}, k_{N_c+d-1}, k_{N_c+d}$}; 

    \end{tikzpicture}
    \caption{\label{open_knots_seq}
    Cubic B-splines defined on $N_c = 4$ uniform cells with an open knot sequence ($\bullet$).}
\end{figure}
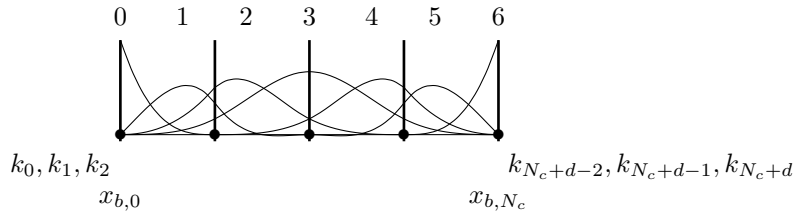

The Fig. \ref{open_knots_seq} shows B-splines we obtain with an open knot sequence. 

\section[]{Proof hints of the explicit formulation of the coefficients} 
\addcontentsline{toc}{section}{Appendix B. Proof hints of the explicit formulation of the coefficients}
\label{Appendix_proof_explicit_formulation}
\normalsize
\paragraph{}
We detail here some steps to demonstrate the equality between the recursive formulations 
\eqref{coeff_forward_n}, \eqref{coeff_backward_m} 
and the explicit formulations \eqref{explicit_relation}, \eqref{explicit_weights} 
of the coefficients in the uniform per patch case. The index $i$ corresponds to the 
index of the interpolation point $x_i$ located on the interface $\mathcal{X}_I$. 

In the uniform per patch case, the coefficients are recursively given by, 
for $n = 2$ and $m = 1$,
\begin{equation*}
    c^i_{1,2} 
        = \frac{1}{4 + a^i_{1,1}} 
            \left[4 c^i_{1,1} +  3 \frac{f^R_{i+2} - f^R_{i}}{\Delta x^R} a^i_{1,1} \right], 
    \qquad
    a^i_{1,2} = \frac{- a^i_{1,1}}{4 + a^i_{1,1}},
    \qquad
    b^i_{1,2} = \frac{4 b^i_{1,1}}{4 + a^i_{1,1}}, 
\end{equation*}
for $n \geq 2$ and $m = 1$,
\begin{equation*}
    \left\{
    \begin{aligned}
        & c^i_{1,n+1} = \frac{1}{4 + \frac{a^i_{1,n}}{a^i_{1,n-1}} } 
            \left[  
                4 c^i_{1,n} 
                + \frac{a^i_{1,n}}{a^i_{1,n-1}} c^i_{1,n-1}
                + 3 \frac{f^R_{i+n+1} - f^R_{i+n-1}}{\Delta x^R} a^i_{1,n}
            \right] \\
        & a^i_{1,n+1} = \frac{- a^i_{1,n}}{4 + \frac{a^i_{1,n}}{a^i_{1,n-1}} }, \\
        & b^i_{1,n+1} = \frac{1}{4 + \frac{a^i_{1,n}}{a^i_{1,n-1}} } 
            \left[
                4 b^i_{1,n}
                + \frac{a^i_{1,n} }{a^i_{1,n-1}} b^i_{1,n-1}
            \right],  \\
    \end{aligned}
    \right.
\end{equation*}
for $n = N^R_I$ and  $m = 2$,
\begin{equation*}
    c^i_{2,N^R} = \frac{1}{4 + b^{i}_{1,N^R}}
        \left[
            4 c^{i}_{1,N^R} +  3 \frac{f^L_{i} - f^L_{i-2}}{\Delta x^L} b^{i}_{1,N^R}
        \right], 
    \qquad
    a^i_{2,N^R} = \frac{4 a^{i}_{1,N^R}}{4 + b^{i}_{1,N^R}}, 
    \qquad
    b^i_{2,N^R} = \frac{- b^{i}_{1,N^R}}{4 + b^{i}_{1,N^R}},
\end{equation*}
for $n = N^R_I$ and  $m \geq 2$,
\begin{equation*}
    \left\{
    \begin{aligned}
        & c^{i}_{m+1,N^R} = \frac{1}{4 + \frac{b^{i}_{m,N^R}}{b^{i}_{m-1,N^R}}}
            \left[
                4 c^{i}_{m,N^R}
                + \frac{b^{i}_{m,N^R}}{b^{i}_{m-1,N^R}} c^{i}_{m-1,N^R}
                + 3 \frac{f^L_{i-m+1} - f^L_{i-m-1}}{\Delta x^L} b^{i}_{m,N^R}  
            \right],
        \\
        & a^{i}_{m+1,N^R} = \frac{1}{4 + \frac{b^{i}_{m,N^R}}{b^{i}_{m-1,N^R}}}
            \left[
                4 a^{i}_{m,N^R} 
                + \frac{b^{i}_{m,N^R}}{b^{i}_{m-1,N^R}} a^{i}_{m-1,N^R} 
            \right].
        \\
        & b^{i}_{m+1,N^R} = \frac{- b^{i}_{m,N^R}}
            {4 + \frac{b^{i}_{m,N^R}}{b^{i}_{m-1,N^R}}}.
        \\
    \end{aligned}
    \right.
\end{equation*}

Looking at the expressions of the $\{a^i_{1,n}\}_{n = 2, ..., N^R_I}$
and $\{b^{i}_{m,N^R}\}_{m = 2, ..., N^L_I}$ sequences,
we start by introducing the $w$ sequence such that $\forall k > 1$, 
\begin{equation}
    w_k = \frac{-1}{4 + w_{k-1}}.
\label{w_sequence_relation}
\end{equation}

The function $x \mapsto \frac{-1}{4+x}$ has two fixed points: $-(2-\sqrt{3})$ and $-(2+\sqrt{3})$. 
Successively injecting the relation \eqref{w_sequence_relation} into itself  gives us 
\begin{equation*}
    w_k = \frac{-(4 + w_{k-2})}{4\times4 -1 + 4 w_{k-2}}, 
    \qquad
    w_k = \frac{-(4\times4-1 + 4 w_{k-3})}{4(4\times4 -1)- 4 + (4\times4-1) w_{k-3}}.
\end{equation*}
We can see appearing the sequence of coefficients in the fraction: 
$\{0, 1, 4, 4\times4 - 1, 4(4\times4 -1)- 4 , ...\}$.
It corresponds to a sequence $u$ such that $u_{k+1} = 4 u_{k} - u_{k-1}$, $\forall k\in\mathbb{Z}$
with $u_0 = 0$  and $u_1 = 1$. It is a linear-recursive sequence whose the characteristic
polynomial $x^2 + 4x -1$ has two zeros: $(2-\sqrt{3})$ and $(2+\sqrt{3})$. 
The sequence $u$ is then in the vector space 
$\spn(\{\{(2-\sqrt{3})^k \}_{k\in\mathbb{N}},\{(2+\sqrt{3})^k \}_{k\in\mathbb{N}}\})$ 
and can be expressed as follows $\forall k > 1$, 
\begin{equation*}
    u_k = (2+\sqrt{3})^k - (2-\sqrt{3})^k,
\end{equation*}
and we can prove that $w$ can be expressed from $u$ as following
\begin{equation*}
    w_k = - \frac{u_{k-1} + u_{k-2}w_1}{u_{k} + u_{k-1}w_1}.
\end{equation*}
Indeed, $\left\{- \frac{u_{k-1} + u_{k-2}w_1}{u_{k} + u_{k-1}w_1}\right\}_{k\geq 1}$
verifies the relation \eqref{w_sequence_relation}.
For $k = 1$, $\left\{- \frac{u_{k-1} + u_{k-2}w_1}{u_{k} + u_{k-1}w_1}\right\}_{k\geq 1}$
is equal to the same initial condition $w_1$ as $w$. 

It is easy to see that by fixing $w_1 = a^i_{1,1}$, 
the $\left\{w_n\right\}_{n = 2, ..., N^R_I}$  sequence gives the same values as 
the $\left\{\frac{a^i_{1,n}}{a^i_{1,n-1}}\right\}_{n = 2, ..., N^R_I}$ sequence.
By fixing $w_1 = b^i_{1, N^R_I}$, 
the $\left\{w_m\right\}_{m = 2, ..., N^L_I}$  sequence  gives the same values as 
the $\left\{\frac{b^i_{m,N^R_I}}{b^i_{m-1,N^R_I}}\right\}_{m = 2, ..., N^L_I}$ sequence.
 
By factoring the two equalities, we deduce explicit formulas for
$\left\{a^i_{1,n}\right\}_{n = 2, ..., N^R_I}$
and $\left\{b^i_{m,N^R_I}\right\}_{m = 2, ..., N^L_I}$, 
\begin{equation*}
    a^i_{1,n}  = \frac{(-1)^{n-1}u_1 a^i_{1,1}}{u_n + u_{n-1}a^i_{1,1}}, 
    \qquad
    b^i_{m,N^R_I} = \frac{(-1)^{m-1}u_1 b^i_{1,N^R_I}}{u_m + u_{m-1} b^i_{1,N^R_I}}.
\end{equation*}

Manipulating the first terms of the 
$\left\{b^i_{1,n}\right\}_{n = 1, ..., N^R_I}$ 
and $\left\{a^i_{m,N^R_I}\right\}_{m = 1, ..., N^L_I}$ sequences
helps us to identify the following expressions, 
\begin{equation*}
    b^i_{1,n}  = \frac{u_n b^i_{1,1}}{u_n + u_{n-1}a^i_{1,1}}, 
    \qquad
    a^i_{m,N^R_I}  = \frac{u_m a^i_{1,N^R_I}}{u_m + u_{m-1} b^i_{1,N^R_I}},
\end{equation*}
which can be shown by recursion. 

We then inject the expressions for $a^i_{1,N^R_I}$ and $b^i_{1,N^R_I}$
into $a^i_{m,n}$ and $b^i_{m,n}$ for $n = N^R_I$,
\begin{equation*}
    a^i_{m,n} = \frac{(-1)^{n-1}u_{1}a^i_{1,1}u_{m}}
                     {u_{m}u_{n} + u_{m}u_{n-1}a^i_{1,1} + u_{n} u_{m-1}b^i_{1,1}}, 
    \qquad
    b^i_{m,n} = \frac{(-1)^{m-1} u_{1} b^i_{1,1} u_{n}}
                     {u_{m}u_{n} + u_{m}u_{n-1}a^i_{1,1} + u_{n} u_{m-1}b^i_{1,1}}.
\end{equation*}

The expressions for the weights $\{\omega^i_{k,N^L_I,N^R_I}\}_{k = -N^L_I, ..., N^R_I}$ in 
the coefficient $c^i_{N^L_I,N^R_I}$ are also demonstrated by recursion.
We start by the first recursion which can be written as follows thanks to the results on 
$a^i_{m,n}$ and $b^i_{m,n}$ for $n \in \mathbb{N}^*$, $m = 1$,
\begin{equation*}
    c^i_{1,n+1} = 4\frac{u_{n} + u_{n-1}a_{1,1}}{u_{n+1} + u_{n}a_{1,1}} c^i_{1,n} 
            - \frac{u_{n-1} + u_{n-2}a_{1,1}}{u_{n+1} + u_{n}a_{1,1}} c^i_{1,n-1}
            - \frac{3(-1)^{n-1}u_{1}a_{1,1}}{u_{n+1} + u_{n}a_{1,1}} \frac{f^R_{i-1+n}}{\Delta x^R}  
            + \frac{3(-1)^{n-1}u_{1}a_{1,1}}{u_{n+1} + u_{n}a_{1,1}} \frac{f^R_{i+1+n}}{\Delta x^R}.
\end{equation*}

We identify six cases where $k = n+1$, $k = n$, $k = n-1$, $k = 1, ..., n-2$, 
$k = 0$ and $k = -1$. Manipulating the first terms 
helps us to identify the following expressions, 
\begin{equation*}
\left\{
\begin{aligned}
    & \omega_{k,1,n} 
    = \frac{3(-1)^{k}}{u_{n} + u_{n-1}a_{1,1}} \frac{a_{1,1}}{\Delta x^R} u_{1}
    ,\qquad & k = n 
    \\
    & \omega_{k,1,n} 
    = \frac{3(-1)^{k}}{u_{n} + u_{n-1}a_{1,1}} \frac{a_{1,1}}{\Delta x^R} (u_{n-k+1} - u_{n-k-1})
    ,\qquad & k = 1, ..., n-1 
    \\
    & \omega_{k,1,n} 
    = \frac{3(-1)^{k}}{u_{n} + u_{n-1}a_{1,1}}
        \left(\frac{a_{1,1}}{\Delta x^R}(u_{n} - u_{n-1}) -\frac{b_{1,1}}{\Delta x^L} u_n\right) 
    ,\qquad & k = 0
    \\
    & \omega_{k,1,n} 
    = \frac{3(-1)^{k+1}}{u_{n} + u_{n-1}a_{1,1}} \frac{b_{1,1}}{\Delta x^L} u_{n}
    ,\qquad & k = -1
    \\
\end{aligned}
\right.
\end{equation*}
which can be shown by recursion.
It is not detailed here, there  are no particular difficulties.

We inject these results in the second recursion. Rewritten, it gives 
\begin{equation*}
\left\{
\begin{aligned}
    & c^i_{2,n} 
    = \frac{4(u_n + u_{n-1}a^i_{1,1})}{u_n (4 + b^i_{1,1}) + u_{n-1}4 a^i_{1,1}}c^i_{n,1} 
    + \frac{3 u_n b^i_{1,1}}{u_n (4 + b^i_{1,1}) + u_{n-1}4 a^i_{1,1}}\frac{f^R_{i} - f^L_{i-2}}{\Delta x^L}
    \\
    & c^i_{m+1, n} = \frac{1}{u_{n}u_{m+1} + u_{n}u_{m} b^i_{1,1} + u_{n-1}u_{m+1}a^i_{1,1}}
    \left[
        4(u_{n}u_{m} + u_{n}u_{m-1} b^i_{1,1} + u_{n-1}u_{m}a^i_{1,1}) c^i_{m,n}
        \right.\\ &\hspace*{7cm}
        + 3 (-1)^{m-1} u_1 b^i_{1,1} u_n \frac{f^L_{i+1-m} - f^L_{i-1-m}}{\Delta x^L}
        \\ &\hspace*{7cm}\left.
        -(u_{n}u_{m-1} + u_{n}u_{m-2} b^i_{1,1} + u_{n-1}u_{m-1}a^i_{1,1}) c^i_{m-1,n}
    \right].
\end{aligned}
\right.
\end{equation*}

We also identify different cases where $k = n$, $k = 1, ..., n-1$, 
$k = 0$, $k = -(m-2), ..., -1$, $k = -(m-1)$, $k = -m$ and $k = -(m+1)$.
It leads us to the followings expressions, 
\begin{equation*}
\left\{
\begin{aligned}
    & \omega^i_{k,m,n} = 3(-1)^{k}
        \frac{\frac{a^i_{1,1}}{\Delta x^R} u_{m}u_{1}}
        {u_{m}u_{n} + u_{m}u_{n-1}a^i_{1,1} + u_{n} u_{m-1}b^i_{1,1}}
        , &k = n, \\
    & \omega^i_{k,m,n} = 3(-1)^{k}
        \frac{\frac{a^i_{1,1}}{\Delta x^R} u_{m}(u_{n-k+1} - u_{n-k-1})}
        {u_{m}u_{n} + u_{m}u_{n-1}a^i_{1,1} + u_{n} u_{m-1}b^i_{1,1}}
        , &k = 1, ..., n-1, \\
    & \omega^i_{k,m,n} = 3(-1)^k
        \frac{\frac{a^i_{1,1}}{\Delta x^R}u_{m} (u_{n} - u_{n-1}) 
            - \frac{b^i_{1,1}}{\Delta x^L}u_{n}(u_{m} - u_{m-1})}
        {u_{m}u_{n} + u_{m}u_{n-1}a_{1,1} + u_{n} u_{m-1}b^i_{1,1}}
        , &k = 0, \\
    & \omega^i_{k,m,n} = 3(-1)^{k+1}
        \frac{\frac{b^i_{1,1}}{\Delta x^L} u_{n}(u_{m+k+1} - u_{m+k-1})}
        {u_{m}u_{n} + u_{m}u_{n-1}a^i_{1,1} + u_{n} u_{m-1}b^i_{1,1}}
        , &k = -(m-1), ..., -1, \\
    & \omega^i_{k,m,n} = 3(-1)^{k+1}
        \frac{\frac{b^i_{1,1}}{\Delta x^L} u_{n}u_{1}}
        {u_{m}u_{n} + u_{m}u_{n-1}a^i_{1,1} + u_{n} u_{m-1}b^i_{1,1}}
        , &k = -m. \\
\end{aligned}
\right.
\end{equation*}
The demonstration is also made by recursion. 

\end{appendices}
\end{document}